\documentclass[12pt,english]{article}

\usepackage{hyperref}
\hypersetup{
    unicode      = false,
    pdftoolbar   = true,
    pdfmenubar   = true,
    pdffitwindow = true,
    pdfnewwindow = true,
    colorlinks   = true,
    linkcolor    = blue,
    citecolor    = blue,
    filecolor    = blue,
    urlcolor     = blue,
    breaklinks   = true
}

\usepackage[utf8]{inputenc}
\usepackage[T1]{fontenc}
\usepackage{amsthm,amssymb,amsmath}
\usepackage{cleveref}
\usepackage[english]{babel}
 \usepackage{lmodern}
\usepackage{graphicx}
\usepackage{fullpage}
\usepackage{caption}
\usepackage{version}
\usepackage{subcaption}

\usepackage[textwidth=1.2\marginparwidth,backgroundcolor=yellow,linecolor=orange,textsize=scriptsize]{todonotes}

\usepackage{orcidlink}

\usepackage{lipsum}
\usepackage{geometry}
\usepackage{systeme}
\usepackage{array,multirow,makecell}
\setcellgapes{1pt}
\makegapedcells
\newcolumntype{R}[1]{>{\raggedleft\arraybackslash }b{#1}}
\newcolumntype{L}[1]{>{\raggedright\arraybackslash }b{#1}}
\newcolumntype{C}[1]{>{\centering\arraybackslash }b{#1}}
\usepackage{graphicx}
\usepackage{tikz}
\usetikzlibrary{positioning,decorations.pathreplacing,shapes}
\usepackage{pict2e}
\usepackage{amsfonts}
\usepackage{multirow}
\usepackage[normalem]{ulem}
\usepackage{epsfig}
\usepackage{soul}
\usepackage{pifont}
\usepackage{colortbl}
\usepackage{enumitem}
\usepackage{optidef}
\usepackage{epic}
\usepackage{xspace}
\usepackage{pgf}
\usepackage{rotating}
\usepackage{array}
\usepackage{geometry}
\usepackage{times}
\usepackage{algorithm}
\usepackage{algpseudocode}
\usepackage[justification=centering]{caption}
\usepackage[flushleft]{threeparttable}

\usetikzlibrary{patterns}
\usetikzlibrary{arrows}

\makeatletter
\tikzset{
        hatch distance/.store in=\hatchdistance,
        hatch distance=5pt,
        hatch thickness/.store in=\hatchthickness,
        hatch thickness=5pt
        }
\newcounter{subsubsubsection}[subsubsection]
\renewcommand\thesubsubsubsection{\@roman\c@subsubsubsection}
\newcommand\subsubsubsection{\@startsection{subsubsubsection}{4}{\z@}%
                                     {-3.25ex\@plus -1ex \@minus -.2ex}%
                                     {1.5ex \@plus .2ex}%
                                     {\normalfont\small\bfseries}}
\newcommand*\l@subsubsubsection{\@dottedtocline{3}{5.2em}{1em}}
\newcommand*{\subsubsubsectionmark}[1]{}

\makeatother

\def\M{\mathbf{M}}
\def\g{\mathbf{g}}
\def\v{\mathbf{v}}
\def\U{\mathbf{U}}
\def\b{\mathbf{b}}

\def\p{\mathbf{p}}
\def\x{\mathbf{x}}
\def\u{\mathbf{u}}
\def\s{\mathbf{s}}
\def\t{\mathbf{t}}
\def\y{\mathbf{y}}
\def\X{\mathbf{X}}
\def\R{{\mathbb{R}}}
\def\N{{\mathbb{N}}}
\def\M{{\mathbf{M}}}
\def\V{{\mathbf{V}}}

\newcommand{\xmark}{\ding{55}}
\newcommand{\argmin}{\mathop{\mathrm{argmin}}}

\newcommand{\cvar}{\mathop{\mathrm{CVaR}}}

\newcommand{\var}{\mathop{\mathrm{VaR}}}

\newtheorem{Def}{Definition}[section]
\newtheorem{Lem}[Def]{Lemma}

\newtheorem{Th}[Def]{Theorem}
\newtheorem{Cor}[Def]{Corollary}
\newtheorem{assumption}{Assumption}

\title{Risk averse constrained blackbox optimization under mixed aleatory/epistemic uncertainties}

\author{
    \addtocounter{footnote}{1}
	\href{mailto:Charles.Audet@gerad.ca}{Charles Audet \orcidlink{0000-0002-3043-5393}}
	\thanks{
		{GERAD}
		and D\'epartement de math\'ematiques et g\'enie industriel,
		\'Ecole Polytechnique de Montr\'eal,
		C.P. 6079, Succ. Centre-ville,
		Montr\'eal, Qu\'ebec, Canada H3C~3A7. \newline
		Mail: charles.audet@gerad.ca 
	}
	\and
	\href{mailto:jean.bigeon@ls2n.fr}{Jean Bigeon \orcidlink{0000-0002-6112-6913}}
	\thanks{
	Nantes Université, École Centrale Nantes, CNRS, LS2N, UMR 6004, F-44000 Nantes, France. \newline
	Mail: jean.bigeon@ls2n.fr 
	}
	\and 
	\href{mailto:romain.couderc@grenoble-inp.fr} 
 {Romain Couderc \orcidlink{0000-0002-8696-6361}} \footnotemark[2] 
	\thanks{Univ. Grenoble Alpes, CNRS, Grenoble INP*, G-SCOP, 38000 Grenoble, France. \newline
	Mail: romain.couderc@grenoble-inp.fr \newline
    *Institute of Engineering Univ. Grenoble Alpes
    }
    \and
    \href{mailto:michael.kokkolaras@mcgill.ca}{Michael Kokkolaras \orcidlink{0000-0003-1546-3393}} 
    \thanks{Gerad and department of Mechanical Engineering, McGill University, Montreal, Canada  \newline
    Mail: michael.kokkolaras@mcgill.ca}    
}

\date{\today}

\begin{document}
\maketitle

\begin{abstract}
This paper addresses risk averse constrained optimization problems where the objective and constraint functions can only be computed by a blackbox subject to unknown uncertainties. To handle mixed aleatory/epistemic uncertainties, the problem is transformed into a conditional value-at-risk (CVaR) constrained optimization problem. General inequality constraints are managed through Lagrangian relaxation. A convolution between a truncated Gaussian density and the Lagrangian function is used to smooth the problem. A gradient estimator of the smooth Lagrangian function is derived, possessing attractive properties: it estimates the gradient with only two outputs of the blackbox, regardless of dimension, and evaluates the blackbox only within the bound constraints. This gradient estimator is then utilized in a multi-timescale stochastic approximation algorithm to solve the smooth problem. Under mild assumptions, this algorithm almost surely converges to a feasible point of the CVaR-constrained problem whose objective function value is arbitrarily close to that of a local solution. Finally, numerical experiments are conducted to serve three purposes. Firstly, they provide insights on how to set the hyperparameter values of the algorithm. Secondly, they demonstrate the effectiveness of the algorithm when a truncated Gaussian gradient estimator is used. Lastly, they show its ability to handle mixed aleatory/epistemic uncertainties in practical applications.
\end{abstract}


\newpage
\section{Introduction}

Blackbox optimization (BBO) is concerned with optimization problems where the functions used to compute the objective and the constraints are blackboxes. In optimization, a blackbox is any process that returns an output when an input is provided, but the inner workings of that process are not analytically available~\cite{AuHa2017}. This type of problem is common in signal processing~\cite{curtis2020adaptive}, machine learning~\cite{papernot2017practical}, and engineering design~\cite{alarie2021two, Kokko06}. In the presence of uncertainties, a constrained blackbox optimization problem may be formulated as follows
\begin{align}
\label{problem_ref}
\begin{split}
\begin{array}{cl}
       \displaystyle \min_{\x \in \mathcal{X} \subset \R^n}   &\Xi_0 [C_0( \x, \boldsymbol{\xi})] \\
     \mbox{s.t.} &\Xi_j [C_j(\x, \boldsymbol{\xi})] \leq 0, \; \forall j \in [1, m],
\end{array}
\end{split}
\end{align}
where $\x$ is the vector of the design variables, $ \mathcal{X} := [\b_\ell, \b_u]$ is a hyperectangle, and $\boldsymbol{\xi}$ is the vector modelling the uncertainties. The source of uncertainties may arise from the design variables, the parameters, the inner processes of the blackbox (for example, when Monte Carlo simulation is used in the blackbox), or even combinations of these factors. Uncertainties may or may not depend on $\mathbf{x}$. $C_0(\cdot, \boldsymbol{\xi})$ denotes the version of the objective function $c_0: \mathcal{X} \to \R$ subject to uncertainties, while for all $j \in \{1,2,\ldots, m\}, C_j(\cdot, \boldsymbol{\xi})$ denotes the version of the constraint $c_j : \mathcal{X} \to \R$ subject to the uncertainties (also called the limit state function in the reliability community). 
Since the objective function and the constraints depend on the uncertainty vector, the measures $\Xi_j, j \in  \{0,1,\ldots, m\}$ are used to map them into $\R$. It follows from this formulation that the key factor is the selection of the uncertainties model, which in turn determines the choice of the measures $\Xi_j$. In the following, various methods commonly found in the literature are presented, depending on the assumptions made, the chosen uncertainty model, and the level of information available about these uncertainties.

\subsection{Related work}
 In probabilistic reliability-based design optimization (RBDO), uncertainties are considered as random vectors with known probabilistic distributions. In this field~\cite{deb2009reliability}, Problem (\ref{problem_ref}) is transformed into the following 
 \begin{align}
\label{problem_rbdo}
\begin{split}
\begin{array}{cl}
    \displaystyle \min_{\x \in \mathcal{X} \subset \R^n}   &C_0(\x, \p ) \\
     \mbox{s.t.}  &\mathbb{P} [C_j(\x, \p) \leq 0] \geq \alpha_j, \; \forall j \in [1,m],
\end{array}
\end{split}
\end{align}
where $\alpha_j, j \in \{1, \ldots, m\}$ are the desired reliability levels and $\x$ and $\p$ are the means of the noised design variables and parameters respectively. In this reformulation, the expectation is utilized to handle the uncertainties in the objective function, and a linear approximation is employed to derive the deterministic objective function~\footnote{For a differentiable function $C_0$ perturbed only by uncertainties in its design variables. These uncertainties can be written as $\x + \boldsymbol{\xi}_\x$ where $\x = \mathbb{E}_{\boldsymbol{\xi}_\x}[\x + \boldsymbol{\xi}_\x]$. Then a first-order Taylor approximation of the function gives that $\mathbb{E}_{\boldsymbol{\xi}_\x}[C_0(\x + \boldsymbol{\xi}_\x)] \approx \mathbb{E}_{\boldsymbol{\xi}_\x}[C_0(\x) + \nabla C_0(\x)^T \boldsymbol{\xi}_\x] = C_0(\x)$. A similar observation holds for the parameters.}. To address the uncertainties in the constraints, a probability measure is employed. The conventional approach to solving Problem (\ref{problem_rbdo}) involves two nested loops: the outer loop searches for an optimal design, while the inner loop evaluates the feasible probability of the optimal candidate.  

The inner loop is often computationally demanding due to the time-consuming estimation of feasible probabilities. To address this challenge, numerically efficient methods for RBDO problems have been developed. In a first set of methods, the inner loop involves solving a deterministic optimization problem. The fundamental idea behind this class of methods is to identify a point on the constraint boundary that is closest to the solution, known as the "most probable point" (MPP) of failure. Then, the task consists in finding this point efficiently. Typically, first or second-order reliability methods (FORM/SORM)~\cite{cizelj1994application} are utilized. These methods transform the uncertainty vectors into uncorrelated Gaussian random vectors using the Rosenblatt or Nataf transformation~\cite{lebrun2009generalization}, then the constraints are approximated linearly or quadratically. Therefore, the probabilistic constraints in Problem (\ref{problem_rbdo}) are reformulated as a deterministic optimization problem, reducing the task of solving Problem (\ref{problem_rbdo}) to two nested deterministic optimizations. Various approaches have been employed to solve it with a double loop, such as the Performance Measure Approach (PMA) or the Reliability Index Approach (RIA)~\cite{Aoues10}, a single loop, such as the Single Loop Approach (SLA)~\cite{liang2007single}, or decoupled approaches like the Sequential Optimization and Reliability Assessment (SORA) approach~\cite{du2004sequential} or the Sequential Approximate Programming (SAP) approach~\cite{cheng2006sequential}. These methods prove to be efficient even when dealing with nonlinear problems, and when gradients are approximated using finite differences~\cite{Aoues10}. Additionally, methods known as reliability-based robust design optimization (RBRDO) have been developed to handle uncertainties in the objective function by employing a bi-objective formulation of the problem~\cite{motta2016efficient}.

However, a major drawback of FORM-based methods is their reliance on linear approximations of the objective and constraint functions. These approximations can be inaccurate in practice if the underlying problem is not smooth. Therefore, other methods have been developed that do not rely on linear approximations. Similar to FORM-based methods, these approaches generally use a double-loop strategy. In the inner loop, a reliability analysis  estimates the feasible probability. Examples of such methods include important sampling~\cite{chaudhuri2020information, yuan2014efficient}, line sampling~\cite{de2015advanced}, subset simulation algorithms~\cite{au2007application}, or surrogate modeling strategies~\cite{li2010evaluation, peherstorfer2018multifidelity}. Subsequently, the estimation of the feasible probability is incorporated into the RBDO problem, resulting in a deterministic problem if the objective function is unnoised or a linear approximation of the objective function can be made.

In addition to the linear approximation, FORM-based methods suffer from another major drawback: they depend on the precise characterization of the uncertainty model of the variables and parameters (required for applying the Nataf transformation). However, the Nataf transformation cannot always be applied, especially when the blackbox inherently contains noise. Even when applicable, the Nataf transformation assumes a specific dependence structure of the uncertainties~\cite{lebrun2009innovating}. Nevertheless, in the absence of sufficient data, justifying and enforcing a specific dependency assumption becomes challenging and unwarranted due to its biasing effect on the final solution. The papers of R. Lebrun and A. Dutfoy~\cite{lebrun2009generalization, lebrun2009innovating} provide a detailed discussion of these issues related to using Nataf's transformation in FORM-based methods.

Uncertainties are commonly classified into two categories: aleatory uncertainties and epistemic uncertainties~\cite{rocchetta2018we}. Aleatory uncertainties represent the stochastic behavior and randomness of events and variables. Epistemic uncertainty is generally associated with a lack of knowledge about phenomena, imprecision in measurements, and poorly designed models. Aleatory uncertainties can be modeled by random variables, while epistemic uncertainties can be represented by interval or point data. Using probabilistic models for epistemic uncertainties may lead to infeasible designs in practice~\cite{rocchetta2018we}. Even for aleatory uncertainties, selecting an appropriate probabilistic model can be challenging, especially when the dimension of the uncertainties is large or when dependencies are unknown due to data scarcity~\cite{rocchetta2018we}. A poorly chosen model can result in underperforming designs or designs with significant failures~\cite{rockafellar2010buffered}. When epistemic uncertainties are involved in reliability analysis, non-probabilistic approaches based on evidence theory~\cite{shafer1976mathematical}, possibility theory~\cite{dubois2001possibility}, or fuzzy sets~\cite{liu2019evidence, zadeh1978fuzzy} may be used. 

Recently, some approaches have utilized ellipsoidal sets to model uncertainties~\cite{meng2018new, wang2019structural}. When both types of uncertainties are present, combining probabilistic and non-probabilistic models to address these uncertainties may be an interesting option~\cite{eldred2011mixed, meng2020new}. Alternatively, distributionally robust chance-constrained programming~\cite{xie2021distributionally} or a Bayesian probabilistic approach using Gaussian processes~\cite{amri2021sampling, nannapaneni2016reliability} also appear promising. Finally, scenario optimization, that tackles the problem (\ref{problem_ref}) using available data without prescribing a specific model (or a set of models) for the uncertainty, has been explored~\cite{rocchetta2021scenario}. Unfortunately, the described approaches are primarily used for reliability analysis, and they do not handle uncertainties in the objective function, except in the work in~\cite{amri2021sampling}, which is limited to parameter uncertainties. Another significant drawback is the lack of a convergence proof to an optimal point of the problem.~\Cref{table_methods} summarizes the different methods based on several criteria. The first two criteria assess whether the methods may deal with nonsmooth problems, while the third evaluates the ability of the method to handle noise in the objective function as well as in the constraints. The fourth criterion examines whether the method requires a precise characterization of the distribution that models the aleatory uncertainties, (e.g. for applying the Nataf transformation). Finally, the last criterion assesses the capability of the method to handle uncertainties in the absence of perfect knowledge of the data.

\begin{table}[ht!]
    \centering
    \caption{Summary of the different methods and their limits}
    \footnotesize
    \begin{threeparttable}
    \begin{tabular}{|C{2.7cm}||C{1cm}|C{2.5cm}|C{2.5cm}|C{2.5cm}|C{2.5cm}|C{2.5cm}|}
    \hline
         Methods & Type \tnote{1}   & Handles nonsmooth constraints
         & Handles nonsmooth objective & Handles noisy objective  &Allows unknown aleatory uncertainty & Allows lack of data \tnote{2} \\
    \hline
    \hline
    FORM-based~\cite{Aoues10, cheng2006sequential, du2004sequential, liang2007single} & O & \xmark  &\xmark & \xmark &  \xmark& \xmark \\
    \hline
    RBRDO~\cite{motta2016efficient} & O & \xmark  &\checkmark & \checkmark & \xmark& \xmark \\
    \hline 
    Importance Sampling~\cite{chaudhuri2020information, yuan2014efficient} & O & \checkmark &\checkmark & \xmark& \xmark& \xmark \\
    \hline
     Line Sampling~\cite{de2015advanced} & RA & \checkmark  & N/A & N/A & \xmark& \checkmark\\
    \hline
     Subset simulation~\cite{au2007application} & RA  & \checkmark  & N/A & N/A & \checkmark & \xmark\\
     \hline 
     Surrogate modelling~\cite{li2010evaluation, peherstorfer2018multifidelity} & RA  & \checkmark & N/A & N/A & \checkmark &\xmark \\
     \hline
     Mixed approaches~\cite{eldred2011mixed, meng2020new}  & O  & \xmark & \xmark & \xmark & \xmark & \checkmark \\
     \hline
     Ellipsoidal set~\cite{meng2018new, wang2019structural} & O  & \xmark & \xmark& \xmark & \xmark & \checkmark\\
     \hline
     Bayesian approach (I)~\cite{ nannapaneni2016reliability} & RA & \checkmark & N/A & N/A & \checkmark & \checkmark \\
     \hline
     Bayesian approach (II)~\cite{amri2021sampling} & O & \checkmark & \checkmark & \checkmark \tnote{3} & \checkmark & \xmark\\
     \hline 
     Scenario Optimization~\cite{rocchetta2021scenario} & O & \xmark & \xmark & \xmark & \checkmark & Only point data\\
     \hline
     \textbf{This work} & O & \checkmark & \checkmark &  \checkmark & \checkmark & \checkmark \tnote{4} \\
     \hline
    \end{tabular}
    \begin{tablenotes}
    \item[1] The type indicates if the method handle the whole stochastic constrained optimization problem (O) or is limited to reliability analysis (RA).
    \item[2] Only points or interval data are available.
    \item[3] Only parameters uncertainties.
    \item[4] For interval data, the method allows only to obtain worst-case solution.
  \end{tablenotes}
    \label{table_methods}
\end{threeparttable}
\end{table}

\subsection{Contributions}

To account for the uncertainties in both the objective and constraint functions, methods utilizing the conditional value-at-risk ($\cvar$) have been developed~\cite{li2020risk,Rock99}. $\cvar$ is a coherent risk measure that evaluates the risk associated with a design solution by combining the probability of undesired events with a measure of the magnitude or severity of those events. $\cvar$ methods have found extensive applications in risk averse optimization like in trust-region algorithms~\cite{menhorn2017trust}, in engineering design problems~\cite{heinkenschloss2018conditional,li2020new, tyrrell2015engineering, xu2021cvar}, and in constrained reinforcement learning ~\cite{chow2017risk, tamar2015optimizing}.

One of the main interest of the $\cvar$ measure lies in the flexibility provided by the parameter $\alpha$. When $\alpha = 0$, the $\cvar$ measure corresponds to the expectation, whereas as $\alpha$ approaches $1$, it corresponds to the supremum of the function over the support of the uncertainties~\cite{rockafellar2015risk}. This versatility allows to handle both aleatory and epistemic uncertainties, albeit in a worst-case scenario only. However, substituting failure probability constraints with $\cvar$ constraints is a conservative approach~\cite[chapter 6]{sha21} that might render the problem infeasible in the worst case. Moreover, the closer the value of $\alpha$ is to $1$, the more sensitive the measure becomes to the uncertainty model, particularly in the tails. Managing this heightened sensitivity necessitates an untractable number of samples. While the former issue is challenging to avoid a priori, the latter can be partially addressed by employing a multi-timescale stochastic approximation algorithm to estimate the $\cvar$ value~\cite{chow2017risk, prashanth2014policy}. Unfortunately, the methods utilized in the referenced papers cannot be directly applied to solve a $\cvar$ formulation of the problem (\ref{problem_ref}). In fact, these methods cleverly leverage the properties of the Markov Decision Process to compute estimates of the gradients, a strategy that is impossible to use in the context of the present study. The contributions of this work are outlined as follows.


First, in~\Cref{sec_smoo_prob}, the process of smoothing the problem and obtaining analytical gradient estimates from noisy measurements of the  blackbox is described. A smooth approximation of the gradient~\cite{bhatnagar2013stochastic, nesterov2017random} is employed. The concept involves approximating the original function by its convolution with a multivariate density function. The resulting approximation possesses several desirable properties: it is infinitely differentiable even if the original function is only piecewise continuous, it preserves the structural properties (such as convexity and Lipschitz constant) of the original function, and an unbiased estimator of the gradient of the smooth approximation can be calculated from only two measurements of the blackbox. In most studies~\cite{ghadimi2013stochastic, nesterov2017random}, Gaussian or uniform density functions are utilized for the approximation. However, in this paper, a truncated Gaussian density function is developed to satisfy the bound constraints of the problem (\ref{problem_ref}). The properties of this new approximation and its associated unbiased gradient estimator are provided.

Second, Problem (\ref{problem_ref}) is reformulated as a $\cvar$-constrained problem, wherein the objective function and the constraints are approximated by their smooth truncated Gaussian counterparts. The quality of this approximation is theoretically examined and depends on several parameters such that the value of $\alpha$, the dimension and the value of the smoothing parameter. Subsequently, following the approach in~\cite{chow2017risk}, a Lagrangian relaxation is applied to the problem. The method used to solve the relaxed problem is developed in~\Cref{sec_algo}. It involves a four-timescale stochastic approximation algorithm. The first timescale aggregates information about the gradient, the second estimates the quantile of the objective and constraint functions, the third updates the design variables in a descent direction, and the last one updates the Lagrange multiplier in the ascent direction. The convergence analysis of this algorithm is studied in~\Cref{sec_analysis} and is conducted using an Ordinary Differential Equation (ODE) approach. Under mild assumptions, this algorithm almost surely converges to a feasible point of the $\cvar$-constrained problem whose objective function value is arbitrarily close to that of a local solution.

Finally, in~\Cref{sec_exp}, practical implementation details are provided to minimize the number of hyperparameters in the developed algorithm. Numerical experiments are conducted to estimate the values of the remaining hyperparameters. Then, comparisons are made between the algorithm using the Gaussian gradient estimator and its truncated counterpart. In the last subsection, the efficiency of the algorithm is demonstrated on problems involving mixed aleatory/epistemic uncertainties. Conclusions are drawn in~\Cref{sec_concl}.

\section{Problem formulation}
\label{sec_prob_stat}

 In order to formally settle the problem and to develop the convergence analysis, the following assumptions are made on the functions $C_j$ and used throughout the paper.
 \begin{assumption}
 \label{assum1}
     Let $(\Omega, \mathcal{F}, \mathbb{P})$ be a probability space and consider $C_j(\x, \boldsymbol{\xi}) : \R^n \times \R^d \to \R, j \in [0, m]$ where $\boldsymbol{\xi}: \Omega \to \Xi \subset \R^d $ is the vector modelling the uncertainties. Then, the following hold for all $j \in [0, m]$.
     \begin{enumerate}
         \item There exists a measurable function $\kappa_1(\boldsymbol{\xi}) : \Xi \to \R$  such that $\mathbb{E}_{\boldsymbol{\xi}}[\kappa_1(\boldsymbol{\xi})] \leq L_1 < \infty$ and for which 
         $$|C_j(\x, \boldsymbol{\xi})| \leq \kappa_1(\boldsymbol{\xi}), \  \forall \x \in \mathcal{X}  \mbox{ and } \boldsymbol{\xi} \in \Xi.$$
         \item There exists a measurable function $\kappa_2(\boldsymbol{\xi}_1, \boldsymbol{\xi}_2) : \Xi\times\Xi \to \R$  where $\boldsymbol{\xi}_1$ and $\boldsymbol{\xi}_2$ are i.i.d. random vectors such that $\mathbb{E}_{\boldsymbol{\xi}}[\kappa_2(\boldsymbol{\xi}_1, \boldsymbol{\xi}_2)] \leq L_2 <\infty$ and for which 
         \begin{equation*}
             |C_j(\x, \boldsymbol{\xi}_1)  - C_j(\y, \boldsymbol{\xi}_2) | \leq \kappa_2(\boldsymbol{\xi})||\x - \y||, \; \forall (\x, \y) \in \mathcal{X} \times \mathcal{X} \mbox{ and }  (\boldsymbol{\xi}_1, \boldsymbol{\xi}_2 ) \in \Xi \times \Xi.
         \end{equation*}
         \item The function $C_j(\cdot, \boldsymbol{\xi}) $ has a continuous cumulative distribution function and there exists a measurable function $\kappa_3(\boldsymbol{\xi}_1, \boldsymbol{\xi}_2) : \Xi \times \Xi \to \R$, where $\boldsymbol{\xi}_1$ and $\boldsymbol{\xi}_2$ are i.i.d. random vectors such that $\mathbb{P}_{\boldsymbol{\xi}}(\kappa_3(\boldsymbol{\xi}_1, \boldsymbol{\xi}_2) \leq L_3 ) = 1$ with $L_3 < \infty$ and for which
         \begin{equation*}
             | C_j(\x, \boldsymbol{\xi}_1 ) - C_j(\x, \boldsymbol{\xi}_2 ) | \leq \kappa_3 ( \boldsymbol{\xi}_1, \boldsymbol{\xi}_2) ||\x - \y ||,  \; \forall (\x, \y) \in \mathcal{X} \times\mathcal{X} \mbox{ and } ( \boldsymbol{\xi}_1, \boldsymbol{\xi}_2 ) \in \Xi \times \Xi.
         \end{equation*}
     \end{enumerate} 
\end{assumption}
Three comments on these assumptions. First, note that no assumptions are made about the differentiability of the functions $C_j$. Second,~\Cref{assum1}.1 will be made throughout this paper because it allows the value-at-risk ($\var$) and the $\cvar$ of the functions $C_j$ to be well defined. The other assumptions are used in Section \ref{sec_smoo_prob} to bound the approximation of the constrained $\cvar$ blackbox problem and in Section \ref{sec_analysis} to study the convergence of the proposed method. Finally, the assumptions are increasingly strong, i.e.,~\Cref{assum1}.3 implies~\Cref{assum1}.2, which implies~\Cref{assum1}.1. 

Now, the $\var$ at level $\alpha \in (0,1)$ of the objective and constraint functions may be defined. It is originally derived from the left-side quantile of level $\alpha$ of a given random variable. Given $j \in  \{0,1,\ldots, m\}$ and a reliability level $\alpha_j \in (0, 1)$, the $\var$ of a function $C_j(\x, \boldsymbol{\xi})$ is defined as 
  \begin{equation*}
     \mathrm{VaR}_{\alpha_j}(\x) \ := \ 
        \inf \{ t \,|\, \mathbb{P}(C_j(\x,\boldsymbol{\xi}) \leq t) \geq  \alpha_j \}.
 \end{equation*}
 The $\var$ of a function has several interesting properties. When the cumulative distribution function $\mathbb{P}(C_j(\x,\boldsymbol{\xi}) \leq u)$ is right continuous with respect to $t$, the infemum is a minimum and if it is, in addition, continuous and strictly increasing, then $\mathrm{VaR}_{\alpha_j}$ is the unique $t$ such that $\mathbb{P}(C_j(\x,\boldsymbol{\xi}) \leq t) = \alpha$. However, the $\var$ of a function is computationally intractable, is not a coherent risk measure~\cite{Art97} and does not take into account the magnitude/severity of the undesired events. Therefore, in practice another measure is used: the Conditional Value-at-Risk. The $\cvar$ of a function $C_j(\cdot, \boldsymbol{\xi})$, for a level $\alpha_j \in (0, 1)$ at a point, $\x$ may be defined as~\cite{Rock99}
\begin{equation}
\label{cvar}
    \mathrm{CVaR}_{\alpha_j}(\x) \ := \ \min_{t \in \mathbb{R}} V_{\alpha_j}(\x, t),
\end{equation}
where 
\begin{equation}
     V_{\alpha_j}(\x, t) \ = \ t + \frac{1}{1-\alpha_j} \mathbb{E}_{\boldsymbol{\xi}}[(C_j(\x, \boldsymbol{\xi}) - t)^+], 
\end{equation}
 where the superscript plus denotes the function $(t)^+ := \max \{0, t \}$. The level $\alpha_j$ gives the possibility to choose the desired degree of reliability. Choosing a level close to $0$ is tantamount to taking the expectation measure into account, i.e. adopting a "risk neutral" approach. On the other hand, choosing a level close to $1$ is tantamount to taking a "worst-case" approach. In this way, different values of $\alpha_j$ can be used for the different objective and constraint functions, depending on the degree of reliability desired for each of them. Now, problem (\ref{problem_ref}) can be reformulated as a  $\cvar$-constrained blackbox optimization problem:
\begin{align}
\label{problem_cvar}
\begin{split}
\begin{array}{cl}
      \displaystyle \min_{\x \in \mathcal{X}}   &\mathrm{CVaR}_{\alpha_0}(\x)\\
      \mbox{s.t.} &\mathrm{CVaR}_{\alpha_j} (\x) \leq 0, \; \forall j \in [1,m].
\end{array}
\end{split}
\end{align}
This formulation is a convex program if the
objective and constraint functions are convex in the design space. This convexification of the design space makes Problem (\ref{problem_cvar}) a conservative approximation of Problem (\ref{problem_rbdo}) [Chapter 6,~\cite{sha21}]. 
Thus, this formulation guarantees a conservative result in terms of failure probability, see e.g.~\cite{rockafellar2010buffered}. To solve Problem (\ref{problem_cvar}), it is usually reformulated with the function $V_\alpha$ as follows
\begin{align}
\label{problem_cvar_bis}
\begin{split}
\begin{array}{cl}
     \displaystyle \min_{(\x, \t) \in \mathcal{X} \times \R^{m+1}}   & V_{\alpha_0} (\x, t_0) \\
     \mbox{s.t.}   & V_{\alpha_j}(\x, t_j) \leq 0, \; \forall j \in [1,m].
\end{array}  
\end{split}
\end{align}
The equivalence between Problem (\ref{problem_cvar}) and Problem (\ref{problem_cvar_bis}) is shown in the following lemma.
\begin{Lem}
   Suppose the solution sets of Problem (\ref{problem_cvar}) and Problem (\ref{problem_cvar_bis}) are not empty. Then these problems are equivalent in the sense that, 
   $\x^*$ is a solution of Problem (\ref{problem_cvar}) if and only if there exist $\t^* \in \R^{m+1}$ such that $(\x^*, \t^*)$ is a solution of Problem (\ref{problem_cvar_bis}), and the optimal values  are the same.
\end{Lem}
\begin{proof}
    By the definition of the Conditional Value-at-Risk given in~\Cref{cvar}, Problem (\ref{problem_cvar}) may be reformulated as follows
    \begin{align}
    \label{prob_cvar_reform}
    \begin{split}
    \begin{array}{cl}
         \displaystyle \min_{\x \in \mathcal{X}}  \Big( \min_{t_0 \in \R}   & V_{\alpha_0} (\x, t_0)  \Big) \\
         \mbox{s.t.}    \displaystyle  \Big( \min_{t_j \in \R} &V_{\alpha_j}(\x, t_j) \Big) \leq 0, \; \forall j \in [1,m].
    \end{array}  
    \end{split}
    \end{align}

    Now, the following relations hold
    \begin{align*}
        &\min_{\x \in \mathcal{X}} \Big( \min_{t_0 \in \R}   V_{\alpha_0} (\x, t_0) \Big) = \min_{(\x, t_0) \in \mathcal{X} \times \R}    V_{\alpha_0} (\x, t_0) \\
        & \Big( \min_{t_j \in \R} V_{\alpha_j}(\x, t_j) \Big) \leq 0, \; \forall j \in [1,m] \iff \forall j, \; \exists t_j \mbox{ s.t. } V_{\alpha_j}(\x, t_j) \leq 0.
    \end{align*}
    Therefore, the Problems (\ref{prob_cvar_reform}) and (\ref{problem_cvar_bis}) are equivalent. Now, let $\x^*$ be a solution of Problem (\ref{problem_cvar}), it is possible to construct the associated vector $\t^*(\x^*)$ where $t_j^*(\x^*) = \mathrm{VaR}_{\alpha_j}(\x^*), \forall j \in [0, m]$. The tuple $(\x^*, \t^*(\x^*))$ is then solution of Problem (\ref{prob_cvar_reform}) and as a consequence of Problem (\ref{problem_cvar_bis}) which ends the proof.
\end{proof}

Despite this property, Problem (\ref{problem_cvar}) is difficult to solve for two main reasons. First, since the functions $C_j$ are the outputs of a blackbox, the gradients of these functions may not exist, and even if they do, their analytic formulations are not available. 
Second, the problem is highly sensitive to the values of $\alpha_j$, and the closer the values are to $1$, the harder the problem is to solve. 
The next section describes the strategy used in this paper to overcome these difficulties.

\section{Smooth approximation and Lagrangian relaxation of the problem}
\label{sec_smoo_prob}
This section introduces a method for solving the Problem (\ref{problem_cvar}). To obtain a more tractable problem, the original problem is approximated by a smooth problem using truncated Gaussian smoothing. The quality of the approximation is then studied and a Lagrangian relaxation of the smooth problem is given.

\subsection{Truncated Gaussian smooth approximation}
In a blackbox optimization framework, all we know is that for any given input, the blackbox will return an output, which may be subject to uncertainties.  To obtain a more tractable problem, a smooth approximation may be used~\cite[pp.~263]{rubinstein_simulation_1981}.
The principle of this method is to approximate the function by its convolution with a kernel density function. Formally, if $c$ is an integrable function, $\beta > 0$ is a scalar, and $\u$ is a random vector with distribution $\phi$, the smooth approximation of $c$ can be defined as
\begin{equation}
    c^{\beta}(\x)\  := \ \int_{-\infty}^{+\infty} c(\x -\beta \u) \phi(\u) d\u
    \ = \ \mathbb{E}_\u[c(\x + \beta \u)].
\end{equation}
The smooth approximation benefits from several attractive properties. First, it can be interpreted as a local weighted average of the function values in the neighborhood of $\x$. If $c$ is continuous at $\x$, it is possible to obtain a value of $c^\beta(\x)$ that is arbitrarily close to the value of $c(\x)$ by using an appropriate value of $\beta$. Second, it inherits the degree of smoothness of the density function as a consequence of the convolution product. Finally, depending on the chosen kernel, stochastic gradient estimators can be computed. They are unbiased estimators of the gradient of $c^\beta$ and can be constructed only from values of $c(\x)$ and $c(\x + \beta \u)$.

The most commonly used kernels are the Gaussian distribution and the uniform distribution on a sphere~\cite{nesterov2017random, ghadimi2013stochastic}. However, if the problem has bound constraints, a significant drawback of these distributions is that the random vector $\x + \sigma \u$  may fall outside the bound constraints. For instance, if $\u \sim \mathcal{N}(0, 1)$, $\x + \sigma \u$ might be sampled outside the bounds. This issue persists even with a uniform distribution if $\x$ is near the bounds. However, the bound constraints are usually non-relaxable in the sense of~\cite{digabel2015taxonomy}, meaning that the output of the blackbox lacks significance for optimization outside the bound constraints. This can occur due to physical phenomena or when the blackbox is undefined beyond the bounds. In such cases, the gradient estimate of $c^\beta$, computed from the values of the function $c$ at the points $\x$ and $\x+ \beta\x$, becomes unreliable. To address this issue, a truncated Gaussian estimator is developed in this paper, and its main properties are summarized in the following lemma.

\begin{Lem}
\label{lem_trunc_gaus}
Let $c$ be an integrable function on $\mathcal{X}$, the smooth approximation $c^\beta$ is defined as
\begin{equation*}
    c^\beta(\x)\  = \ \mathbb{E}_\u[c(\x + \beta \u) ],
\end{equation*}
 where $\u \sim \mathcal{TN}(\mathbf{0}, \mathbf{I}, \frac{\b_\ell-\x}{\beta}, \frac{\b_u -\x}{\beta})$,  $\b_\ell$ and $\b_u$ are respectively the lower and the upper bounds of the problem. In what follows, $\phi$ and $\Phi$ denote respectively the probability density function (p.d.f.) and the cumulative density function (c.d.f.) of the standard Gaussian distribution. Now, the following holds.
\begin{enumerate}
    \item $c^\beta$ is  infinitely differentiable: $c^\beta \in \mathcal{C}^\infty$.
    \item A one-sided unbiased estimator of $\nabla c^\beta$ is 
    \begin{equation}
        \Tilde{\nabla} c^\beta(\x)\  =\  \frac{(\u - \boldsymbol{\mu})c(\x + \beta \u) - (\u -  \boldsymbol{\mu})c(\x)}{\beta},
        \label{estim1}
    \end{equation}
    where $\boldsymbol{\mu}$ is the mean of the truncated Gaussian vector, i.e,
    \begin{equation*}
        \mu_i\  =\  \frac{\phi \left(\frac{b_{\ell_i} -x_i }{\beta} \right) - \phi\left(\frac{b_{u_i} -x_i }{\beta}\right)}{\Phi\left(\frac{b_{u_i} - x_i}{\beta} \right) - \Phi \left(\frac{b_{\ell_i} - x_i}{\beta}\right)},\quad  \forall i \in [1,n].
    \end{equation*}
    \item Let $\u_1 \sim \mathcal{TN}(\mathbf{0}, \mathbf{I}, \frac{\b_\ell - \x}{\beta}, \frac{\b_u - \x}{\beta})$ and $\u_2 \sim \mathcal{TN}(\mathbf{0}, \mathbf{I}, \frac{\x - \b_u}{\beta}, \frac{\x - \b_\ell}{\beta})$, a two-sided unbiased estimator of $\nabla c^\beta$ is 
     \begin{equation}
        \Tilde{\nabla} c^\beta(\x) \ =\  \frac{(\u_1 - \boldsymbol{\mu}_1)(c(\x + \beta \u_1) - c(\x)) - (\u_2 -  \boldsymbol{\mu}_2)(c(\x - \beta \u_2) - c(\x))}{2 \beta},
        \label{estim2}
    \end{equation}
    \item In addition, if $c$ is a L-Lipschitz continuous function, let $\beta \geq 0$, then $\forall \x \in \mathbb{R}^n$
    \begin{equation*}
        |c^\beta(\x) - c(\x)| \ \leq \ L \beta \sqrt{n}.
    \end{equation*}
\end{enumerate}
\end{Lem}

\begin{proof}
1.) 
This can be shown by noting that the truncated Gaussian kernel is infinitely differentiable within the bounds. 
However, to obtain the above estimators, the calculation must be done. Therefore, using the above notation, and given that the components $u_i$ of $\u$ are mutually independent, it follows that
\begin{align*}
    \mathbb{E}_\u[c(\x+\beta \u)] &= \int_{\frac{\b_\ell-\x}{\beta}}^{\frac{\b_u-\x}{\beta}} c(\x + \beta \u) \prod_{i = 1}^n \frac{\phi(\u_i)}{\Phi \left(\frac{\b_{u_i} - \x_i}{\beta}\right) - \Phi\left(\frac{\b_{\ell_i} - \x_i}{\beta}\right)} d \u \\
      &= \int_{\frac{\b_\ell-\x}{\beta}}^{\frac{\b_u-\x}{\beta}} \frac{1}{(2 \pi)^{\frac{n}{2}}}c(\x + \beta \u) \prod_{i = 1}^n \frac{e^{-\frac{u_i^2}{2}}}{\Phi \left(\frac{\b_{u_i} - \x_i}{\beta}\right) - \Phi\left(\frac{\b_{\ell_i} - \x_i}{\beta}\right)} d \u \\
    &= \frac{1}{(2 \pi)^{\frac{n}{2}}} \left( \prod_{i = 1}^n\frac{1}{\Phi \left(\frac{\b_{u_i} - \x_i}{\beta}\right) - \Phi\left(\frac{\b_{\ell_i} - \x_i}{\beta}\right)} \right) \int_{-\infty}^\infty \mathbf{1}_{\left[\frac{\b_\ell-\x}{\beta},\frac{\b_u-\x}{\beta}\right] } (\u) c(\x + \beta \u) \prod_{i = 1}^n e^{-\frac{u_i^2}{2}} d \u ,
\end{align*}
where $\mathbf{1}_{[\cdot]}(\cdot)$ denotes the indicator function. 
Substituting $\v = \x + \beta \u$ leads to:
\begin{equation*}
    \mathbb{E}_\u[c(\x+\beta \u)] = \frac{1}{(2 \pi )^{\frac{n}{2}} \beta^n } \left(\prod_{i = 1}^n\frac{1}{\Phi \left(\frac{\b_{u_i} - \x_i}{\beta}\right) - \Phi\left(\frac{\b_{\ell_i} - \x_i}{\beta}\right)} \right) \int_{-\infty}^\infty \mathbf{1}_{\left[\b_\ell, \b_u \right] } (\v) c(\v) \prod_{i = 1}^n e^{-\frac{(x_i - v_i)^2}{2 \beta^2}} d \v.
\end{equation*}
By setting
\begin{align*}
    h_1(\x) &= \frac{1}{(2 \pi)^{\frac{n}{2}} \beta^n} \prod_{i = 1}^n\frac{1}{\Phi \left(\frac{\b_{u_i} - \x_i}{\beta}\right) - \Phi\left(\frac{\b_{\ell_i} - \x_i}{\beta}\right)}, \\
    h_2(\x) &= \mathbf{1}_{\left[\b_\ell, \b_u \right] } (\x) c(\x)
    \qquad \mbox{ and } \qquad 
    h_3(\x) = \prod_{i = 1}^n e^{-\frac{(x_i)^2}{2 \beta^2}},
\end{align*}
 $c^\beta(\x)$ may be compactly written as
\begin{equation*}
    c^\beta(\x) = h_1(\x)(h_2*h_3)(\x),
\end{equation*}
where $*$ is the convolution product between two functions. As $h_3 \in \mathcal{C}^\infty(\mathbb{R}^n)$ and $h_2 \in \mathcal{L}^1(\Omega,\mathcal{F}, \mathbb{P})$ then $(h_2*h_3) \in \mathcal{C}^\infty(\mathbb{R}^n)$ (property of convolution product). 
Moreover, $h_1 \in \mathcal{C}^\infty(\mathbb{R}^n)$ as well, therefore $c^\beta(\x) \in \mathcal{C}^\infty(\mathbb{R}^n)$ as it is the product of infinitely continuously differentiable functions. \\

2.) By using the same notation as above, the partial derivative of $c^\beta$ may be computed, for $j \in [1, n]$ as
\begin{align*}
    \frac{\partial c^\beta(\x)}{\partial x_j} = \frac{\partial h_1(\x)}{\partial x_j} \left(h_2*h_3 \right)(\x) + h_1(\x)\left(h_2*\frac{\partial h_3}{\partial x_j}\right)(\x).
\end{align*}
Yet, we have
\begin{align*}
    \frac{\partial h_1(\x)}{\partial x_j} &= \frac{1}{(2 \pi)^{\frac{n}{2}} \beta^n} \left( \displaystyle\prod_{i = 1}^n\frac{1}{\Phi \left(\frac{\b_{u_i} - \x_i}{\beta}\right) - \Phi\left(\frac{\b_{\ell_i} - \x_i}{\beta}\right)} \right) \frac{\phi\left(\frac{b_{u_j} -x_j }{\beta}\right) - \phi\left(\frac{\b_{\ell_j} -x_j }{\beta}\right)}{\beta\left(\Phi \left(\frac{\b_{u_j} - \x_j}{\beta}\right) - \Phi\left(\frac{\b_{\ell_j} - \x_j}{\beta}\right) \right)} = -\frac{\mu_j h_1(\x)}{\beta} \\
    \frac{\partial h_3(x)}{\partial x_j} &= -\frac{x_j }{\beta^2} h_3(\x).
\end{align*}
Thus, we obtain
\begin{equation*}
    \frac{\partial c^\beta(\x)}{\partial x_j} = \mathbb{E}_\u \left[\frac{u_j - \mu_j}{\beta} c(\x + \beta \u) \right].
\end{equation*}
From this result, an unbiased estimator of the gradient of $c^\beta$ is
\begin{equation*}
    \Tilde{\nabla} c^\beta(\x) = \frac{\u - \boldsymbol{\mu}}{\beta} c(\x + \beta \u).
\end{equation*}
As the variance of this estimator gets unbounded as $\beta$ goes to $0$, in practice the following estimator is used
\begin{equation*}
     \Tilde{\nabla} c^\beta(\x) = \frac{(\u - \boldsymbol{\mu})c(\x + \beta \u) - (\u -  \boldsymbol{\mu})c(\x)}{\beta}. 
\end{equation*}
This estimator is still unbiased since $\mathbb{E}_\u[(\u - \boldsymbol{\mu} )c(\x)] = 0$.\\

3.) Symmetrically, if $\u \sim \mathcal{TN}(\mathbf{0}, \mathbf{I}, \frac{\x - \b_u}{\beta}, \frac{\x - \b_\ell}{\beta})$, an unbiased estimator is 
\begin{equation*}
     \Tilde{\nabla} c^\beta(\x) = \frac{(\boldsymbol{\mu} -  \u )c(\x - \beta \u) - (\boldsymbol{\mu} -  \u)c(\x)}{\beta}. 
\end{equation*}
thus, by summation of the two one-sided estimator, the two-sided estimator is obtained. \\

4.) Finally, we have, with $\u_1 \sim \mathcal{TN}(\mathbf{0}, \mathbf{I}, \frac{\b_\ell - \x}{\beta^1},\frac{\b_u - \x}{\beta^1})$ and $\u_2 \sim \mathcal{TN}(\mathbf{0}, \mathbf{I}, \frac{\b_\ell - \x}{\beta^2},\frac{\b_u - \x}{\beta^2})$
\begin{align*}
    |c^\beta(\x) - c(\x)| &=  |\mathbb{E}_{\u} [ c(\x + \beta \u)] - c(\x)| \leq  \mathbb{E}_{\u}[ |c(\x + \beta \u) - c(\x)|] \leq  L \beta \mathbb{E}_{\u}[|| \u||] .
\end{align*}
where the first inequality comes from the Jensen's inequality and the second one comes from the L-Lipschitz continuity of $c$. It remains to bound $\mathbb{E}_\u[||\u||]$ when $\u$ is a truncated Gaussian vector, for this purpose, the proof of Lemma 1 of~\cite{nesterov2017random} is adapted for truncated Gaussian distribution. The following identity is used:
\begin{equation*}
    \int_{\frac{\b_\ell-\x}{\beta}}^{\frac{\b_u -\x}{\beta}}  e^{-\frac{||u||^2}{2}} d\u  \ = \ 
    (2\pi)^{n/2} \prod_{i = 1}^n \left( \Phi \left(\frac{b_{u_i} - x_i}{\beta}\right) -\Phi \left(\frac{b_{\ell_i}- x_i}{\beta}\right) \right) := \kappa.
\end{equation*}
By setting $\v = \x + \beta \u$ and multiplying by $\beta^n$,
 the last equalities become
\begin{equation*}
 \int_{\b_\ell}^{\b_u}  e^{-\frac{||\v -\x||^2}{2\beta^2}} d\v \ = \ 
 (2\pi)^{n/2} \prod_{i = 1}^n \left(\Phi\left(\frac{b_{u_i} - x_i}{\beta}\right) -\Phi\left(\frac{b_{\ell_i} - x_i}{\beta}\right) \right) \beta^n \ = \ \kappa \beta^n.
\end{equation*}
Taking the logarithm yields
\begin{equation}
    \ln \left(  \int_{\b_\ell}^{\b_u}  e^{-\frac{||\v -\x||^2}{2\beta^2} } d\v \right) \ = \ 
    n \ln(\beta) + \frac{n}{2} \ln(2 \pi) + \sum_{i = 1}^n \ln \left( \Phi\left(\frac{b_{u_i} - x_i}{\beta}\right) - \Phi\left(\frac{b_{\ell_i} - x_i}{\beta}\right) \right).
    \label{equa_ln}
\end{equation}
Now, the derivative of the left-hand-side of ~\Cref{equa_ln} with respect to $\beta$ is given by
\begin{align*}
    \frac{\partial}{\partial \beta}\ln \left(  \int_{\b_\ell}^{\b_u}  e^{-\frac{||\v -\x||^2}{2\beta^2} } d\v \right) &=
    \frac1{\kappa \beta^n}\int_{\b_\ell}^{\b_u} \frac{||\v -\x||^2}{\beta^3} e^{-\frac{||\U -\x||^2}{2\beta^2} } d\v \\
    &= \frac{1}{\kappa \beta} \int_{\frac{\b_\ell-\x}{\beta}}^{\frac{\b_u-\x}{\beta}} ||\u||^2 e^{-\frac{||\u||^2}{2} } d\u \qquad \qquad 
    \mbox{ since } \frac{\v - \x}{\beta} = \u\\
    &= \frac{1}{\beta} \mathbb{E}_\u[||\u||^2]
\end{align*}

and the derivative of the right-hand-side of~\Cref{equa_ln} is given by
\begin{equation*}
    \frac{n}{\beta} + \sum_{i = 1}^n \frac{\frac{b_{\ell_i} - x_i}{\beta^2} \phi \left(\frac{b_{\ell_i} - x_i}{\beta}\right) - \frac{b_{u_i} - x_i}{\beta^2} \phi\left(\frac{b_{u_i} - x_i}{\beta}\right)}{\Phi\left(\frac{b_{u_i} - x_i}{\beta}\right) - \Phi\left(\frac{b_{\ell_i} - x_i}{\beta}\right)}.
\end{equation*}
Thus,
\begin{equation}
    \mathbb{E}_\u[||\u||^2] \ = \  
    n + \sum_{i = 1}^n \frac{\frac{b_{\ell_i} - x_i}{\beta} \phi\left(\frac{b_{\ell_i} - x_i}{\beta}\right) - \frac{b_{u_i} - x_i}{\beta} \phi\left(\frac{b_{u_i} - x_i}{\beta}\right)}{\Phi\left(\frac{b_{u_i} - x_i}{\beta}\right) - \Phi\left(\frac{b_{\ell_i} - x_i}{\beta}\right)}
    \ \leq \ n ,
    \label{bound}
\end{equation}
where the inequality holds because the sum is negative for $\x \in \mathcal{X}$. Finally, with the result in~\Cref{bound} and the results of Lemma 1 of~\cite{nesterov2017random}, the following bound appears
\begin{equation*}
    \mathbb{E}_\u[||\u ||] \leq \sqrt{n}.
\end{equation*}
\end{proof}
When only noisy outputs of the blackbox are available, the following estimator is used 
\begin{align}
    \label{estim_grad_noised}
    \Tilde{\nabla} C^\beta(\x, \boldsymbol{\xi}) \ = \
    \frac{(\u - \boldsymbol{\mu}) \left(C(\x + \beta \u, \boldsymbol{\xi}_1 ) - C(\x, \boldsymbol{\xi}_2)\right)}{\beta},
\end{align}
where $\boldsymbol{\xi}_1$ and  $\boldsymbol{\xi}_2$ are two independent identically distributed realizations of a random vector $\boldsymbol{\xi}$. This estimator is still unbiased because
\begin{equation*}
    \mathbb{E}_{\u, \boldsymbol{\xi}}[\Tilde{\nabla} c^\beta(\x, \boldsymbol{\xi})]\ = \ 
    \mathbb{E}_{\u}[\mathbb{E}_{\boldsymbol{\xi}}[\Tilde{\nabla} c^\beta(\x, \boldsymbol{\xi}) | \u] ] = \nabla c^\beta(\x).
\end{equation*}

\subsection{Smooth approximation of CVaR-constrained blackbox optimization problem}
The non-smoothness of a $\cvar$-constrained blackbox optimization problem arises from two elements: the potential non-smoothness of the functions $C_j$ and the non-smoothness introduced by the function $\max$ in the $\cvar$ formulation. The concept of smoothing a $\cvar$-constrained optimization problem is not novel; it has been explored in prior works~\cite{meng2011smoothing, soma2020statistical}. In this study, this concept is applied to both sources of non-smoothness using the aforementioned truncated Gaussian smoothing. As $\t$ is an unconstrained vector, arbitrarily large bounds are introduced for this vector. Let $\beta_1, \beta_2 > 0$ be two scalars, $\u \sim \mathcal{TN}(\mathbf{0}, \mathbf{I}, \frac{\b_\ell-\x}{\beta_1}, \frac{\b_u -\x}{\beta_1})$ a random vector of size $n$ and $\v \sim \mathcal{TN}(\mathbf{0}, \mathbf{I}, \frac{-\t_{\max}}{\beta_2}, \frac{\t_{\max}}{\beta_2})$, a random vector of size $m+1$, where $\t_{\max}$ is chosen to be sufficiently large, the smooth approximation of $V_{\alpha_j}$ and $\mathrm{CVaR}_{\alpha_j}$ for all $j \in [0, m]$ are defined respectively as
\begin{align*}
    V_{\alpha_j}^\beta(\x, t_j) &= \mathbb{E}_{\u, \v}[V_{\alpha_j}(\x + \beta_1 \u, t_j + \beta_2 v_j)], \; \text{ and }\\
    \mathrm{CVaR}_{\alpha_j}^\beta(\x) &= \min_{t_j \in \R} \mathbb{E}_{\u, \v}[V_{\alpha_j}(\x + \beta_1 \u, t_j + \beta_2 v_j)].
\end{align*}
 Then, the smooth approximation of the Problem (\ref{problem_cvar_bis}) may be formulated as follows
\begin{align}
\label{problem_cvar_smooth}
\begin{split}
\begin{array}{cl}
     \displaystyle \min_{(\x, \t) \in \mathcal{X} \times \R^{m+1}}   & V_{\alpha_0}^\beta(\x, t_0)   \\
     \mbox{s.t.}  &V_{\alpha_j}^\beta(\x, t_j) \leq 0, \; \forall j \in [1,m].
\end{array}
\end{split}
\end{align}
Now, the quality of this smooth approximation is studied.  The following Lemma states properties of the truncated Gaussian smoothing approximation applied with the $\cvar$ measure. 
\begin{Th}
\label{th_quali_approx}
    Under~\Cref{assum1}.2, the following holds.
    \begin{enumerate}
        \item $| \mathrm{CVaR}_{\alpha_j}^\beta(\x) -  \mathrm{CVaR}_{\alpha_j}(\x)| \leq \frac{L_2 \beta_1 \sqrt{n} + \beta_2}{1 - \alpha_j} $ for all $j \in [0, m]$ and $\x \in \mathcal{X}$; 
        \item $| \mathrm{CVaR}_{\alpha_0}^\beta(\Tilde{\x}^*) - \mathrm{CVaR}_{\alpha_0}(\x^*)| \leq \frac{L_2 \beta_1 \sqrt{n} + \beta_2}{1 - \alpha_j} $, where $ \Tilde{\x}^*$ and $\x^*$ are solutions of Problem (\ref{problem_cvar_smooth}) and (\ref{problem_cvar}) respectively.
        \item If~\Cref{assum1}.3 holds, then there exists a threshold $\Bar{\alpha}_j \in (0, 1]$ such that for all $\alpha_j \geq \Bar{\alpha}_j$
    \begin{equation*}
        \mathrm{CVaR}_{\alpha_j}(\x) \ \leq \ 
        \mathrm{CVaR}_{\alpha_j}^\beta(\x) \ \leq \ 
        \mathrm{CVaR}_{\alpha_j}(\x) +  \frac{L_2 \beta_1 \sqrt{n} + \beta_2}{1 - \alpha_j}.
    \end{equation*}
    \end{enumerate}
    Thus, for $\alpha_j \geq \Bar{\alpha}_j$, if $\Tilde{x}^*$ is a solution of Problem (\ref{problem_cvar_smooth}), then it is a feasible point for Problem (\ref{problem_cvar}).
\end{Th}

\begin{proof}
    1. Under Assumption~\Cref{assum1}.2, it follows that for all $(\x, t_j) \in \mathcal{X}\times \R$
    \begin{align*}
        |V_{\alpha_j}^\beta(\x, t_j) - V_{\alpha_j}(\x, t_j)| &= \frac{1}{1-\alpha_j} \left|\mathbb{E}_{\u, \v, \boldsymbol{\xi}} [(C_j(\x + \beta_1 \u, \boldsymbol{\xi}_1) - (t_j+ \beta_2 v_j))^+ - (C_j(\x, \boldsymbol{\xi}_2) - t_j)^+] \right| \\
        &\leq \frac{1}{1-\alpha_j} \mathbb{E}_{\u, \v, \boldsymbol{\xi}} [ |(C_j(\x + \beta_1 \u, \boldsymbol{\xi}_1) - (t_j+ \beta_2 v_j))^+ - (C_j(\x, \boldsymbol{\xi}_2) - t_j)^+|] \\
        &\leq \frac{1}{1-\alpha_j} \mathbb{E}_{\u, \v, \boldsymbol{\xi}} [ |C_j(\x + \beta_1 \u, \boldsymbol{\xi}_1) -  \beta_2 v_j - C_j(\x, \boldsymbol{\xi}_2) |] \\
        &\leq \frac{1}{1-\alpha_j} \mathbb{E}_{\u, \v, \boldsymbol{\xi}} [  \kappa_2(\boldsymbol{\xi}_1, \boldsymbol{\xi}_2) \beta_1 ||\u|| + \beta_2 |v_j|] \\
        &\leq \frac{L_2 \beta_1 \sqrt{n} + \beta_2}{1 - \alpha_j},
    \end{align*}
    where the first inequality follows from Jensen's inequality, the second from the following inequality $|\max(0, a) - \max(0, b)| \leq |a-b| $, the third from~\Cref{assum1}.2 and the last one from the independence of $\u$ and $\kappa_2(\boldsymbol{\xi})$ and the bound on the expectation of the norm of (truncated) Gaussian random vectors. This is true for all tuples $(\x, t_j) \in \mathcal{X}\times \R$, in particular for $t_j^* \in \argmin V(\x, t_j)$ and  $\Tilde{t}_j^* \in \argmin \mathbb{E}_{\u, \v}[V_{\alpha_j}(\x + \beta_1 \u, t_j + \beta_2 v_j)]$. Therefore, it follows that for any $j \in [0, m]$ and any $\x \in \mathcal{X}$
    \begin{equation*}
        V_{\alpha_j}^\beta(\x, \Tilde{t}_j^*) \ \leq \ 
        V_{\alpha_j}^\beta(\x , t_j^*) \ \leq \ 
        V_{\alpha_j}( \x , t_j^*) + \frac{L_2 \beta_1 \sqrt{n} + \beta_2}{1 - \alpha_j}.
    \end{equation*}
    Conversely, it also follows that
    \begin{equation*}
        V_{\alpha_j}(\x, t_j^*) \ \leq \ 
        V_{\alpha_j}(\x , \Tilde{t}_j^*) \ \leq \ 
        V_{\alpha_j}^\beta( \x , \Tilde{t}_j^*) + \frac{L_2 \beta_1 \sqrt{n} + \beta_2}{1 - \alpha_j}.
    \end{equation*}
    Recalling that $\mathrm{CVaR}_{\alpha_j}(\x) =  V_{\alpha_j}(\x, t_j^*)$ and $\mathrm{CVaR}_{\alpha_j}^\beta(\x) = V_{\alpha_j}^\beta(\x, \Tilde{t}_j^*)$, we obtain that 
    \begin{equation*}
        | \mathrm{CVaR}_{\alpha_j}^\beta(\x) -  \mathrm{CVaR}_{\alpha_j}(\x)| \ \leq \ 
        \frac{L_2 \beta_1 \sqrt{n} + \beta_2}{1 - \alpha_j} \; \forall j \in [1, m].
    \end{equation*} \\
    
    2. Using the same previous argument but with respect to $\x$ instead of $t$ allows to obtain the second inequality. \\

    3. Consider $\x \in  \mathcal{X}$ and suppose that~\Cref{assum1}.3 holds. It follows that for all $j \in [0, m]$ and $\boldsymbol{\xi} \in \Xi$ 
    \begin{equation*}
        |C_j(\x, \boldsymbol{\xi} )|-|C_j(\mathbf{0}, \mathbf{0})|\  \leq \ 
        |C_j(\x, \boldsymbol{\xi}) - C_j(\mathbf{0}, \mathbf{0})| \ \leq \ \kappa_3(\boldsymbol{\xi}, \mathbf{0}) ||\x||,
    \end{equation*}
    which implies that $|C_j(\x, \boldsymbol{\xi})| $ is almost surely bounded by a function depending on $\x$. Now, for all $\x \in \mathcal{X}$, $M_j(\x)$ is defined as the essential supremum of $C_j(\x, \boldsymbol{\xi})$, i.e,
    \begin{equation*}
         M_j(\x) := \inf \{ t \in \R \ | \ C_j(\x, \boldsymbol{\xi}) \leq t \text{ for almost every}\ \boldsymbol{\xi} \in \Xi \}.
    \end{equation*} 
    Now, we have by definition
    \begin{equation*}
        \mathrm{VaR}_{\alpha_j = 1} (\x)\  =  \ \inf \{t \; | \; \mathbb{P}(C_j(\x, \boldsymbol{\xi}) \leq t) = 1 \}\  =\  M_j(\x). 
    \end{equation*}
     As the c.d.f. of $C_j(\cdot, \boldsymbol{\xi})$ is assumed continuous, then it follows by~\cite{rockafellar2014random} that
    \begin{equation*}
        \mathrm{CVaR}_{\alpha_j}(\x) \ = \ 
        \frac{1}{1-\alpha_j} \int_{\alpha_j}^1 \mathrm{VaR}_\tau (\x) d \tau.
    \end{equation*}
As for $\tau \in [\alpha_j, 1]$, the $\mathrm{VaR}_\tau$ function is continuous with respect to $\tau$ with $\mathrm{VaR}_{\alpha_j} (\x) \leq \mathrm{VaR}_\tau(\x) \leq \mathrm{VaR}_{\alpha_j = 1}(\x) = M_j(\x)$, the mean value theorem ensures
\begin{equation*}
    \mathrm{VaR}_{\alpha_j} (\x) \ \leq \ 
    \mathrm{CVaR}_{\alpha_j}(\x) \ \leq \ 
    M_j(\x).
\end{equation*}
Thus, for all $\x \in \mathcal{X}$, $\lim_{\alpha_j \to 1} \mathrm{CVaR}_{\alpha_j}(\x) = M_j(\x)$ and we can set $\mathrm{CVaR}_{\alpha_j = 1}(\x) = M_j(\x) $ which ensures continuity of the $\mathrm{CVaR}_{\alpha_j}$ function with respect to $\alpha_j$ for $\alpha_j \in (0, 1]$. Now,
\begin{equation*}
    \mathrm{CVaR}_{\alpha_j = 1}^\beta(\x) \ = \ 
    \mathrm{VaR}_{\alpha_j = 1}^\beta(\x) \ = \ 
    \inf \{t \; | \; \mathbb{P}(C_j(\x + \beta_1 \u, \boldsymbol{\xi}) - \beta_2 v_j \leq t) = 1 \},
\end{equation*}
where the probability measure is taken with respect to $\boldsymbol{\xi}$, $\u \sim \mathcal{TN}(\mathbf{0}, \mathbf{I}, \frac{\b_\ell-\x}{\beta_1}, \frac{\b_u -\x}{\beta_1})$ and \linebreak $\v \sim \mathcal{TN}(\mathbf{0}, \mathbf{I}, \frac{-\t_{\max}}{\beta_2}, \frac{\t_{\max}}{\beta_2}) $. It follows that 
\begin{equation*}
    \mathrm{VaR}_{\alpha_j = 1}^\beta(\x) \ =  \ 
    \displaystyle  \sup_{ \substack{\boldsymbol{\xi} \in \Xi, \;  \u \in [\frac{\b_\ell-\x}{\beta_1}, \frac{\b_u -\x}{\beta_1}], \\ 
    v_j \in [(\frac{-\t_{\max}}{\beta_2})_j,(\frac{\t_{\max}}{\beta_2})_j]}} C_j(\x + \beta_1 \u, \boldsymbol{\xi}) - \beta_2 v_j \ = \ 
    \sup_{\x \in \mathcal{X}} M_j(\x) + (\t_{\max})_j,
\end{equation*}
where the $\sup$ is understood as the essential supremum of the function.  
Thus, for any $\x \in \mathcal{X}$, \linebreak $\mathrm{CVaR}_{\alpha_j = 1}(\x) <  \mathrm{CVaR}_{\alpha_j = 1}^\beta(\x)$. Therefore, by continuity of $\mathrm{CVaR}_{\alpha_j}$ with respect to $\alpha_j$, there exists  $\Bar{\alpha}_j \in (0, 1]$ such that for all $\alpha_j \geq \Bar{\alpha}_j$
    \begin{equation*}
        \mathrm{CVaR}_{\alpha_j}(\x) \ \leq \ 
        \mathrm{CVaR}_{\alpha_j}^\beta(\x) \ \leq \ 
        \mathrm{CVaR}_{\alpha_j}(\x) +  \frac{L_2 \beta_1 \sqrt{n} + \beta_2}{1 - \alpha_j},
    \end{equation*}
where the second inequality comes from the first part of the theorem. 
\end{proof}

\Cref{th_quali_approx}.2 shows that the difference in the values of objective function of Problem (\ref{problem_cvar_smooth}) and Problem (\ref{problem_cvar}) is bounded by a constant that depends on the values of $\alpha_j$, $\beta_1$, and $\beta_2$.~\Cref{th_quali_approx}.3 demonstrates that, with additional mild conditions, if $\alpha_j$ is chosen sufficiently close to $1$, the solution obtained in Problem (\ref{problem_cvar_smooth}) is feasible for Problem (\ref{problem_cvar}).  Therefore, the solution of Problem (\ref{problem_cvar_smooth}) may be feasible for Problem (\ref{problem_cvar}) and its value can be arbitrarily close to that of Problem (\ref{problem_cvar}) with sufficiently small values of $\beta_1$ and $\beta_2$. However, it is important to note that in practice, if $\beta_1$ and $\beta_2$ are chosen too small, the difference between the empirical values of the function will also be too small to represent the function differential~\cite{chen2019zo}.

To solve Problem (\ref{problem_cvar_smooth}) and to avoid the use of inner loops, which are computationally intractable, a Lagrangian relaxation is employed. This approach leads to the following unconstrained problem.

\begin{equation}
\label{problem_lagrangian}
    \max_{0 \leq \boldsymbol{\lambda} \in \R^m} \min_{(\x, \t) \in \mathcal{X}\times \R^{m+1}} L^\beta(\x, \t, \boldsymbol{\lambda}) :=  V_{\alpha_0}^\beta(\x, t_0) + \sum_{j = 1}^m \lambda_j V_{\alpha_j}^\beta(\x, t_j)  ,
\end{equation}
where $\t = (t_0, \dots, t_m) \in \R^{m+1}$. The next section describes a method allowing convergence to a saddle point of the Problem (\ref{problem_lagrangian}) whose the definition is recalled here.
\begin{Def}[Saddle point]
A saddle point of $L(\x, \t, \boldsymbol{\lambda})$ is a point $(\x^*, \t^*, \boldsymbol{\lambda}^*)$ such that for some $r > 0$, $\forall (\x, \t) \in \mathcal{X} \times \R^{m+1} \bigcap \mathcal{B}_{(\x^*, \t^*)}(r)$ and for all $\boldsymbol{\lambda} \geq 0$, we have
\begin{equation*}
    L(\x, \t, \boldsymbol{\lambda}^*) \geq L(\x^*, \t^*, \boldsymbol{\lambda}^*) \geq L(\x^*, \t^*, \boldsymbol{\lambda}),
\end{equation*}
where $\mathcal{B}_{(\x^*, \t^*)}(r)$ is a hyper-dimensional ball centred at $(\x^*, \t^*)$ with radius $r > 0$.
\end{Def}

\section{A Risk Averse Multi-timescale Stochastic Approximation Algorithm}
\label{sec_algo}
\Cref{sec-multi} presents the multi-timescale stochastic approximation methods,  and~\Cref{sec-algo} describes the complete algorithm used to solve the Problem (\ref{problem_lagrangian}).
\subsection{Multi-timescale Stochastic approximation methods} \label{sec-multi}
Multi-timescale is used to address the second difficulty raised at the end of Section 2, i.e., to avoid using nested loops to estimate a quantile of the level $\alpha$ and to compute the probabilistic constraints. Multi-timescale  stochastic approximation~\cite{bhatnagar2013stochastic, borkar_stochastic_2008} is a method that utilizes updates with different step-size schedules. Multi-timescale algorithms are useful when, between two successive updates of the algorithm, an inner-loop procedure must be performed recursively until it converges. Employing a multi-timescale algorithm allows both updates (for the inner and outer loops) to run together and converge to the desired point. In conditional value-at-risk ($\cvar$) optimization, this is typically the case for updating the additional variable $\t$ that could have been updated in an inner loop procedure. For example, the work~\cite{chow2017risk, prashanth2014policy} use a multi-timescale algorithm to update the additional variable. Other cases where multi-timescale can be applied include aggregating information about the gradient through an exponential moving average and updating the Lagrangian multipliers in the case of a Lagrangian relaxation. For more details on multi-timescale stochastic approximation, readers may refer to~\cite[Chapter 6]{borkar_stochastic_2008} or  ~\cite[Section 3.3]{bhatnagar2013stochastic}.

In this work, four different timescales are used. The four different step sizes $s_1^k, s_2^k$, $s_3^k$ and $s_4^k$ are chosen so that ~\Cref{assum2} holds.
\begin{assumption}
\label{assum2}
For $ k \geq 0$, the step sizes sequences $s_1^k, s_2^k$, $s_3^k$ and $s_4^k$ are strictly positive and satisfy the requirements:
\begin{align*}
    &\sum s_1^{k} =   \sum s_2^{k} =  \sum s_3^{k} = \sum s_4^{k} = +\infty, \\
    &\sum \left( (s_1^{k})^2 + (s_2^{k})^2 + (s_3^{k})^2 + (s_4^k)^2\right) < \infty,\\
    &\lim_{k \to \infty} \frac{s_1^{k}}{s_2^{k}} = \lim_{k \to \infty} \frac{s_2^{k}}{s_3^{k}} = \lim_{k \to \infty} \frac{s_3^{k}}{s_4^{k}} =  0.
\end{align*}
\end{assumption}
These four step sizes differ by their speed to reach the infinity. In fact, under the previous assumption, there exists an integer $k_0$ such that, for every $K \geq k_0$, the partial sums satisfy 
\begin{equation*}
    \sum_{k= 0}^K s_1^k < \sum_{k= 0}^K s_2^k
\end{equation*}
\begin{algorithm}[ht!]
\footnotesize
	\caption{Risk Averse Multi-timescale Stochastic Approximation (RAMSA) algorithm} 
	\begin{algorithmic}[1]
        \State{\textbf{Input:} $\x^{0}$, $\mathcal{X},\mathcal{T}, \mathcal{L}, K^{\max}$. }
        \State{Set $k = 0$ be an iteration counter}
        \State{Define stepsize sequences $(s_1^{k}), (s_2^{k}), (s_3^k) $ and $(s_4^{k})$ having the following form :}
        \begin{equation*}
            s_i^k = \frac{s_i^0}{(k+1)^{\tau_i}}, \, \forall i \in \{1, 2, 3, 4\}
        \end{equation*}
        \State{where the exponential decays $\tau_i, i = 1, \hdots, 4$ are chosen such that the~\Cref{assum2} are satisfied.}
        \State{Set $\M^0 = \Tilde{\g}^0, \V^0 = (\M^0)^2$ and $\t^0 = 0$}
        \While{ $k \leq K^{\max}$}
        \State{Draw samples $\u^k \sim \mathcal{TN}(\mathbf{0}, \mathbf{I}, \frac{\b_\ell - \x^k}{\beta_1}, \frac{\b_u- \x^k}{\beta_1})$ and $\v^k  \sim \mathcal{TN}(\mathbf{0}, \mathbf{I}, \frac{-\t_{\max}- \t^k}{\beta_2}, \frac{\t_{\max}- \t^k}{\beta_2})$.}
        \State{Recall that an unbiased output of the Lagrangian is given by:}
        \begin{equation}
            \Tilde{L}(\x, \t, \boldsymbol{\lambda, \boldsymbol{\xi}})= \Tilde{V}_{\alpha_0}(\x , t_0 , \boldsymbol{\xi}) + \sum_{j = 1}^m \lambda_j \Tilde{V}_{\alpha_j}(\x , t_j , \boldsymbol{\xi})
        \end{equation}
        \State{where $\Tilde{V}_{\alpha_j}(\x, t_j, \boldsymbol{\xi})  = t_j + \frac{1}{1-\alpha_j}\left(C_j(\x, \boldsymbol{\xi}) - t_j\right)^+ $.}
        \State{Calculate the gradient estimate $\Tilde{\g} := (\Tilde{\g}_\x, \Tilde{\g}_\t, \Tilde{\g}_{\boldsymbol{\lambda}}) \in \R^n \times \R^{m+1}\times \R^m$ with respect to $\x$, $\t$ and $\boldsymbol{\lambda}$ with:}
        \begin{align}
        \label{approx_grad_lagrangian}
        \begin{split}
                        &\Tilde{\g}_\x^k = \frac{\left( \Tilde{L}(\x^k + \beta_1 \u^k, \t^k + \beta_2 \v^k, \boldsymbol{\lambda}^k, \boldsymbol{\xi}_1^k) - \Tilde{L}(\x^k, \t^k, \boldsymbol{\lambda}^k, \boldsymbol{\xi}_2^k) \right) (\u^k - \mu_1^k) }{\beta_1}, \\
                        &\Tilde{\g}_\t^k = \frac{\left( \Tilde{L}(\x^k + \beta_1 \u^k, \t^k + \beta_2 \v^k, \boldsymbol{\lambda}^k, \boldsymbol{\xi}_1^k) - \Tilde{L}(\x^k, \t^k, \boldsymbol{\lambda}^k, \boldsymbol{\xi}_2^k) \right) (\v^k - \mu_2^k)}{\beta_2},\\
                        &\Tilde{\g}_{\lambda_j}^k = \Tilde{V}_{\alpha_j}(\x^k, t_j^k, \boldsymbol{\xi}_1^k) \; \forall j \in [1, m].
        \end{split}
        \end{align}
        \State{Update the long term gradient estimators:}
        \begin{align}
            \label{update_grad}
            \begin{split}
            &\M^{k+1} = s_4^{k} \Tilde{\g}^k + (1-s_4^{k})\M^{k}\\
            &\V^{k+1} = s_4^{k} ( \Tilde{\g}^k)^2 + (1-s_4^{k})\V^{k}
            \end{split}
        \end{align}
        \State{Update the current iterates $\x^k$, $\t^k$ and $\boldsymbol{\lambda}^k$;}
            \begin{align}
            \label{update_t}
            &\t^{k +1} = \Pi_{\mathcal{T}} \left[ \t^{k} - s_3^k \frac{\M_\t^{k+1}}{\sqrt{\V_\t^{ k+1}} + \epsilon} \right]  \\
            \label{update_x}
            &\x^{ k+1} = \Pi_{\mathcal{X}} \left[x^{k} - s_2^k \frac{\M_\x^{ k+1}}{\sqrt{\V_\x^{k+1}} + \epsilon} \right] \\
            \label{update_lambda}
            & \boldsymbol{\lambda}^{ k+1} = \Pi_{\mathcal{L}} \left[\boldsymbol{\lambda}^{k} + s_1^k \frac{\M_{\boldsymbol{\lambda}}^{ k+1}}{\sqrt{\V_{\boldsymbol{\lambda}}^{k+1}} + \epsilon} \right] 
        \end{align}
        \State{$k \leftarrow k +1$}
        \EndWhile
    \State{Return $\x^{k}$}
	\end{algorithmic} 
	\label{algo1}
\end{algorithm}
and the gap between the above two summations increases with $K$. Thus, the time scale associated with $s_2$ is said to be faster than the time scale associated with $s_1$. In this work, the fastest timescale is used to aggregate information about the gradient, the first intermediate timescale is used to update the additional variable $t$ and the second intermediate timescale is used to update the design vector $\x$, and the slowest timescale is used to update the Lagrangian multipliers $\boldsymbol{\lambda}$.

\subsection{The RAMSA algorithm} \label{sec-algo}
 
Algorithm \ref{algo1} summarizes the different updates. Note that when the square $(\cdot)^2$, the square root $\sqrt{\cdot}$ or division $\frac{\cdot}{\cdot}$ operators are applied to a vector, it is elementwise. Further remarks about algorithm \ref{algo1} are outlined: 
\begin{itemize}
    \item The updates (\ref{update_grad}) are the updates used to aggregate information about the gradient and are computed from the unbiased estimator defined in~\Cref{estim_grad_noised}. It will be shown later in the convergence proof that in fact $||\M^{k} - \nabla L(\x^{k}, \t^{k}, \boldsymbol{\lambda}) || \to 0$ and $||\V^{k} - (\nabla L(\x^{k}, \t^{k}, \boldsymbol{\lambda}))^2 || \to 0$ almost surely when $k \to \infty$. The $\M^k$ iterates can be thought of as an exponential moving average of the gradient estimators and aim to aggregate information about the direction of the gradient. The $\V^k$ iterates aim to avoid exploding gradient updates and aggregate information about the magnitude of the gradient.
    \item The update of the variable $\t$ is done in the update (\ref{update_t}). The interest of updating $\t$ with a faster timescale than those of $\x$ is that $\x$ will be quasi static compared to $\t$. Thus, for a given $\x$, the updates of $\t$ will appear to have converged to a point $\t^{*}(\x)$, where $\t^{*}$ is an estimate of the $\var$ at the point $\x$ of the objective and constraint functions. 
    \item A projection is employed in the updates  of the variables $\x$, $\t$ and $\boldsymbol{\lambda}$. This projection is required in the case of $\x$ because the space of the design variables is bounded. For $\t$ and $\boldsymbol{\lambda}$, the projection is required for convergence analysis. Since the bounds on $\t$ and $\boldsymbol{\lambda}$ can be arbitrarily large, this is not a problem in practice. In the algorithm, the sets $\mathcal{X}$, $\mathcal{T}$, and $\mathcal{L}$ are all hyperrectangles, i.e., sets of type $[\b_\ell, \b_u] \subset \R^d$ where $d$ is a given dimension. Furthermore, the projection operator $\Pi_{\mathcal{X}}(\x)$ is defined as $\Pi_{\mathcal{X}}(\x) = (\Pi_1(x_1), \dots,\Pi_d(x_d))$, where the individual projection operators $\Pi_j : \R \to \R$ are defined by $\Pi_j(x_j) = \min( (\b_u)_j, \max( (\b_\ell)_j, x_j))$ for all $j \in [1, d]$. The projection operators for the variables $\t$ and $\boldsymbol{\lambda}$ are defined in the same way.
\end{itemize}

\section{Convergence analysis}
\label{sec_analysis}

The convergence of the RAMSA algorithm is stated in the following theorem. 
\begin{Th}
\label{main_th}
    Under~\Cref{assum1}.3 and~\Cref{assum2}, let further assume that the problem given in~\Cref{problem_cvar_smooth} is strictly feasible and there exists $K \in \N$  such that $\x^K$ and $\boldsymbol{\lambda}^K$ are in the domain of attraction of $\x^*$ and $\boldsymbol{\lambda}^*$ with $\boldsymbol{\lambda}^* \in \mathcal{L}^\circ$ respectively. Then, the iterates $(\x^k, \t^k, \boldsymbol{\lambda}^k)$, produced by the RAMSA algorithm, converge almost surely to a saddle point of the Lagrangian function $L^\beta$ and $(\x^*, \t^*)$ is a locally optimal solution for the smooth CVaR-constrained problem given in~\Cref{problem_cvar_smooth}.
\end{Th}

While the technical details of the proof of this theorem are given in~\Cref{appendix_proof}, a high-level overview of the proof steps is given below.
\begin{itemize}
    \item First, for each timescale, a discrete stochastic approximation analysis is used to prove the almost sure convergence of the iterates $(\M^k, \V^k, \x^k, \t^k, \boldsymbol{\lambda}^k )$ to a stationary point $(\M^*, \V^*, \x^*, \t^*, \boldsymbol{\lambda}^*)$ of the corresponding continuous-time system. 
    \item Then, to show that the continuous-time system is locally asymptotically stable at the stationary point, a Lyapunov analysis is performed.
    \item Finally, considering the iterates $(\x^k, \t^k, \boldsymbol{\lambda}^k )$, the Lyapunov function used in the above analysis is the Lagrangian function $L(\x, \t, \boldsymbol{\lambda})$. Therefore, the stationary point $(\x^*, \t^*, \boldsymbol{\lambda}^*)$ is a saddle point. Thus, by the saddle point theorem, we deduce that $\x^*$ is a locally optimal solution to the smooth  $\cvar$-constrained blackbox optimization problem given in~\Cref{problem_cvar_smooth}.
\end{itemize}
This convergence proof procedure is standard for multi-timescale stochastic approximation algorithms, see~\cite[chapter 10]{bhatnagar2013stochastic},~\cite[chapter 6]{borkar_stochastic_2008} or~\cite{bhatnagar2013stochastic, chow2017risk}, for further references. Note that this procedure must be done for each timescale, requiring four similar proof steps. This is due to the different speeds of the timescales. Here, the updates $(\M^k, \V^k)$ converge on a faster timescale than $\t^k$, which converges on a faster timescale than $\x^k$, while $\boldsymbol{\lambda}^k$ converges on the slowest timescale. The idea of multi-timescale convergence analysis is then to assume that, given a timescale, the updates made on faster timescales are quasi-equilibrated, i.e. have already converged to an equilibrium point. The updates made on slower timescales are quasi-static, i.e. fixed with respect to the given timescale. Therefore, the convergence analysis of the updates of the given timescale is done by considering all other updates as fixed. To illustrate the mathematical meaning of this assumption, consider two updates $\x^k, \x_2^k \in \mathcal{X}_1 \times \mathcal{X}_2$ such that
\begin{align}
    \x_1^{k+1} &= \x_1^k + s_1^k \left( f_1(\x_1^k, \x_2^k) + \delta_1^{k+1} \right),\\
    \x_2^{k+1} &= \x_2^k + s_2^k \left( f_2(\x_1^k, \x_2^k) + \delta_2^{k+1} \right),
\end{align}
where $f_1$ and $f_2$ are Lipschitz continuous function and $\delta_1$, $\delta_2$ are square integrable martingale difference sequence with respect to the $\sigma$-field $\sigma( \x_1^i, \x_2^i, \delta_1^i; i \leq k)$ and $\sigma( \x_1^i, \x_2^i, \delta_2^i; i \leq k)$. If $s_1^k$ and $s_2^k$ are non-summable and square summable step sizes with $s_2^k$ which is a faster timescale than $s_1^k$, i.e., $s_1^k = o(s_2^k)$. Then, the previous recursion may be rewritten as follows 
\begin{align}
    \x_1^{k+1} &= \x_1^k + s_2^k \left( \frac{s_1^k}{s_2^k} \left( f_1(\x_1^k, \x_2^k) + \delta_1^{k+1} \right) \right),\\
    \x_2^{k+1} &= \x_2^k + s_2^k \left( f_2(\x_1^k, \x_2^k) + \delta_2^{k+1} \right).
\end{align}
As $s_1^k = o(s_2^k)$, this recursion may be seen as a noisy discretization of the  ODEs $\Dot{\x_1} = 0 $ and $\Dot{\x_2} = f_2(\x_1, \x_2)$. Since $\Dot{\x_1} = 0 $, $\x_1$ is a constant and the second ODE may be replace with $\Dot{\x_2} = f_2(\x_1^0, \x_2)$, where $\x_1^0$ is a constant. Finally it can be proved~\cite[Chapter 6, Theorem 2]{borkar_stochastic_2008} that $(\x_1^k, \x_2^k)$ converge  $(\x_1^*, \mu(\x_1^*))$, where $\mu$ is a Lipschitz continuous function, $\mu(\x_1^*)$ is a locally stable equilibrium of the ODE $\Dot{\x_2} = f_2(\x_1^*, \x_2)$ and $\x_1^*$ is a locally stable equilibrium of the ODE $\Dot{\x_1} = f_1(\x_1, \mu(\x_1))$.

In~\Cref{main_th}, it is proved that the iterations converge to a locally optimal solution of the problem given in~\Cref{problem_cvar_smooth}. It is possible to obtain a result for the original CVaR-constrained problem given in~\Cref{problem_cvar} by utilizing ~\Cref{th_quali_approx}. This is the subject of the following corollary. 
\begin{Cor}
    Under the same assumptions as~\Cref{main_th}, it follows that there exists a threshold $\Bar{\alpha} \in (0,1]$ such that if $ \alpha_j \geq \alpha$ for all $j \in [0, m]$, then the iterates $\x^k$  converge almost surely to a feasible solution $\x^*$ of Problem (\ref{problem_cvar}) whose 
    the objective function value is within $\frac{L_2 \beta_1 \sqrt{n} + \beta_2}{1 - \max_{j \in [1,m]} \alpha_j}$ of that of a local solution of  Problem (\ref{problem_cvar}).
\end{Cor}
\begin{proof}
    The proof is straightforward, considering the result of Theorems \ref{main_th} and \ref{th_quali_approx}.
\end{proof}
This corollary is particularly interesting because it ensures the almost sure convergence of~\Cref{algo1} to  a feasible point of the $\cvar$-constrained problem whose objective function value is arbitrarily close to that of a local solution. To the best of our knowledge, this result is the first of its kind in the area of derivative-free RBDO with unknown uncertainty distribution. 
\section{Computational implementations and numerical experiments}
\label{sec_exp}
This section is divided into five parts: details of the numerical implementation are given in~\Cref{sec-num}.~\Cref{sec-setup} describes the setup of the experiments.~\Cref{sec-comp-a} presents the experiments aimed at finding relations between the hyperparameters and the problems to be solved. Finally,~\Cref{sec-comp-b} exhibits the results obtained using the truncated Gaussian gradient estimator instead of its classical counterpart, while~\Cref{sec-comp-c} shows the results when the problem is subject to mixed aleatory/epistemic uncertainties.
\subsection{Computational implementation} \label{sec-num}
In this section, practical details of the implementation of~\Cref{algo1} are given. They aim to reduce the number of hyperparameters required by the algorithm and improved its practical efficiency.

The first difficulty the algorithm faces is when the bounds of the decision variables differ in magnitude. A first approach is then to adjust the initial step sizes according to each coordinate. However, this increases the number of hyperparameter values to be set. Another approach, which requires only one step size for all coordinates $j \in [0, m]$, is to map the initial hyperrectangle to the hypercube $[0, 1]^n$. The output of the blackbox $C_j : \mathcal{X} \to \R$ is simply replaced by $C_j^1 : [0, 1]^n \to \R$, where 
\begin{equation*}
    C_j^1(\x, \boldsymbol{\xi}) = C_j( \b_\ell + (\b_u - \b_\ell) \x, \boldsymbol{\xi}).
\end{equation*}

The algorithm encounters a second difficulty related to the Lagrangian relaxation, where the values of the objective function and constraints are added together. When constraint magnitudes differ, the algorithm is biased towards the larger ones. To mitigate this bias, a solution consists of choosing different step sizes for updating $\boldsymbol{\lambda}$ but that increases the number of hyperparameters. Alternatively, a transformation may be applied to normalize the values, allowing the use of a single step size. In this method, the $\arctan(\cdot)$ function is employed to map the blackbox output values to the range of $[-\frac{\pi}{2}, \frac{\pi}{2}]$. However, there is an issue when the bounds of the $\arctan$ function are approached because the gradient estimator is computed from the difference between the values returned by the $\arctan$ function. If this difference is too small, especially in the presence of noisy blackbox outputs, the quality of the gradient estimator decreases. To address this issue, the cubic root function is applied beforehand to increase the difference between these values. That leads to the following transformation
\begin{equation*}
    C_j^2(\x, \boldsymbol{\xi}) = \arctan \left( \sqrt[3]{C_j(\x, \boldsymbol{\xi}) } \right), \forall j \in [0, m].
\end{equation*}
 In the rest of the paper, we refer to  $ \Tilde{C}_j : [0, 1]^n \to [-\frac{\pi}{2}, \frac{\pi}{2}], \, \forall j \in [0,m]$, the map corresponding to the two previous transformations applied to the outputs of the blackbox.

Finally, in practical applications, it appears that initiating the process directly at the intended reliability level can be counterproductive~\cite{zhu_simulation_2018}. To overcome this difficulty, the values of $\alpha_j, \forall j \in [0, m]$ are initially set to $0$. Then, these values are gradually increased until the desired reliability levels are reached. This is done by inserting reliability level setting 
\begin{equation*}
        \alpha_j^{k+1} = \alpha_j^* +  \gamma \left( \alpha_j^k - \alpha_j^* \right) 
\end{equation*}
 for every index $j \in [0, m]$ in between lines $12$ and $13$ of~\Cref{algo1}. Here, $\alpha_j^*$ are the desired reliability levels and $\gamma  \in [0, 1)$ is a fixed threshold.

\subsection{Numerical experiments} \label{sec-setup}
Before proceeding to the numerical experiments, this section describes the test problems chosen, the way the experiments are performed, and the objectives of the different experiments.

First, four analytical test problems, each with a known practical optimum, are chosen from existing literature. These problems include a Steel Column Design (SCD) problem~\cite{yang2020enriched}, a Welded Beam Design (WBD) problem~\cite{yang2020enriched}, a Vehicle Side Impact (VSI) problem~\cite{yang2020enriched}, and a Speed Reducer Design (SRD) problem~\cite{chen2013optimal}. These problems are decribed in~\Cref{prob-def}, and further information regarding their physical interpretations can be found in the associated references. Except in the last subsection, the goal is to solve the following standard RBDO problem
\begin{align}
\label{prob_stand}
\begin{split}
\begin{array}{cl}
    \displaystyle \min_{\x \in [0, 1]^n }   &\mathbb{E}_{\boldsymbol{\xi}}[\Tilde{C}_0(\x, \boldsymbol{\xi} )]\\
     \mbox{s.t.}  &\mathbb{P}(\Tilde{C}_j(\x, \boldsymbol{\xi}) \leq 0) \geq 0.99, \; \forall j \in [1,m].
\end{array}
\end{split}
\end{align}
It is important to note that Problem (\ref{prob_stand}), unlike the classical FORM-based problem, incorporates uncertainties not only in the constraints but also in the objective function. Moreover, despite the analytical expressions of the problems are available and the  uncertainty distributions are known, the RAMSA algorithm operates without utilizing these information. As outlined in Sections \ref{sec_prob_stat} and \ref{sec_smoo_prob} it solves formally a smooth Lagrangian relaxation of Problem (\ref{prob_stand}).

In order to make comparisons, it is essential to devise a strategy for evaluating the quality of solutions generated by the RAMSA algorithm. As both the problem and the algorithm are subject to uncertainties, multiple runs of the RAMSA algorithm are necessary, and the values of the proposed solutions need to be estimated using Monte Carlo simulations. In this work, a trial consists of running the algorithm $100$ times with the same set of hyperparameters values. For each run, a maximum budget of $5000$ function evaluations is allocated. At the end of these $100$ runs, the final solution points are recorded. For each solution point, the mean of the objective function and the probabilisty to satisfy the constraints are estimated through $10 000$ Monte Carlo simulations. A run is deemed successful if all constraints are satisfied with a probability greater than $0.99$. Moreover, the mean solution point over the 100 runs, denoted as $\Bar{\x}^*$, is calculated as well as its standard deviation. That allows to check that the RAMSA algorithm consistently converges to the same neighborhood of an optimal point. To further validate the results, this point is also compared with the solution obtained by the SORA algorithm in~\cite{chen2013optimal, yang2020enriched}. Note that the aim is not to directly compare the RAMSA and SORA algorithms since the SORA algorithm takes advantage of the analytical expressions of the problems and knowledge of uncertainty distributions. When a trial is consistent for a set of hyperparameter, the set and the trial are said to be satisfactory. 

Now, the objectives of the upcoming experimental sections are threefold. First, despite the transformations introduced in the previous section, there are still some hyperparameters that need to be configured.~\Cref{sec-comp-a} provides guidelines on how to set these hyperparameters. Second, a critical aspect is the selection of the kernel density used to estimate gradients during the optimization process. In~\Cref{sec-comp-b}, a comparison is made between the classical Gaussian gradient estimator and the truncated Gaussian gradient estimator introduced in~\Cref{sec_smoo_prob}.1. Third, the VSI problem is described slightly differently in~\cite{yang2020enriched}, allowing the means of the uncertainty variables $\boldsymbol{\xi}_8$ and $\boldsymbol{\xi}_9$ to take two values: $0.192$ and $0.345$. This is an opportunity to employ the RAMSA algorithm for solving the VSI problem under mixed aleatory/epistemic uncertainties. In fact, the uncertainty in distribution parameters can be regarded as a source of epistemic uncertainty~\cite{nannapaneni2016reliability}. Detailed descriptions of the conducted experiments are presented in~\Cref{sec-comp-c}.
   
\subsection{Hyperparameters setting rules} \label{sec-comp-a}
The RAMSA algorithmn involves four types of hyperparameters: the exponential decays of the step sizes $\tau \in (\frac{1}{2}, 1)^4$, the threshold for the adaptive reliability level $\gamma$, the initial step sizes $s^0 \in \R_+^4$, and the smoothing parameters $\beta \in \R_{+*}^2$. Two strategies can be employed to determine the values of these hyperparameters.

On the one hand, theoretical considerations are employed to set some hyperparameter values. This approach is employed to set the values of the exponential decays. These values must satisfy ~\Cref{assum2} to ensure the convergence of the algorithm. Moreover, they must be distinct enough to achieve the desired multi-timescale effect, but also not too different, otherwise, either the fastest timescale is too fast (leading to increased noise) or the slowest timescale is overly slow (impeding the convergence process)~\cite[Chapter 6]{bhatnagar2013stochastic}. Thus, the decays are arbitrarily set to $\tau = (0.8, 0.7, 0.6, 0.501)$.  The threshold for the adaptive reliability level $\gamma$ can be determined similarly. This hyperparameter depends only on the value of $K^{\max}$ because for $j \in [0, m]$, it follows that $\forall k \in \N, \, \alpha_j^{k+1}  = \alpha_j^*(1-\gamma^k)$. Thus $\gamma$ can be chosen such that $\alpha_j^{K^{\max}} \approx \alpha_j^*$. However, if $\alpha_j^*$ is chosen close to $1$, the problem given in~\Cref{problem_cvar_smooth} is particularly conservative for Problem (\ref{problem_cvar}), as shown in~\Cref{th_quali_approx}, and even more so for Problem (\ref{prob_stand}). Therefore, to avoid overly conservative results, $\gamma$ is chosen to be equal to $ 1 - \frac{5}{2 K^{\max}}$ so that $\alpha_j^{K_{\max}} \approx 0.9$ provided that $K_{\max} = 2500$ and $\alpha^* = 0.99$.

On the other hand, there are some hyperparameters values that cannot be determined theoretically. In this case, they have to be computed experimentally. This is achieved through a two-step strategy. The set of test problems is divided into two groups: the experimental test problems and the validation test problems. In the first step, for each experimental problem, a set of hyperparameters, that gives satisfactory results on this test problem, is identified. By analyzing the results obtained on the different problems and the associated hyperparameter values, a distinction may be deduced between the hyperparameters which are problem-dependent and which are not. For problem-dependent hyperparameters, we try to establish correlations between the hyperparameter values and relevant problem-related quantities. Examples of such quantities include the objective function value, the gradient norm, or its variance at the starting point. Then, the validation step is undertaken to check the rules derived from the experimental step. During this phase, the rules are applied to the validation test problems to determine the  hyperparameter values of the RAMSA algorithm. If the results obtained with this set of hyperparameters are satisfactory, the rules are  deemed effective.

In this study, the two-step strategy is applied as follows. The experimental test problems selected are the VCD, WBD, and VSI problems. Trials of~\Cref{algo1} are conducted with different sets of hyperparameter values and the classical Gaussian gradient estimator~\cite[Equation (26)]{nesterov2017random}. For the sake of brevity, only one set of satisfactory hyperparameters and its associated results are presented for each problem. The values of this set are listed in~\Cref{table_hyper},  while in~\Cref{table_res1} the associated average results of the trials are presented. Detailed results from the $100$ runs of the trials are provided in~\Cref{detailed-res} in the form of boxplots.
\begin{table}[ht!]
    \centering
    \begin{tabular}{|C{2cm}||C{2cm}|C{2cm}|C{2cm}|C{2cm}|C{2cm}|C{2cm}|}
    \hline
           Problem& 
          $\beta_1$ & $\beta_2$  & $s_1^0$ & $s_2^0$ & $s_3^0$ & $s_4^0$  \\
    \hline
    \hline
    SCD & $0.05$ & $0.0001$  & $0.01$ & $0.05$ &  $0.001$& $0.2$ \\
    \hline
    WBD & $0.002$ & $0.0001$  & $0.01$ & $0.001$ &  $0.001$& $0.4$ \\
    \hline
    VSI & $0.1$ & $0.0001$  & $0.01$ & $0.5$ &  $0.001$& $0.5$ \\
    \hline
    \end{tabular}
    \caption{Satisfactory set of hyperparameter values found for each problem}
    \label{table_hyper}
\end{table}

\begin{table}[ht!]
    \centering
    \small
    \begin{tabular}{|C{2cm}||C{2cm}|C{4cm}|C{4cm}|C{2.cm}|C{1.5cm}|}
    \hline
           Problem/ Algo & Average of $\mathbb{E}[C(\x^*, \boldsymbol{\xi})]$
           &  Average of $ \mathbb{P}(C_j(\x^*, \boldsymbol{\xi}) \leq 0)$ & Average result point $\Bar{\x}^*$ (and standard deviation)& Number of successful runs & Function queries  \\
    \hline
    \hline
    SCD & $3967$ & $[0.9938]$ & $\underset{[\pm 4.4, \, \pm 0.25, \, \pm 3.8]}{[229.7, 15.03, 103.1]}$ &  $100$& $5000$ \\
    \hline
    SORA & $3989$ & $[0.9947]$ & $[258, 13.5, 100]$  & N/A & $216$ \\
    \hline \hline
    WBD & $2.53$ & $[1.0, 1.0, 0.9995, 1.0, 1.0]$ & $\underset{[\pm 0.01, \, \pm 0.29, \, \pm 0.29, \, \pm0.02]}{[6.36, 158, 211, 6.59]}$ &  $100$& $5000$ \\
    \hline
    SORA & $2.49$ & $[1.0, 1.0, 1.0, 1.0, 1.0]$ & $[5.92, 181, 211, 6.22]$  & N/A & $505$ \\
    \hline \hline
        \multirow{4}{*}{VSI } & \multirow{4}{*}{$28.38$} & \multirow{1}{*}{$[1.0, 1.0, 1.0, 1.0, 0.9993$} & \multirow{1}{*}{$\underset{[\pm 0.03, \, \pm 0.004, \, \pm 0.02, \, \pm 0.009]}{[0.88, 1.34, 0.51, 1.49,}$} &  \multirow{4}{*}{$95$}& \multirow{4}{*}{$5000$} \\
     & &\multirow{1}{*}{ $1.0, 1.0, 0.9925, 1.0,$} &  & &\\
    & & \multirow{1}{*}{$ 0.9996]$} &  \multirow{1}{*}{$\underset{\pm 0.07, \, \pm 0.01, \, \pm 0.08]}{1.29, 1.19, 0.45]}$} & &\\
    & & & & &\\ 
    \hline
     \multirow{3}{*}{SORA} &  \multirow{3}{*}{$29.55$} &  \multirow{1}{*}{$[1.0, 1.0, 1.0, 1.0, 0.9987,$} &  \multirow{1}{*}{$[0.78, 1.35, 0.69, 1.5,$}  &  \multirow{3}{*}{N/A} &  \multirow{3}{*}{$8054$} \\
     & &  \multirow{1}{*}{$1.0, 0.9987, 0.9983, 1.0,$} &  \multirow{1}{*}{$1.07, 1.2, 0.78]$} & &\\
     & &\multirow{1}{*}{$ 0.9993]$} & & &\\
    \hline 
    \end{tabular}
    \caption{Average result over $100$ runs obtained for each problem}
    \label{table_res1}
\end{table}

\Cref{table_res1} shows that the RAMSA algorithm achieves satisfactory results in all three problems. Interestingly, it appears to perform better on problems with higher dimensions and more constraints. This phenomenon can be attributed to the approximation of the gradient used in the RAMSA algorithm. This approximation estimates the gradient of the Lagrangian function with only two blackbox evaluations, regardless of the dimension or number of constraints. Upon analyzing~\Cref{table_hyper}, it seems that $\beta_2$, $s_1^0$, and $s_3^0$ are problem-independent. Moreover, the value of $s_4^0$ falls within a relatively narrow interval of $[0.1, 0.6]$. In contrast, the smoothing parameter $\beta_1$ and the initial step size $s_2^0$, both associated with the design vector $\x$, exhibit variations from one problem to another. This variability suggests the problem dependency of these hyperparameters.

The first claim to be proven experimentally is the following: an appropriate order of magnitude of $\beta_1$ is so that the variance of the gradient estimator at the starting point is minimal. A such value should reduce the variability during the initial stages of the optimization process and thus improve the convergence rate. To validate this assertion, the gradient is approximated by computing $N$ Lagrangian gradient estimators given in~\Cref{approx_grad_lagrangian}, at the point $(\x^0, \mathbf{0}, \mathbf{0})$. The gradient is approximated for only $6$ different values of $\beta_1$ 
to prevent excessive computations. The values chosen are $ [0. 001, 0.005, 0.01, 0.05, 0.1, 0.2]$. Then, the variance of the first $n$ components of the gradient (i.e., the components of $\Tilde{\g}_\x $) is computed, and the average of these variances is calculated for each value of $\beta_1$. The value of $\beta_1$ is finally chosen as the one leading to the smallest average variance. If the minimum is reached for two different values of $\beta_1$, the larger value is selected. The results for the three different problems are presented in~\Cref{table_beta1}. It is observed that, selecting $\beta_1$ to minimize the average variance and halving it, yields to similar results to those of~\Cref{table_hyper}.
\begin{table}[ht!]
    \caption{Average variance of $N$ gradient approximations for different values of the smoothing parameter~$\beta_1$}
    \centering
    \begin{tabular}{|C{4cm}||C{1.5cm}|C{1.5cm}|C{1.5cm}|C{1.5cm}|C{1.5cm}|C{1.5cm}|}
    \hline 
     Value of $\beta_1$&$0.001$  & $0.005$ & $0.01$ & $0.05$ & $0.1$ & $0.2$   \\
    \hline \hline
    Average variance for SCD problem  & $4.2$  & $0.16$ & $0.04$ & $0.004$ & $0.003$ & $0.94$ \\
    \hline 
    Average variance for WBD problem  & $2.1$  & $1.75$ & $1.78$ & $5.2$ & $15$ & $8.9$ \\
    \hline
    Average variance for VSI problem  & $0.67$  & $0.03$ & $0.008$ & $0.0018$ & $0.0016$ & $0.0016$ \\
    \hline
    \end{tabular}
    \label{table_beta1}
\end{table}

The second claim to be experimentally shown is that: there is a correlation between the norm of the stochastic gradient and the value of the initial step size $s_2^0$. Intuitively, that means that the smaller the gradient norm, the larger the initial step size should be, and vice versa. To validate this hypothesis, $N$ stochastic gradients with $\beta_1 = 0.1$ are computed, and the norm of their mean is calculated. The result, normalized by the square root of the dimension, is presented in the third line of ~\Cref{table_norm} for each problem. The second line displays the result obtained in~\Cref{table_hyper}, and the last line shows the corresponding correlation coefficients. Based on these results, it can be deduced that the correlation coefficient should be around $10^{-3}$.
\begin{table}[ht!]
\caption{Correlation between the norm of the gradient and the initial step size $s_2^0$}
    \centering
    \begin{tabular}{|C{5cm}||C{2.5cm}|C{2.5cm}|C{2.5cm}|}
    \hline
            & SCD  & WBD  & VSI  \\
    \hline \hline
    \vspace{2pt}
    Value of $s_2^0$ & $0.05$ & $0.001$ & $0.5$ \\
    \hline
    \vspace{2pt}
    Estimated value of
    $\frac{||\nabla_{\x} L^\beta(\x^0, \boldsymbol{\xi})||_2}{\sqrt{n}} $& $\approx 0.02$& $\approx 1.4$ & $\approx 0.01$ \\
    \hline
    $s_2^0 \times \frac{ ||\nabla_{\x} L^\beta(\x^0, \boldsymbol{\xi})||_2}{\sqrt{n}} $ & $\approx 0.001$ & $\approx 0.001$ & $\approx 0.005$ \\
    \hline
    \end{tabular}
    \label{table_norm}
\end{table}

In the conducted experiments, the value of $N$ is set to $10000$. It is worth noting that while this large sample size is suitable for these experiments, in a BBO context, such a number might be intractable due to its computational cost. However, the methodology employed here can be adapted to work with smaller sample sizes. The goal of this approach is to provide only an order of magnitude for the hyperparameter values. Thus, a reduced number of samples can be used in a BBO context. Additionally, it is worth mentioning that the calculated gradients used to estimate the value of $\beta_1$ can also be used to estimate the value of $s_2^0$, reducing the computational cost of the method.

To validate the experimental step, the claims previously stated are applied to compute the hyperparameter values for solving the SRD problem. For this problem, the minimum value of the average variance occurs for $\beta_1 = 0.1$, and the norm of the Lagrangian gradient (normalized by the dimension) is estimated to be $0.006$. These values are then utilized to set the values of $\beta_1 = 0.05$ and  $s_2^0 = 0.15$. The values of the others hyperparameters are set as in~\Cref{table_hyper} and $s_4^0 = 0.2$. The results obtained with this set of values are shown in~\Cref{table_res2}.
\begin{table}[ht!]
    \caption{Average result over $100$ runs for Speed Reducer design problem}
    \centering
    \small
    \begin{tabular}{|C{2cm}||C{2cm}|C{4cm}|C{4cm}|C{2.cm}|C{1.5cm}|}
    \hline
           Problem/ Algo & Average of $\mathbb{E}[C(\x^*, \boldsymbol{\xi})]$
           &  Average of $ \mathbb{P}(C_j(\x^*, \boldsymbol{\xi}) \leq 0)$ & Average result point $\Bar{\x}^*$ (and standard deviation)& Number of successful runs & Function queries  \\
    \hline \hline
        \multirow{4}{*}{SRD} & \multirow{4}{*}{$3148$} & \multirow{1}{*}{$[1.0, 1.0, 1.0, 1.0, 1.0$} & \multirow{1}{*}{$\underset{[\pm 0.0, \, \pm 0.0, \, \pm 0.06, \, \pm 0.04]}{[3.6, 0.7, 17.0, 7.41,}$} &  \multirow{4}{*}{$100$}& \multirow{4}{*}{$5000$} \\
     & &\multirow{1}{*}{ $1.0, 1.0, 0.9996, 1.0,$} &  & &\\
    & & \multirow{1}{*}{$1.0, 1.0]$} &  \multirow{1}{*}{$\underset{\pm 0.04, \, \pm 0.01, \, \pm 0.01]}{7.99, 3.51, 5.37]}$} & &\\
    & & & & &\\ 
    \hline
     \multirow{3}{*}{SORA} &  \multirow{3}{*}{$3038$} &  \multirow{1}{*}{$[1.0, 1.0, 1.0, 1.0, 0.9975,$} &  \multirow{1}{*}{$[3.57, 0.7, 17, 7.3,$}  &  \multirow{3}{*}{N/A} &  \multirow{3}{*}{$2486$} \\
     & &  \multirow{1}{*}{$0.9986, 1.0, 0.9986, 1.0,$} &  \multirow{1}{*}{$7.75, 3.36, 5.3]$} & &\\
     & &\multirow{1}{*}{$1.0, 0.9986]$} & & &\\
    \hline 
    \end{tabular}
    \label{table_res2}
\end{table}

Based on these results, it appears that the rules established for setting the hyperparameter values lead to satisfactory solutions. The consistency observed in the solution points, as indicated by the small standard deviations obtained, suggests that the algorithm consistently converges to the same vicinity. Moreover, this solution is relatively close to the optimal point found by the SORA algorithm. Note, however, that these rules do not guarantee to find the best possible set of hyperparameters. For example, by retaining all hyperparameter values but adjusting $\beta_1$ to $0.01$,  similar values of probabilistic constraints can be achieved, with an average objective function value of $3066$. 

In summary, the rules established in this section provide valuable insights into obtaining a satisfactory set of hyperparameter values for the RAMSA algorithm. However, they must be used with caution due to the limited number of problems used to derive them, especially for  the value of $\beta_1$. It is known~\cite{chen2019zo} that setting the appropriate $\beta_1$
value is a challenging task in practice. One potential approach to address this challenge is to dynamically decrease the value of $\beta_1$ during the optimization process, as done in ~\cite{audet2023sequential}. Nevertheless, this topic falls beyond the scope of the present paper and is not explored further here.

\subsection{Truncated Gaussian vs Gaussian gradient estimator}\label{sec-comp-b}

In this section, the focus is on investigating the behavior of the algorithm when the bound constraints are unrelaxable~\cite{digabel2015taxonomy}, meaning that the outputs of the blackbox are not meaningful for the optimization process. This situation can arise when the blackbox is not defined outside its bounds or due to physical phenomenon. In this section, the uncertainties specified in~\Cref{prob-def} are truncated, ensuring that $\x + \boldsymbol{\xi} \in \mathcal{X}$ for every realization of $\boldsymbol{\xi}$. Moreover, to solve the constrained problem, the algorithm is executed using the truncated Gaussian gradient estimator instead of the classical Gaussian gradient estimator utilized in the previous section. This modification guarantees that all the candidate points are evaluated inside the bound constraints $\mathcal{X}$.

To determine the hyperparameter values for the algorithm using the truncated Gaussian gradient estimator, the methodology introduced in the previous section is applied. The values of $\beta_1$ that minimize the variance of the truncated Gaussian estimator are found to be $0.2, 0.005, 0.2$ and $0.01 $for the SCD, WBD, VSI, and SRD problems, respectively. Consequently, the values of $\beta_1$ are set to $0.1, 0.0025, 0.1$ and $0.025$. Furthermore, the correlation coefficient between the norm of the approximate gradient and the initial step size $s_2^0$ 
is approximately $5 \times 10^{-4}$. Thus, the values of $s_2^0$ are set to $0.1, 0.0008, 0.6$ and $0.01 $ for the SCD, WBD, VSI, and SRD problems, respectively. Finally, the values of $s_4^0$ are set to $0.25, 0.4, 0.6$, and $0.2$.
The results of these experiments are presented in~\Cref{table_res3}, and the detailed results from the 100 runs are depicted in boxplots in~\Cref{detailed-res}. 
\begin{table}[ht!]
    \centering
    \small
        \caption{Best average result over $100$ runs obtained for each problem with truncated Gaussian gradient estimator with $15000$ function evaluations by run}
    \begin{tabular}{|C{2cm}||C{2cm}|C{4cm}|C{4cm}|C{2.5cm}|}
    \hline
           Problem & Average of $\mathbb{E}[C(\x^*, \boldsymbol{\xi})]$
           &  Average of $ \mathbb{P}(C_j(\x^*, \boldsymbol{\xi}) \leq 0)$ & Average result point $\Bar{\x}^*$ (and standard deviation)& Number of successful runs   \\
    \hline
    \hline
    SCD & $3957$ & $[0.9958]$ & $\underset{[\pm 10, \, \pm 0.6, \, \pm 6]}{[226, 15, 106]}$ &  $97$ \\
    \hline
    WBD & $2.53$ & $[0.999, 1.0, 1.0, 1.0, 1.0]$ & $\underset{[\pm 0.01, \, \pm 0.24, \, \pm 0.24, \, \pm0.01]}{[6.37, 158, 211, 6.59]}$ &  $99$ \\
    \hline
     \multirow{4}{*}{VSI} & \multirow{4}{*}{$28.97$} & \multirow{1}{*}{$[1.0, 1.0, 1.0, 1.0, 0.9998$} & \multirow{1}{*}{$\underset{[\pm 0.03, \, \pm 0.004, \, \pm 0.03, \, \pm 0.008]}{[1.00, 1.35, 0.54, 1.49,}$} &  \multirow{4}{*}{$91$} \\
    & &\multirow{1}{*}{ $1.0, 1.0, 0.9923, 1.0,$} &  & \\
    & & \multirow{1}{*}{$ 0.9977]$} &  \multirow{1}{*}{$\underset{\pm 0.1, \, \pm 0.008, \, \pm 0.1]}{1.21, 1.19, 0.48]}$} & \\
    & & & & \\ 
    \hline
     \multirow{4}{*}{SRD } & \multirow{4}{*}{$3093$} & \multirow{1}{*}{$[1.0, 1.0, 1.0, 1.0, 1.0$} & \multirow{1}{*}{$\underset{[\pm 0.002, \, \pm 0.0, \, \pm 0.03, \, \pm 0.004]}{[3.58, 0.7, 17.0, 7.30,}$} &  \multirow{4}{*}{$100$} \\
     & &\multirow{1}{*}{ $1.0, 1.0, 0.994, 1.0,$} &  & \\
    & & \multirow{1}{*}{$1.0, 0.999]$} &  \multirow{1}{*}{$\underset{\pm 0.003, \, \pm 0.002, \, \pm 0.0008]}{7.78, 3.42, 5.31]}$} & \\
    & & & & \\ 
    \hline
    \end{tabular}
    \label{table_res3}
\end{table}

In~\Cref{table_res3}, it is shown that utilizing the truncated Gaussian gradient approximation leads to satisfactory results. However, the algorithm convergence is significantly slower than with classical Gaussian gradient approximation, requiring three times more function queries. This phenomenon cannot be attributed to the chosen hyperparameter values, as experiments with different sets of values do not significantly improve the results. Our main hypothesis is that this phenomenon may come from a side effect of using the truncated Gaussian distribution. However, a comprehensive investigation of this issue requires dedicated research, left for future work.

\subsection{Solving problems under mixed aleatory/epistemic uncertainties} \label{sec-comp-c}

In this section, the behavior of the algorithm in the presence of mixed aleatory and epistemic uncertainties is examined. Epistemic uncertainties may arise from uncertainties about distribution parameters~\cite{nannapaneni2016reliability}. In the VSI problem presented in ~\cite{yang2020enriched}, it is noted that the mean of the uncertainty variables $\xi_8$ and $\xi_9$  can take two different values: $0.192$ and $0.345$. While both values were fixed to $0.345$ in~\cite{yang2020enriched} and in the previous experiments, in this section, these means are treated as epistemic uncertainties. Two types of epistemic uncertainty are studied: points epistemic uncertainty where the means $\mu_{\xi_8}$ and $\mu_{\xi_9}$ of $ \xi_8$ and  $\xi_9$  belong to $\{(0.192, 0.192), (0.192, 0.345), (0.345, 0.192),$$ (0.345, 0.345) \}$ and interval epistemic uncertainty where  $\mu_{\xi_8}$ and $\mu_{\xi_9}$ belong to the same interval $[0.192, 0.345]$. The others uncertain variables remain the same (no truncated) and are considered as aleatory uncertainties.

In this type of problems, a solution is deemed feasible if, for any values $\mu_{\xi_8}$ and $\mu_{\xi_9}$ , the probabilistic constraints are satisfied with a probability greater than $0.99$. Checking solution feasibility is more complex than in the previous section. In the case of points epistemic uncertainty, checking feasibility remains relatively straightforward since it involves evaluating the solution for the four possible pairs of means. However, when dealing with interval epistemic uncertainty, there is no ideal method for this verification. The approach adopted in this paper involves seeking the worst possible values of the epistemic uncertainties, $\mu_{\xi_8}$ and $\mu_{\xi_9}$, at a candidate solution $\x^*$. To achieve this, the following problem is solved for each constraint  $C_j, \, j \in [1, m]$
\begin{equation}
\label{prob_inter}
    \max_{(\mu_{\xi_8}, \mu_{\xi_9}) \in [0.192, 0.345]^2 } C_j(\x^*, \mathbb{E}[\boldsymbol{\xi}]).
\end{equation}
This problem aims to find the most challenging combination of $\mu_{\xi_8}$ and $\mu_{\xi_9}$. In this problem, all the uncertainties are fixed to their means and therefore the problem is deterministic. For each constraint, the couples solution of Problem (\ref{prob_inter}) are recorded. Next, the aleatory uncertainties are introduced. For each pair of $\mu_{\xi_8}$ and $\mu_{\xi_9}$ obtained , the probabilities of satisfying the constraints at $\x^*$ are computed using the original distribution of the aleatory uncertainties. If these probabilities are all larger than $0.99$, then the candidate solution is considered feasible. This approach provides a robust assessment of feasibility under interval epistemic uncertainty.
It is noteworthy that applying this methodology to the solution point obtained by the SORA algorithm reveals that this point is infeasible in the presence of epistemic uncertainty.  For instance, if the means $\mu_{\xi_8}$ and $\mu_{\xi_9}$ are taken to be equal to $(0. 192, 0.345)$, the probability of satisfying the $7\textup{th}$ constraint is $\mathbb{P} (C_7(\x_{SORA}^*, \boldsymbol{\xi}) \leq 0) \approx 0.88$. 
\begin{table}[ht!]
    \centering
        \caption{Average result over $100$ runs obtained with mixed aleatory/points epistemic uncertainty in the VSI problem with $15000$ function evaluations by run}
    \small
    \begin{tabular}{|C{2.5cm}||C{1.8cm}|C{5cm}|C{4cm}|C{2.5cm}|}
    \hline
           Value of $(\mu_{\xi_8}, \mu_{\xi_9})$ & Average of $\mathbb{E}[C(\x^*, \boldsymbol{\xi})]$
           &  Average of $ \mathbb{P}(C_j(\x^*, \boldsymbol{\xi}) \leq 0)$ & Average result point $\Bar{\x}^*$ (and standard deviation)& Number of successful runs   \\
    \hline  \hline
        \multirow{4}{*}{$(0.192, 0.192)$ } & \multirow{10}{*}{$30.38$} & \multirow{1}{*}{$[1.0, 1.0, 1.0, 0.9986, 1.0$} & \multirow{6}{*}{$\underset{[\pm 0.05, \, \pm 0.0, \, \pm 0.02, \, \pm 0.007}{[1.27, 1.35, 0.51, 1.49,}$} &  \multirow{3}{*}{$98$} \\
     & &\multirow{1}{*}{ $1.0, 0.9993, 0.9930, 1.0, 1.0]$} & & \\
    & &  &  \multirow{7}{*}{$\underset{\pm 0.08, \, \pm 0.01, \, \pm 0.1]}{1.26, 1.19, 0.47]}$} & \\
   
    \cline{1-1}  \cline{3-3} \cline{5-5}
     \multirow{2}{*}{$(0.192, 0.345)$} &  &  \multirow{1}{*}{$[1.0, 1.0, 1.0, 1.0, 1.0,$} &  &  \multirow{2}{*}{98}  \\
     & &  \multirow{1}{*}{$1.0, 0.9993, 0.9930, 1, 0.9995]$} &  & \\
    \cline{1-1}  \cline{3-3} \cline{5-5}
         \multirow{2}{*}{$(0.345, 0.192)$} &   &  \multirow{1}{*}{$[1.0, 1.0, 1.0, 1.0, 1.0,$} &    &  \multirow{2}{*}{98}  \\
     & &  \multirow{1}{*}{$1.0, 1.0, 0.9930, 1.0, 1.0]$} &  & \\
    \cline{1-1}  \cline{3-3} \cline{5-5}
         \multirow{2}{*}{$(0.345, 0.345)$} &   &  \multirow{1}{*}{$[1.0, 1.0, 1.0, 1.0, 1.0,$} &    &  \multirow{2}{*}{99}  \\
     & &  \multirow{1}{*}{$1.0, 1.0, 0.9929, 1.0, 0.9995]$} &  & \\
    \hline 
    \end{tabular}
    \label{table_res4}
\end{table}

To address this type of problems with the RAMSA algorithm, it is necessary to associate a probability distribution with the mean of $\xi_8$ and $\xi_9$. It is important to underline that this does not imply making an assumption about the distribution of the epistemic uncertainty itself. The distribution is just utilized to generate blackbox outputs. That allows to approach the problem from a worst-case perspective, leveraging the $\cvar$ properties when the values of $\alpha_j$ are taken sufficiently close to $1$. In the algorithm, the Bernoulli distribution is employed to generate the means for points epistemic uncertainty, while the uniform distribution is used to generate the means for interval epistemic uncertainty. The results for  mixed aleatory/points epistemic uncertainties are presented in~\Cref{table_res4}, and for mixed aleatory/interval epistemic uncertainties in~\Cref{table_res5}.
\begin{table}[ht!]
    \centering
        \caption{Average result over $100$ runs obtained with mixed aleatory/interval epistemic uncertainty in the VSI problem with $10000$ function evaluations by run.}
    \small
    \begin{threeparttable}
    \begin{tabular}{|C{2.5cm}||C{1.8cm}|C{5cm}|C{4cm}|C{2.5cm}|}
    \hline
           Solution of Problem (\ref{prob_inter}) \tnote{1} & Average of $\mathbb{E}[C(\x^*, \boldsymbol{\xi})]$
           &  Average of $ \mathbb{P}(C_j(\x^*, \boldsymbol{\xi}) \leq 0)$ & Average result point $\Bar{\x}^*$ (and standard deviation)& Number of successful runs  \\
    \hline
    \hline
        \multirow{4}{*}{$(0.192, 0.345)$ } & \multirow{4}{*}{$29.71$} & \multirow{1}{*}{$[1.0, 1.0, 1.0, 1.0, 1.0$} & \multirow{1}{*}{$\underset{[\pm 0.04, \, \pm 0.00005, \, \pm 0.01, \, \pm 0.01}{[1.15, 1.35, 0.51, 1.49,}$} &  \multirow{4}{*}{$99$} \\
     & &\multirow{1}{*}{ $1.0, 0.9974, 0.9929, 1.0, 0.9994]$} &  & \\
    & &  &  \multirow{1}{*}{$\underset{\pm 0.07, \, \pm 0.02, \, \pm 0.1]}{1.23, 1.19, 0.48]}$} & \\
    & & & & \\ 
    \hline
    \end{tabular}
    \begin{tablenotes}
    \item[1] There is only one point because for each run, the solutions $(\mu_{\xi_8}, \mu_{\xi_9})$ of (\ref{prob_inter}) are always the same.
  \end{tablenotes}
    \label{table_res5}
    \end{threeparttable}
\end{table}

In both cases, the RAMSA algorithm achieves satisfactory results. An interesting observation is that the results obtained with mixed aleatory/interval epistemic uncertainties are better to those with mixed aleatory/points epistemic uncertainties. This observation might appear counterintuitive since, in this experiment, points epistemic uncertainty is a subset of interval epistemic uncertainty. However, this phenomenon could be explained because the algorithm is better at handling continuous distributions than discrete distributions. The continuous nature of interval epistemic uncertainty could potentially make it more amenable for the gradient estimator, leading to enhanced performance in these cases.

\section{Concluding remarks} \label{sec_concl}

This work targets the constrained blackbox optimization problem given in~\Cref{problem_ref}, where the output of the blackbox is subject to uncertainties. To deal with the uncertainties, a $\cvar$-constrained problem formulation is adopted. This formulation allows the selection of the desired level of reliability. A smooth approximation of the $\cvar$-constrained problem is then derived by convolving the objective and constraint functions with a truncated multivariate Gaussian density. The use of the truncated Gaussian density, as opposed to the classical Gaussian density, ensures that sampling points are drawn within the bound constraints. Consequently, this approach avoids numerical failures that may occur when functions are undefined outside their bounds. Then, a Lagrangian relaxation is applied to handle the constraints. The resulting Lagrangian function possesses several appealing properties for optimization. First, it is infinitely differentiable since it is a sum of smooth approximations of the objective and constraint functions. Second, gradient estimators of the Lagrangian function can be computed with only two noisy blackbox outputs, making it computationally efficient. Theoretical bounds on the quality of the approximation have been derived. These bounds depend on the size of the problem, the value of the smoothing parameters, and the desired level of reliability. It is worth noting that it has been proved that for a reliability level sufficiently close to $1$, a feasible solution of the approximated problem remains a feasible solution of the original $\cvar$-constrained problem.

A new algorithm has been proposed to find a saddle point of the Lagrangian function. This algorithm is based on multi-timescale stochastic approximation updates. In this work, four different timescales are used. On the fastest timescale, the updates aggregate information about the gradient of the smooth Lagrangian function. On a first intermediate timescale, they estimate the value-at-risk of the objective and constraint functions. On a second intermediate timescale, the updates compute the optimal solution with respect to $\x$, while on the slowest timescale, the updates compute the optimal values of the Lagrangian multipliers. A convergence analysis based on Lyapunov theory shows that the different updates almost surely converge to a saddle point of the Lagrangian function. This point is locally optimal for the smooth approximation of the $\cvar$-constrained problem. Furthermore, using the previous result on the quality of the approximation, we prove that for reliability level values sufficiently close to one, this point is feasible and its value may be arbitrarily close to an optimal value of the $\cvar$-constrained problem.

Once theoretical results have been stated, details of the numerical implementations are given. These details mainly concern two transformations: one mapping the design variables into $[0, 1]^n$ and another mapping the blackbox outputs into $[-\frac{\pi}{2}, \frac{\pi}{2}]^{m+1}$. These transformations are designed to scale the design variables and the blackbox outputs, effectively reducing the number of hyperparameters. Then, numerical experiments are performed. In these experiments, the primary objective is to establish rules for selecting the values of the remaining hyperparameters. The results reveal that all hyperparameter values, except two, are independent of the problem and can be pre-specified using the values determined in this work. The first problem-dependent hyperparameter identified is the initial value of the step size for updating $\x$. It is determined that this value can be estimated from the norm of the gradient estimator at the starting point. The second problem-dependent hyperparameter is the value of the smoothing parameter. It is found that this parameter can be chosen in such a way that its value minimize the variance of the gradient estimator at the starting point. 

The secondary objective is to compare the effectiveness of the methods when truncated Gaussian gradient estimators are used instead of the classical Gaussian gradient estimator.
The proposed strategy for setting the hyperparameters is applied to experiments conducted with the truncated Gaussian gradient estimator. However, its use come at a cost. In the conducted experiments, it is observed that the truncated estimator is approximately three times less efficient than the classical Gaussian gradient estimator in terms of blackbox evaluations.

The tertiary objective of the experiments is to apply the algorithm to problems involving mixed aleatory/epistemic uncertainties. In these experiments, the epistemic uncertainties are related to the parameter distribution of the uncertainty variables. Two types of epistemic uncertainty are explored: points epistemic uncertainty and interval epistemic uncertainty. The algorithm demonstrated significant efficacy in handling both types of uncertainties. Notably, it performed particularly well in cases involving interval uncertainties, yielding promising results.

Future work will focus on validating these results using real-world industrial test cases. Additionally, there are plans to compare the RAMSA algorithm with other state-of-the-art algorithms to further assess its performance and competitiveness on problems subject to mixed aleatory/epistemic uncertainties.

\clearpage
\bibliographystyle{spmpsci} 
\bibliography{manuscript}

\clearpage
\appendix
\section{Proof of~\Cref{main_th}}
\label{appendix_proof}
First, two technical lemmas are stated to show that the iterates $\M^k$ and $\V^k$ are uniformly bounded almost surely. For this purpose, properties about the random gradient estimator must be shown.
\begin{Lem}
\label{bound_g}
    Under~\Cref{assum1}.3, the random gradient estimator $\Tilde{\g} := (\Tilde{\g}_\x, \Tilde{\g}_\t, \Tilde{\g}_{\boldsymbol{\xi}})$ is almost surely Lipschitz continuous with respect to $\x, \t$ and $\boldsymbol{\lambda}$. Moreover, $||\Tilde{\g}||$ is almost surely bounded.
\end{Lem}
\begin{proof}
    Let $(\x, \y) \in \mathcal{X}^2$, $(\t, \s) \in \R^{m+1} \times \R^{m+1}$ and consider any fixed realization of $\u, \v, \boldsymbol{\xi}_1$ and $\boldsymbol{\xi}_2$, it follows that for $\alpha_j \in (0, 1)$
    \begin{align*}
        &|\Tilde{V}_{\alpha_j}(\x + \beta_1 \u, (t_j + \beta_2 v_j, \boldsymbol{\xi}_1) - \Tilde{V}_{\alpha_j}(\y + \beta_1 \u, s_j + \beta_2 v_j, \boldsymbol{\xi}_2) | \\ &\leq |t_j - s_j| + \left| \big( C_j(\x + \beta_1 \u, \boldsymbol{\xi}_1) - (t_j + \beta_1 v_j) \big)^+ - \big( C_j(\y + \beta_1 \u, \boldsymbol{\xi}_2) - (s_j + \beta_1 v_j) \big)^+ \right|\\
        &\leq 2 |t_j - s_j| + |C_j(\x + \beta_1 \u, \boldsymbol{\xi}_1) - C_j(\x + \beta_1 \u, \boldsymbol{\xi}_2)| \leq 2 |t_j - s_j| + L_3||\x - \y|| \mbox{ a.s. },
    \end{align*}
    where the second inequality follows from $|\max(a, 0) - \max(b, 0)| \leq |a-b|$ and the third is due to~\Cref{assum1}.3. Therefore, $\Tilde{V}_{\alpha_j}$ is almost surely Lipschitz continuous with respect to $\x \in \mathcal{X}$ and $t_j \in \R$. As, $\Tilde{L}$ is a sum of almost surely Lipschitz continuous functions with respect to $\x$ and $\t$, it is also an almost surely Lipschitz continuous function. Moreover, $\Tilde{L}$ is a linear function with respect to $\lambda$ and thus Lipschitz continuous with respect to $\lambda$.
    
    Finally, by~\Cref{assum1}.3, we have for all $\x \in \mathcal{\X}$ and $\boldsymbol{\xi} \in \Xi$
    \begin{equation*}
        |C_j(\x, \boldsymbol{\xi})|- |C_j(\mathbf{0}, \mathbf{0})| \leq |C_j(\x, \boldsymbol{\xi}) - C_j(\mathbf{0}, \mathbf{0})| \leq \kappa_3(\boldsymbol{\xi}, \mathbf{0}) ||\x||.
    \end{equation*}
    Thus, the function $C_j$ is almost surely bounded. Since $\Tilde{L}$ is a sum of almost surely bounded functions, $\x, \t$ and $\boldsymbol{\lambda}$ are taken in compact sets and $\v$ and $\u$ are truncated Gaussian random vectors, it follows directly that $||\Tilde{\g}||$ is almost surely bounded.
\end{proof}

Once this was shown, $\M^k$ and $\V^k$ may be bounded. 
\begin{Lem}
\label{bound_M_V}
    The sequence of updates $\M^k$ and $\V^k$ are uniformly bounded with probability one. 
\end{Lem}
\begin{proof}
    Let $k \in \N$, we have 
    \begin{equation*}
        \M^{k+1} = s_4^k  \Tilde{\g}^k + \sum_{r = 0}^{k-1} s_4^l \prod_{q = r}^{k-1} (1-s_4^{q+1}) \Tilde{\g}^r + \prod_{q = 0}^{k} (1-s_4^{q}) \Tilde{\g}^0.
    \end{equation*}
    It follows directly by triangular inequality that 
    \begin{equation*}
        ||\M^{k+1}|| \leq s_4^k  ||\Tilde{\g}^k|| + \sum_{r = 0}^{k-1} s_4^l \prod_{q = r}^{k-1} (1-s_4^{q+1}) ||\Tilde{\g}^r|| + \prod_{q = 0}^{k} (1-s_4^{q}) ||\Tilde{\g}^0||.
    \end{equation*}
    Now according to~\Cref{bound_g}, for all $r \in \N$, the random gradient estimator is almost surely bounded. Therefore, we have
     \begin{equation*}
        ||\M^{k+1}|| \leq  \left( s_4^k  + \sum_{r = 0}^{k-1} s_4^l \prod_{q = r}^{k-1} (1-s_4^{q+1})  + \prod_{q = 0}^{k} (1-s_4^{q}) \right)  \sup_{r \in [0, k]} ||\Tilde{\g}^r|| < + \infty.
    \end{equation*}
    The same arguments may be applied for $\V^k$, thus the claim follows directly.  
\end{proof}

The remainder of the section is composed of four steps.\\

\textbf{Step 1: Convergence of $\M$ and $\V$ updates.} Since $\M$ and $\V$ converge on the fastest timescale, according to Lemma 1 in~\cite[chapter 6]{borkar_stochastic_2008}, the convergence properties of the updates in~\Cref{update_grad} may be analyzed for arbitrary quantities of $\x$, $\t$ and $\boldsymbol{\lambda}$ (here $\x = \x^k$, $\t = \t^k$ and $\boldsymbol{\lambda} = \boldsymbol{\lambda}^k$ are used). These updates may be rewritten as follows
\begin{align}
    \label{update_Mk}
    &\M^{k+1} = \M^k + s_4^k \left( \nabla L(\x^k, \t^k , \boldsymbol{\lambda}^k) - \M^k + \delta_\M^{k+1} \right),\\
    \label{update_Vk}
    &\V^{k+1} = \V^k + s_4^k \left( (\nabla L(\x^k, \t^k , \boldsymbol{\lambda}^k) )^2 + \mathbb{V}(\x^k, \t^k, \boldsymbol{\lambda}^k) - \V^k + \delta_\V^{k+1} \right),
\end{align}
where $\delta_\M^{k+1} = \Tilde{\g}^k - \nabla L^\beta(\x^k, \t^k , \boldsymbol{\lambda}^k)$ and $\delta_\V^{k+1} = (\Tilde{\g}^k)^2 - \nabla L^\beta(\x^k, \t^k , \boldsymbol{\lambda}^k)^2 - \mathbb{V}(\x^k, \t^k, \boldsymbol{\lambda}^k)$, with $\mathbb{V}(\x^k, \t^k, \boldsymbol{\lambda}^k) = \mathbb{E}[(\Tilde{\g}^k - \mathbb{E}[\Tilde{\g}^k | \mathcal{F}^k])^2| \mathcal{F}^k]
$ the variance conditioned by the associated sigma field  $\mathcal{F}^k = \sigma(\x^r, \t^r, \boldsymbol{\lambda}^r, \M^r, \V^r; r \leq k)$. Now, the following Lemma may be stated to prove the convergence properties of the updates $\M$ and $\V$.
\begin{Lem}
\label{conv_M_V}
    Consider  the following continuous time system dynamics of the updates,
    \begin{align}
        \label{odeM}
        \begin{split}
        &\Dot{\M} = h_1(\M, \x, \t, \boldsymbol{\lambda}) := \nabla L^\beta(\x, \t , \boldsymbol{\lambda}) - \M, \\
        &\Dot{\V} = h_2(\M, \x, \t, \boldsymbol{\lambda}) := (\nabla L^\beta(\x, \t , \boldsymbol{\lambda}))^2 + \mathbb{V}(\x, \t, \boldsymbol{\lambda}) - \V,\\
        & (\Dot{\x}, \Dot{\t}, \Dot{\boldsymbol{\lambda}})   = (\mathbf{0}, \mathbf{0}, \mathbf{0}).
        \end{split}
    \end{align}
    This o.d.e. has a globally asymptotically stable equilibrium 
    \begin{equation*}
        \big\{ \big(\nabla L^\beta(\x, \t , \boldsymbol{\lambda}), \nabla L^\beta(\x, \t , \boldsymbol{\lambda}))^2 + \mathbb{V}(\x, \t, \boldsymbol{\lambda}), \x, \t, \boldsymbol{\lambda} \big) \; \big| \; (\x, \t, \boldsymbol{\lambda}) \in \mathcal{X}\times \mathcal{T} \times \mathcal{L} \big\},
    \end{equation*}
    and  the sequences $(\M^k, \v^k, \x^k, \t^k ,\boldsymbol{\lambda}^k)$ converge almost surely to this equilibrium. 
\end{Lem}
\begin{proof}
The proof may be decomposed in two parts: the first part consists of analyzing the solutions of the two first o.d.e. given in~\Cref{odeM} and the second part consists of verifying that all the assumptions needed to apply Lemma 1 in~\cite[Chapter 6]{borkar_stochastic_2008} are satisfied. 

First, let $(\x, \t, \boldsymbol{\lambda}) \in \mathcal{X}\times \mathcal{T} \times \mathcal{L}$ be fixed and consider the following functions,
    \begin{align*}
        &\mathcal{L}_{\x, \t , \boldsymbol{\lambda}}^1(\M) = ||\nabla L^\beta(\x, \t, \boldsymbol{\lambda}) - \M||^2, \\
        &\mathcal{L}_{\x, \t , \boldsymbol{\lambda}}^2(\V) = ||(\nabla L^\beta(\x, \t, \boldsymbol{\lambda}))^2 + \mathbb{V}(\x, \t, \boldsymbol{\lambda}) - \V||^2.
    \end{align*}
Let denote $\M^* = \nabla L^\beta(\x, \t, \boldsymbol{\lambda})$ and $\V^* = (\nabla L^\beta(\x, \t, \boldsymbol{\lambda}))^2 + \mathbb{V}(\x, \t, \boldsymbol{\lambda})$ the equilibrium points of the two first equations in~\Cref{odeM}. The both functions satisfy the following conditions:
 \begin{itemize}
        \item They are globally positive definite, i.e, $\mathcal{L}_{\x, \t , \boldsymbol{\lambda}}^1(\M) > 0$, for all $\M \neq \M^*$ and $\mathcal{L}_{\x, \t , \boldsymbol{\lambda}}^2(\V) > 0$, for all $\V \neq \V^*$.
        \item They are radially unbounded since $||\M|| \to \infty \implies \mathcal{L}_{\x, \t , \boldsymbol{\lambda}}^1(\M) \to \infty$ and $||\V|| \to \infty \implies \mathcal{L}_{\x, \t , \boldsymbol{\lambda}}^2(\V) \to \infty$.
        \item The time derivatives of the both functions are globally negative definite since $\frac{d}{d \tau} \mathcal{L}_{\x, \t , \boldsymbol{\lambda}}^1(\M(\tau))  = - 2||\nabla L^\beta(\x, \t, \boldsymbol{\lambda}) - \M(\tau)||^2 $ and $\frac{d}{d \tau} \mathcal{L}_{\x, \t , \boldsymbol{\lambda}}^2(\V(\tau)) = -2||(\nabla L^\beta(\x, \t, \boldsymbol{\lambda}))^2 +\mathbb{V}(\x, \t, \boldsymbol{\lambda}) - \V(\tau)||^2 $.
    \end{itemize}
Thus, both functions are Lyapunov functions associated to the two first o.d.e. given in~\Cref{odeM}. By a corollary of the LaSalle invariance theorem (see for instance~\cite[Corollary 4.2]{Khalil_2002}), the equilibrium points $M^*$ and $V^*$ are globally asymptotically stable. Moreover, $\nabla L^\beta$ is Lipschitz with respect to $\x, \t$ and $\boldsymbol{\lambda}$ since it is a continuously differentiable function defined on a bounded space. The same may be applied for the function $(\nabla L^\beta)^2$. Finally, the function $\mathbb{V}$ is also Lipschitz, since $\Tilde{\g}$ is Lipschitz by~\Cref{bound_g}.

Now, we use the framework of the Lemma 1 in~\cite[Chapter 6]{borkar_stochastic_2008}. 
\begin{enumerate}[label=(\roman*)]
    \item By~\Cref{bound_M_V}, the updates $\M^k$ and $\V^k$ are uniformly bounded almost surely. The same goes for the updates $\x^k, \t^k$ and $\boldsymbol{\lambda}^k$ because of the projection operator.
    \item The functions $h_1$ and $h_2$ are Lipschitz continuous with respect to $\x, \t, \boldsymbol{\lambda}, \M$ and $\V$ by properties of $\nabla L^\beta$  and $\mathbf{V}$.
    \item The sequence $(\delta_\M^{k+1} )$ is a martingale difference sequence with respect to the  increasing sigma fields $\mathcal{F}^k = \sigma(\x^r, \t^r, \boldsymbol{\lambda}^r, \M^r, \V^r; r \leq k)$ since, by properties of truncated Gaussian smoothing, it follows that
    \begin{align*}
        \mathbb{E}[\delta_\M^{k+1} | \mathcal{F}^k] = \mathbb{E}[\Tilde{\g}^k | \mathcal{F}^k] - \nabla L^\beta(\x^k, \t^k, \boldsymbol{\lambda}^k) = 0.
    \end{align*}
    This sequence is also square integrable since
    \begin{align*}
        \mathbb{E}[||\delta_\M^{k+1}||^2 | \mathcal{F}^k]  \leq 2( \mathbb{E}[||\Tilde{g}||^2 | \mathcal{F}^k] + \mathbb{E}[||\nabla L^{\beta}(\x^k, \t^k, \boldsymbol{\lambda}^k)||^2]) < \infty,
    \end{align*}
    because $||a - b||^2 \leq 2(||a||^2 + ||b||^2)$, $\Tilde{g}$ is almost surely bounded by~\Cref{bound_g} and $\nabla L^\beta$ is a continuous function taking inputs in a compact set.
    \item The sequence  $(\delta_\V^{k+1} )$ is a martingale difference sequence with respect to $\mathcal{F}^k$ since 
    \begin{equation*}
        \mathbb{E}[\delta_\V^{k+1} | \mathcal{F}^k] = \mathbb{E}[(\Tilde{\g}^k)^2 | \mathcal{F}^k] - (\nabla L^\beta(\x^k, \t^k, \boldsymbol{\lambda}^k))^2 - \mathbb{V}(\x^k, \t^k, \boldsymbol{\lambda}^k) = 0,
    \end{equation*}
    by definition of conditional variance $\mathbb{V}(\x^k, \t^k, \boldsymbol{\lambda}^k) = \mathbb{E}[(\Tilde{\g}^k)^2|\mathcal{F}^k] - ( \mathbb{E}[\Tilde{\g}^k|\mathcal{F}^k])^2$ and is square integrable
    \begin{equation*}
        \mathbb{E}[||\delta_\V^{k+1}||^2 | \mathcal{F}^k] \leq 2( \mathbb{E}[||(\Tilde{\g})^2||^2] + ||(\nabla L^\beta(\x^k, \t^k, \boldsymbol{\lambda}^k))^2 + \mathbb{V}[\Tilde{\g}| \mathcal{F}^k] ||^2 < + \infty,
    \end{equation*}
    thanks to the same arguments as for $\delta_\M^{k+1}$.
    \item Finally, the step sizes $s_1^k$, $s_2^k$, $s_3^k$ and $s_4^k$ satisfy~\Cref{assum2}.
\end{enumerate}
Under these conditions, Lemma 1 in~\cite[Chapter 6]{borkar_stochastic_2008} may be applied, and the claim follows directly.
\end{proof}

\textbf{Step 2: Convergence of the $\t$-update.}
The $\t$-update converges on a faster timescale than the ones on $\x$ and $\boldsymbol{\lambda}$, while $\M$ and $\V$ converge faster than $\t$, thus, according to Lemma 1 in~\cite[Chapter 6]{borkar_stochastic_2008} the convergence 
of the $\t$ update may be proved for any arbitrary $\boldsymbol{\lambda}$ and $\x$ (here $\x = \x^k$ and $\boldsymbol{\lambda} = \boldsymbol{\lambda}^k$ are taken). Furthermore, in the $\M$-updates and $\V$-updates, as a result of~\Cref{conv_M_V} the following limits hold $|| \M^k - \nabla L^\beta(\x^k, \t^k, \boldsymbol{\lambda}^k) || \to 0$ and $|| \V^k - (\nabla L^\beta(\x^k, \t^k, \boldsymbol{\lambda}^k))^2 - \mathbb{V}(\x^k, \t^k, \boldsymbol{\lambda}^k) || \to 0$ almost surely. Consequently,  by defining
\begin{equation*}
    \nabla_\t^k L^\beta = \nabla_\t L^\beta(\x^k, \t^k, \boldsymbol{\lambda}^k) \text{ and } \mathbb{V}_\t^k =  \mathbb{V}_\t(\x^k, \t^k, \boldsymbol{\lambda}^k),
\end{equation*}
the update on $\t$ may be rewritten as follows 
\begin{align}
\label{update_tk}
    \t^{k+1} = \Pi_{\mathcal{T}} \left[\t^k + s_3^k \left( - \Psi_{\x^k, \boldsymbol{\lambda}^k}( \t^k)   + \delta_\t^{k+1} \right) \right], \text{ where } \left \{ \begin{array}{cll}
          &\Psi_{\x^k, \boldsymbol{\lambda}^k}( \t^k) &= \frac{ \nabla_\t^k L^\beta }{\sqrt{(\nabla_\t^k L^\beta)^2 + \mathbb{V}_\t^k } +  \epsilon}, \\
          &\delta_\t^{k+1} &= \Psi_{\x^k, \boldsymbol{\lambda}^k} (\t^k) - \frac{\M_\t^{k+1}}{\sqrt{\V_\t^{k+1}} + \epsilon}. 
    \end{array} \right.
\end{align}
 Now, the following Lemma may be stated to prove the convergence properties of the update $\t$.
\begin{Lem}
    \label{conv_t}
    Consider the following continuous time system dynamics of the updates,
    \begin{align}
    \label{odet}
        \begin{split}
            &\Dot{\t} = \Gamma_\t \left[ -\Psi_{\x, \boldsymbol{\lambda}}( \t) \right] =\Gamma_\t \left[ \frac{-\nabla_\t L^\beta(\x, \t, \boldsymbol{\lambda})}{\sqrt{\nabla_\t L^\beta(\x, \t, \boldsymbol{\lambda}) + \mathbb{V}_\t(\x, \t , \boldsymbol{\lambda})} + \epsilon}\right], \\
            &(\Dot{\x}, \Dot{\boldsymbol{\lambda}}) = (\mathbf{0}, \mathbf{0}),
        \end{split}
    \end{align}
    where 
    \begin{equation*}
        \Gamma_\t[  -\Psi_{\x, \boldsymbol{\lambda}}( \t) ] := \lim_{0 < \eta \to 0} \frac{\Pi_{\mathcal{T}}[\t - \eta  \Psi_{\x, \boldsymbol{\lambda}}( \t)] - \Pi_{\mathcal{T}}[\t]}{\eta}.
    \end{equation*}
    This o.d.e. has an asymptotically globally stable equilibrium \begin{equation*}
        \big\{ \big( \x, \t^*(\x, \boldsymbol{\lambda}), \boldsymbol{\lambda} \big) \; \big| \; (\x, \boldsymbol{\lambda}) \in \mathcal{X}\times \mathcal{L} \big\},
    \end{equation*}
    where $t^*(\x, \boldsymbol{\lambda}) = \{ \t \; | \; \Gamma_\t \left[   -\Psi_{\x, \boldsymbol{\lambda}}( \t) \right] = 0 \} $ and  the sequences $(\x^k, \t^k ,\boldsymbol{\lambda}^k)$ converge almost surely to this equilibrium. 
\end{Lem}
It is worth noting that $\Gamma_\t[K(\t)]$ is the left directional derivative of the function $\Pi_{\mathcal{\t}}[\t]$ in the direction of $K(\t)$. By using the left directional derivative $\Gamma_\t \left[  -\Psi_{\x, \boldsymbol{\lambda}}( \t) \right]$ in the gradient
descent algorithm for $\t$, the gradient will point in the descent direction along the boundary of $\mathcal{T}$
whenever the $\t$-update hits its boundary.
\begin{proof}
    Similar to the analysis made for the $\M$-update and $\V$-update, the proof is decomposed in two parts. First, the solution of the first o.d.e. given in~\Cref{odet} is described. Let $(\x, \boldsymbol{\lambda}) \in \mathcal{X}\times \mathcal{L}$ be fixed and consider the following function 
    \begin{equation*}
        \mathcal{L}_{\x, \boldsymbol{\lambda}} (\t) = L^\beta(\x, \t, \boldsymbol{\lambda}) - L^\beta(\x, \t^*, \boldsymbol{\lambda}),
    \end{equation*}
    where $\t^*$ is a minimum point (for any $(\x, \boldsymbol{\lambda})$, the function $L^\beta$ is convex in $\t$). This function satisfies the following conditions:
    \begin{itemize}
        \item The function is positive definite since $\mathcal{L}_{\x, \boldsymbol{\lambda}}(\t) > 0$, for all $\t \neq \t^*$ and radially unbounded since $||\t|| \to \infty, \implies \mathcal{L}_{\x, \boldsymbol{\lambda}}(\t) \to \infty$.
        \item The time derivative of the function is 
        \begin{equation*}
            \frac{d \mathcal{L}_{\x, \boldsymbol{\lambda}} (\t) }{d \tau} =  \nabla_\t L^\beta(\x, \t , \boldsymbol{\lambda} )^T \; \Gamma_\t \left[   -\Psi_{\x, \boldsymbol{\lambda}}( \t) \right]
        \end{equation*}
        and the goal is to show that this quantity is negative definite. There are two sets of cases to study:
        \begin{itemize}
            \item The cases where $\t \in \mathcal{T}^\circ = \mathcal{T} \backslash \partial \mathcal{T}$. In all this cases, there exist $\eta > 0$ sufficiently small such that $t - \eta \Psi_{\x, \boldsymbol{\lambda}}( \t) \in \mathcal{T}$, therefore by definition of $\Gamma_\t$ and $\Psi_\t$, it follows that (recall that the operators on the vectors are elementwise):
        \begin{equation*}
             \frac{d \mathcal{L}_{\x, \boldsymbol{\lambda}} (\t) }{d \tau} =  - \sum_{j = 0}^{m} \frac{ \left( \frac{\partial L^\beta(\x, \t, \boldsymbol{\lambda}) }{\partial t_j } \right)^2}{\sqrt{\left( \frac{\partial L^\beta(\x, \t, \boldsymbol{\lambda}) }{\partial t_j } \right)^2 + \mathbb{V}_{t_j}(\x, \t, \boldsymbol{\lambda})} + \epsilon}.
        \end{equation*}
        \item The cases where $\t \in \partial \mathcal{T}$. When $\t \in \partial \mathcal{T}$, the indices $j \in [0, m]$ of the variables of $\t$ may be grouped in three complementary sets : $S^{\min} = \{j \in [0, m]\; | \; t_j = - (\t_{\max})_j \}, S^{\max} = \{j \in [0, m]\; | \; t_j = (\t_{\max})_j \}$ or $S^\circ = \{j \in [0, m] \; | \; t_j = (- (\t_{\max})_j,  (\t_{\max})_j) \}$. Then, for the variables $t_j$ whose the indices are in $S^\circ$, there exists $\eta > 0$, sufficiently small such that $(t - \eta \Psi_{\x, \boldsymbol{\lambda}} (\t))_j \in (- (\t_{\max})_j,  (\t_{\max})_j)$. For the variables $t_j$ whose the variables are in $S^{\min}$, then either  $(\Psi_{\x, \boldsymbol{\lambda}} (\t))_j \leq 0$, so $\Pi_{t_j}[(- \t_{\max} -\eta \Psi_{\x, \boldsymbol{\lambda}} (\t))_j] = (- \t_{\max} -\eta \Psi_{\x, \boldsymbol{\lambda}} (\t))_j$; or $(\Psi_{\x, \boldsymbol{\lambda}} (\t))_j > 0$ so $\Pi_{t_j}[- \t_{\max} -\eta \Psi_{\x, \boldsymbol{\lambda}} (\t))_j] = - (\t_{\max})_j$. For the variables $t_j$ whose the variables are in $S^{\max}$, the symmetric result may be obtained. Therefore, it follows that 
        \begin{align*}
             \frac{d \mathcal{L}_{\x, \boldsymbol{\lambda}} (\t) }{d \tau} &= \lim_{0 < \eta \to 0} \;  \nabla_\t L^\beta(\x, \t , \boldsymbol{\lambda} )^T \; \bigg(\frac{\Pi_\t[\t - \eta \Psi_{\x, \boldsymbol{\lambda}}(\t)] - \t}{\eta} \bigg) \\
             &= \lim_{0 < \eta \to 0} \Bigg( -\sum_{j \in S^\circ } \frac{\partial L^\beta(\x, \t, \boldsymbol{\lambda})}{\partial t_j} (\Psi_{\x, \boldsymbol{\lambda}}(\t))_j  \\
             &\qquad \qquad \; \: - \sum_{j \in S^{\min}  } \frac{\partial L^\beta(\x, \t, \boldsymbol{\lambda})}{\partial t_j} \frac{\Pi_{t_j}[(-\t_{\max} -\eta \Psi_{\x, \boldsymbol{\lambda}} (\t))_j]  + (\t_{\max})_j}{\eta} \\
             & \qquad \qquad \; \: - \sum_{j \in S^{\max}  } \frac{\partial L^\beta(\x, \t, \boldsymbol{\lambda})}{\partial t_j} \frac{\Pi_{t_i}[(\t_{\max} -\eta \Psi_{\x, \boldsymbol{\lambda}} (\t))_j] -  (\t_{\max})_j}{\eta} \bigg)\\
             &\leq   -\sum_{j \in S^\circ }  \frac{ \left( \frac{\partial L^\beta(\x, \t, \boldsymbol{\lambda}) }{\partial t_j } \right)^2}{\sqrt{\left( \frac{\partial L^\beta(\x, \t, \boldsymbol{\lambda}) }{\partial t_j } \right)^2 + \mathbb{V}_{t_j}(\x, \t, \boldsymbol{\lambda})} + \epsilon},
        \end{align*}
        \end{itemize}
    Therefore, $ \frac{d \mathcal{L}_{\x, \boldsymbol{\lambda}} (\t) }{d \tau} < 0$ whenever $\Gamma_\t \left[   -\Psi_{\x, \boldsymbol{\lambda}}( \t) \right] \neq 0$, i.e, is negative definite.
    \end{itemize}
Thus, the function $\mathcal{L}_{\x, \boldsymbol{\lambda}}$ is a Lyapunov function and by~\cite[Corrolary 4.2]{Khalil_2002}, the equilibrium point $t^*(\x, \boldsymbol{\lambda}) = \{ \t \; | \; \Gamma_\t \left[   -\Psi_{\x, \boldsymbol{\lambda}}( \t) \right] = 0 \} $ is  globally asymptotically stable. Moreover, since $\nabla L^\beta$ is Lipschitz continuous with respect to $\x$ and $\boldsymbol{\lambda}$, it follows that $t^*(\x, \boldsymbol{\lambda})$ is Lipschitz continuous with respect to these vectors as well.  Now, the framework of the Lemma 1 and Theorem 2 in~\cite[Chapter 6]{borkar_stochastic_2008} is used. 
\begin{itemize}
    \item The conditions (i) to (v) given in the proof of~\Cref{conv_M_V} are still satisfied.
    \item The function $\Gamma_\t \left[  -\Psi_{\x, \boldsymbol{\lambda}}( \t) \right]$ is Lipschitz continuous by properties of $\nabla L^\beta$.
    \item The random sequence $(\delta_t^{k+1})$ converges asymptotically to 0 by~\Cref{conv_M_V}.
\end{itemize}
Therefore, the $\t$-update is a stochastic approximation with a null martingale difference sequence term and an additional error term $\delta_t^{k+1}$. Then, by applying Theorem 2 in~\cite[Chapter 6]{borkar_stochastic_2008} and the enveloppe theorem~\cite[Theorem 16]{chow2017risk}, the claim follows directly. 
\end{proof}

\textbf{Step 3: Convergence of the $\x$-update.} The convergence of the $\x$-update is very similar to the convergence of the $\t$-update. The $\x$-update converges on a faster timescale than the one of $\boldsymbol{\lambda}$, while $\t$, $\M$ and $\V$ converge on faster timescales than $\x$, thus, according to~\cite[Chapter 6]{borkar_stochastic_2008} the convergence of the $\x$ update may be proved for any arbitrary $\boldsymbol{\lambda}$ (here $\boldsymbol{\lambda} = \boldsymbol{\lambda}^k$ is taken). Furthermore, in the $\t$, $\M$ and $\V$ updates, as a result of~\Cref{conv_M_V} and~\Cref{conv_t} the following limits hold $|| \M^k - \nabla L^\beta(\x^k, \t^k, \boldsymbol{\lambda}^k) || \to 0$, $|| \V^k - (\nabla L^\beta(\x^k, \t^k, \boldsymbol{\lambda}^k))^2 - \mathbb{V}(\x^k, \t^k, \boldsymbol{\lambda}^k) || \to 0$ and $||\t^k - \t^*(\x^k, \boldsymbol{\lambda}^k)|| \to 0$ almost surely. Consequently by defining 
\begin{equation*}
    \nabla_\x^k L^\beta = \nabla_\x L^\beta(\x^k, \t^k, \boldsymbol{\lambda}^k)\text{ and } \mathbb{V}_\x^k =  \mathbb{V}_\x(\x^k, \t^k, \boldsymbol{\lambda}^k),
\end{equation*}
the update on $\x$ may be rewritten as follows
\begin{align}
\label{update_xk}
    \x^{k+1} = \Pi_{\mathcal{X}} \left[ \x^k + s_2^k \left( - \Psi_{\boldsymbol{\lambda}^k}(\x^k)   + \delta_{1,\x}^{k+1} + \delta_{2,\x}^{k+1} \right) \right],
\end{align}
where 
\begin{align*}
    \Psi_{\boldsymbol{\lambda}^k}(\x^k) &= \frac{ \nabla_\x L^\beta(\x^k, \t^*(\x^k, \boldsymbol{\lambda}^k), \boldsymbol{\lambda}^k) }{\sqrt{(\nabla_\x L^\beta(\x^k, \t^*(\x^k, \boldsymbol{\lambda}^k), \boldsymbol{\lambda}^k))^2 + \mathbb{V}_\x(\x^k, \t^*(\x^k, \boldsymbol{\lambda}^k), \boldsymbol{\lambda}^k)} + \epsilon}, \\
    \delta_{1,\x}^{k+1} &= \frac{ \nabla_\x^k L^\beta}{\sqrt{( \nabla_\x^k L^\beta)^2 + \mathbb{V}_\x^k} + \epsilon}- \frac{\M_\x^{k+1}}{\sqrt{\V_\x^{k+1}} + \epsilon},\\
    \delta_{2,\x}^{k+1} &= \Psi_{ \boldsymbol{\lambda}^k}(\x^k) - \frac{ \nabla_\x^k L^\beta}{\sqrt{( \nabla_\x^k L^\beta)^2 + \mathbb{V}_\x^k} + \epsilon}.
\end{align*}
Now, the following Lemma may be stated to prove the convergence properties of the update $\x$.
\begin{Lem}
    \label{conv_x}
    Consider the following continuous time system dynamics of the updates,
    \begin{align}
    \label{odex}
        \begin{split}
            &\Dot{\x} = \Gamma_\x \left[ -\Psi_{ \boldsymbol{\lambda}}( \x) \right] =\Gamma_\x \left[ \frac{-\nabla_\x L^\beta(\x, \t^*(\x, \boldsymbol{\lambda}), \boldsymbol{\lambda})}{\sqrt{\nabla_\x L^\beta(\x, \t^*(\x, \boldsymbol{\lambda}), \boldsymbol{\lambda}) + \mathbb{V}_\x(\x, \t^*(\x, \boldsymbol{\lambda}) , \boldsymbol{\lambda})} + \epsilon}\right], \\
            &\Dot{\boldsymbol{\lambda}} = \mathbf{0},
        \end{split}
    \end{align}
    where 
    \begin{equation*}
        \Gamma_\x[  -\Psi_{\boldsymbol{\lambda}}( \x) ] := \lim_{0 < \eta \to 0} \frac{\Pi_{\mathcal{X}}[\x - \eta  \Psi_{\boldsymbol{\lambda}}( \x)] - \Pi_{\mathcal{X}}[\x]}{\eta}.
    \end{equation*}
    Assume there exists $K_1 \in \N$ such that $\x^{K_1}$ is in the domain of attraction of $\x^*$ where $\x^*$ is some local minimum of $L^\beta$ with respect to $\x$. Then, this o.d.e. has a locally asymptotically stable equilibrium
    \begin{equation}
    \label{sol_x}
        \big\{ \big( \x^*(\boldsymbol{\lambda}),  \boldsymbol{\lambda} \big) \; \big| \;  \boldsymbol{\lambda} \in  \mathcal{L} \big\},
    \end{equation}
    where $\x^*(\boldsymbol{\lambda}) = \{ \x \in \mathcal{X}  \; | \; \Gamma_\x \left[   -\Psi_{ \boldsymbol{\lambda}}( \x) \right] = 0 \}$ is the local minima of the assumption and  the sequences $(\x^k, \boldsymbol{\lambda}^k)$ converge almost surely to the set given in~\Cref{sol_x}.
\end{Lem}
\begin{proof}
     First, the solutions of the first o.d.e. in~\Cref{odex} is described. Let $\boldsymbol{\lambda} \in \mathcal{L}$ be fixed and consider the following function
    \begin{equation*}
        \mathcal{L}_{\boldsymbol{\lambda}}(\x) = L^\beta(\x, \t^*(\x, \boldsymbol{\lambda}), \boldsymbol{\lambda}) - L^\beta(\x^*, \t^*(\x^*, \boldsymbol{\lambda}), \boldsymbol{\lambda}),
    \end{equation*}
    where $\x^*$ is the local minimum in $\mathcal{X}$ defined in the statement of the Lemma. This function is locally positive definite and its time derivatives  is 
    \begin{equation*}
            \frac{d \mathcal{L}_{\boldsymbol{\lambda}} (\x) }{d \tau} =  \nabla_\x L^\beta(\x, \t^*(\x, \boldsymbol{\lambda}) , \boldsymbol{\lambda} )^T \; \Gamma_\x \left[   -\Psi_{ \boldsymbol{\lambda}}( \x) \right],
    \end{equation*}
    which is negative definite (the proof may be done in the exact same way as the one given in~\Cref{conv_t} and is omitted here). Therefore, the function is a Lyapunov function and, by Lyapunov stability theorem~\cite[Theorem 4.1]{Khalil_2002}, $\x^*(\boldsymbol{\lambda}) = \{ \x \; | \;\Gamma_\x[ -\Psi_{\boldsymbol{\lambda}}(\x) ] = 0 \}$ is a locally asymptotically stable equilibrium. Since $\nabla L^\beta$ is Lipschitz continuous with respect to $\boldsymbol{\lambda}$, it follows that $\x^*(\boldsymbol{\lambda})$ is Lipschitz as well.
    Now, the framework in~\cite[Chapter 6]{borkar_stochastic_2008} is used.
    \begin{itemize}
    \item The conditions (i) to (v) given in the proof of~\Cref{conv_M_V} are still satisfied.
    \item The function $\Gamma_\x \left[  -\Psi_{ \boldsymbol{\lambda}}( \x) \right]$ is Lipschitz continuous by properties of $\nabla L^\beta$.
    \item The random sequence $(\delta_{1, \x}^{k+1})$ and $(\delta_{2, \x}^{k+1})$ converges asymptotically to 0 by~\Cref{conv_M_V} and~\Cref{conv_t}.
    \end{itemize}
    By assumption, the iterates $\x^{K_1}$ belongs to the domain of attraction of $\x^*$ for some $K_1 \in \N$. By definition of the domain of attraction, $\x^k$ is in the domain of attraction for all $k \geq K_1$.  Thus, by applying Theorem 2 in~\cite[Chapter 6]{borkar_stochastic_2008} from the iteration $K$, the claim follows directly. 
\end{proof}
At this stage, the results obtained in~\Cref{conv_t} and~\Cref{conv_x} allows concluding that for any fixed $\boldsymbol{\lambda} \in \mathcal{L}$, the following holds:
\begin{equation*}
    (\x^k, \t^k) \to (\x^*(\boldsymbol{\lambda}), \t^*(\x^*(\boldsymbol{\lambda}), \boldsymbol{\lambda})) \in \mathcal{X} \times \mathcal{T}.
\end{equation*}
Moreover, $ \t^*(\x^*(\boldsymbol{\lambda})$ is a minimum of $L^\beta$ with respect to $\t$ while $\x^*(\boldsymbol{\lambda})$ is a local minimum of $L^\beta$ with respect to $\x$. Since we have 
\begin{equation*}
    \min_{\x \in \mathcal{X}} \Big( \min_{\t \in \mathcal{T}}   L(\x, \t, \boldsymbol{\lambda}) \Big) = \min_{(\x, \t) \in \mathcal{X} \times \mathcal{T}}   L(\x, \t, \boldsymbol{\lambda}),
\end{equation*}
it follows that this point is a local minimum for the function $L^\beta$.\\

\textbf{Step 4: Convergence of the $\boldsymbol{\lambda}$-update.} Since the $\boldsymbol{\lambda}$-update converges in the slowest time scale, according to previous analysis, the following limits hold $|| \M^k - \nabla L^\beta(\x^k, \t^k, \boldsymbol{\lambda}^k) || \to 0$, $|| \V^k - (\nabla L^\beta(\x^k, \t^k, \boldsymbol{\lambda}^k))^2 - \mathbb{V}(\x^k, \t^k, \boldsymbol{\lambda}^k) || \to 0$, $||\t^k - \t^*(\x^k, \boldsymbol{\lambda}^k)|| \to 0$ and $||\x^k \to \x^*(\boldsymbol{\lambda})|| \to 0$ almost surely. Therefore, by defining
\begin{align*}
    &\nabla_{\boldsymbol{\lambda}}^k L^\beta = \nabla_{ \boldsymbol{\lambda}} L^\beta(\x^k, \t^k, \boldsymbol{\lambda}^k), \mathbb{V}_{\boldsymbol{\lambda}}^k = \mathbb{V}_{\boldsymbol{\lambda}}(\x^k, \t^k, \boldsymbol{\lambda}^k), \nabla_{\boldsymbol{\lambda}}^* L^\beta = \nabla_{ \boldsymbol{\lambda}} L^\beta(\x^k, \t^*(\x^k, \boldsymbol{\lambda}^k), \boldsymbol{\lambda}^k)  \\
    &\text{ and } \mathbb{V}_{\boldsymbol{\lambda}}^* = \mathbb{V}(\x^k, \t^*(\x^k, \boldsymbol{\lambda}^k),\boldsymbol{\lambda}^k),
\end{align*}
the $\boldsymbol{\lambda}$-update rule can be re-written as follows
\begin{align}
\label{update_lambdak}
    \x^{k+1} = \Pi_{\mathcal{L}} \left[ \boldsymbol{\lambda}^k + s_1^k \left(  \Psi(\boldsymbol{\lambda}^k)   + \delta_{1,\boldsymbol{\lambda}}^{k+1} + \delta_{2,\boldsymbol{\lambda}}^{k+1}  + \delta_{3,\boldsymbol{\lambda}}^{k+1} \right) \right],
\end{align}
where
\small
\begin{align*}
    \Psi(\boldsymbol{\lambda}^k) &= \frac{ \nabla_{\boldsymbol{\lambda}} L^\beta(\x^*(\boldsymbol{\lambda}^k), \t^*(\x^*(\boldsymbol{\lambda}^k), \boldsymbol{\lambda}^k), \boldsymbol{\lambda}^k) }{\sqrt{(\nabla_{\boldsymbol{\lambda}} L^\beta(\x^*(\boldsymbol{\lambda}^k), \t^*(\x^*(\boldsymbol{\lambda}^k), \boldsymbol{\lambda}^k), \boldsymbol{\lambda}^k))^2 + \mathbb{V}_{\boldsymbol{\lambda}}(\x^*(\boldsymbol{\lambda}^k), \t^*(\x^*(\boldsymbol{\lambda}^k), \boldsymbol{\lambda}^k), \boldsymbol{\lambda}^k)} + \epsilon}, \\
    \delta_{1, \boldsymbol{\lambda}}^{k+1} &= \frac{\M_{ \boldsymbol{\lambda}}^{k+1}}{\sqrt{\V_{ \boldsymbol{\lambda}}^{k+1}} + \epsilon} - \frac{\nabla_{\boldsymbol{\lambda}}^k L^\beta}{\sqrt{(\nabla_{\boldsymbol{\lambda}}^k L^\beta)^2 + \mathbb{V}_{\boldsymbol{\lambda}}^k } + \epsilon}, \\
    \delta_{2,\boldsymbol{\lambda}}^{k+1} &= \frac{\nabla_{\boldsymbol{\lambda}}^k L^\beta}{\sqrt{(\nabla_{\boldsymbol{\lambda}}^k L^\beta)^2 + \mathbb{V}_{\boldsymbol{\lambda}}^k } + \epsilon} - \frac{\nabla_{\boldsymbol{\lambda}}^* L^\beta}{\sqrt{(\nabla_{\boldsymbol{\lambda}}^* L^\beta)^2 + \mathbb{V}_{\boldsymbol{\lambda}}^*} + \epsilon},\\
    \delta_{3,\boldsymbol{\lambda}}^{k+1} &= \frac{\nabla_{\boldsymbol{\lambda}}^* L^\beta}{\sqrt{(\nabla_{\boldsymbol{\lambda}}^* L^\beta)^2 + \mathbb{V}_{\boldsymbol{\lambda}}^* }+ \epsilon} - \Psi(\boldsymbol{\lambda}^k).
\end{align*}
\normalsize 
 Now, the following Lemma may be stated to prove the convergence properties of the update $\boldsymbol{\lambda}$.
 \begin{Lem}
    \label{conv_lambda}
    Let consider the following continuous time system dynamics of the updates,
    \footnotesize
    \begin{align}
    \label{odelambda}
        \Dot{\boldsymbol{\lambda}} = \Gamma_{\boldsymbol{\lambda}} \left[ \Psi(\boldsymbol{\lambda}) \right] =\Gamma_{\boldsymbol{\lambda}} \left[  \frac{ \nabla_{\boldsymbol{\lambda}} L^\beta(\x^*(\boldsymbol{\lambda}^k), \t^*(\x^*(\boldsymbol{\lambda}^k), \boldsymbol{\lambda}^k), \boldsymbol{\lambda}^k) }{\sqrt{(\nabla_{\boldsymbol{\lambda}} L^\beta(\x^*(\boldsymbol{\lambda}^k), \t^*(\x^*(\boldsymbol{\lambda}^k), \boldsymbol{\lambda}^k), \boldsymbol{\lambda}^k))^2 + \mathbb{V}_{\boldsymbol{\lambda}}(\x^*(\boldsymbol{\lambda}^k), \t^*(\x^*(\boldsymbol{\lambda}^k), \boldsymbol{\lambda}^k), \boldsymbol{\lambda}^k)} + \epsilon}\right],
    \end{align}
    \normalsize
    where 
    \begin{equation*}
        \Gamma_{\boldsymbol{\lambda}}[  \Psi(\boldsymbol{\lambda}) ] := \lim_{0 < \eta \to 0} \frac{\Pi_{\mathcal{L}}[\boldsymbol{\lambda} - \eta  \Psi(\boldsymbol{\lambda})] - \Pi_{\mathcal{L}}[\boldsymbol{\lambda}]}{\eta}.
    \end{equation*}
    Assume there exists $K_2 \in \N$ such that $\boldsymbol{\lambda}^{K_2}$ is in the domain of attraction of $\boldsymbol{\lambda}^*$ where $\boldsymbol{\lambda}^*$ is some local maximum of $L^\beta$ with respect to $\boldsymbol{\lambda}$. Then, this o.d.e. has a locally asymptotically stable equilibrium
    \begin{equation}
    \label{sol_lambda}
        \boldsymbol{\lambda}^* = \{ \boldsymbol{\lambda} \in \mathcal{L} \; | \; \Gamma_{\boldsymbol{\lambda}} \left[  \Psi( \boldsymbol{\lambda}) \right] = 0 \},
    \end{equation}
   and  the sequences $(\boldsymbol{\lambda}^k)$ converges almost surely to this local maximum given in~\Cref{sol_lambda}.
\end{Lem}
\begin{proof}
    The proof is analog to the proof of convergence for the $\x$-update. First, the solutions of the first o.d.e. in~\Cref{odelambda} is described. Let  consider the following function
    \begin{equation*}
        \mathcal{L}(\boldsymbol{\lambda}) = -L^\beta(\x^*(\boldsymbol{\lambda}), \t^*(\x^*(\boldsymbol{\lambda}), \boldsymbol{\lambda}), \boldsymbol{\lambda}) + L^\beta(\x^*(\boldsymbol{\lambda}^*), \t^*(\x^*(\boldsymbol{\lambda}^*), \boldsymbol{\lambda}^*), \boldsymbol{\lambda}^*),
    \end{equation*}
    where $\boldsymbol{\lambda}^*$ is the local maximum in $\mathcal{L}$ defined in the statement of the Lemma. This function is locally positive definite and its time derivatives  is 
    \begin{equation*}
            \frac{d \mathcal{L} (\boldsymbol{\lambda}) }{d \tau} =  \nabla_{\boldsymbol{\lambda}} L^\beta(\x^*(\boldsymbol{\lambda}), \t^*(\x^*(\boldsymbol{\lambda}), \boldsymbol{\lambda}) , \boldsymbol{\lambda} )^T \; \Gamma_{\boldsymbol{\lambda}} \left[   \Psi( \boldsymbol{\lambda}) \right],
    \end{equation*}
    which is negative definite (the proof may be done in the exact same way as the one given in~\Cref{conv_t} and is omitted here). Therefore, the function is a Lyapunov function and, by Lyapunov stability theorem~\cite[Theorem 4.1]{Khalil_2002}, $\boldsymbol{\lambda}^* = \{ \boldsymbol{\lambda} \; | \;\Gamma_{\boldsymbol{\lambda}}[ \Psi(\boldsymbol{\lambda}) ] = 0 \}$ is a locally asymptotically stable equilibrium. 
    Now, the framework in~\cite[Chapter 6]{borkar_stochastic_2008} is used.
    \begin{itemize}
    \item The conditions (i) to (v) given in the proof of~\Cref{conv_M_V} are still satisfied.
    \item The function $\Gamma_{\boldsymbol{\lambda}} \left[  \Psi( \boldsymbol{\lambda}) \right]$ is Lipschitz continuous by properties of $\nabla L^\beta$.
    \item The random sequence $(\delta_{1,\boldsymbol{\lambda} }^{k+1})$, $(\delta_{2, \boldsymbol{\lambda}}^{k+1})$ and $(\delta_{3, \boldsymbol{\lambda}}^{k+1})$  converges asymptotically to 0 by~\Cref{conv_M_V}, ~\Cref{conv_t} and~\Cref{conv_x}.
    \end{itemize}
    By assumption, the iterates $\boldsymbol{\lambda}^k$ belongs to the domain of attraction of $\boldsymbol{\lambda}^*$ for some $K_2 \in \N$. By definition of the domain of attraction, $\boldsymbol{\lambda}^k$ is in the domain of attraction for all $k \geq K_2$.  Thus, by applying Theorem 2, in~\cite[Chapter 6]{borkar_stochastic_2008} from the iteration $K = \max(K_1, K_2)$, the claim follows directly. 
\end{proof}
\textbf{Main result: convergence to a saddle point.}
By letting $\x^* = \x^*(\boldsymbol{\lambda}^*)$ and $\t^* = \t^*(\x^*(\boldsymbol{\lambda}^*), \boldsymbol{\lambda}^*)$, it  will be shown that $(\x^*, \t^*, \boldsymbol{\lambda}^*)$ is a  saddle point of the Lagrangian function $L^\beta$ if $\boldsymbol{\lambda}^* \in \mathcal{L}^\circ$ and thus by the  saddle point theorem, $\x^*$ is a locally optimal solution for the smooth CVaR-constrained problem given in~\Cref{problem_cvar_smooth}. This result is formally settled in~\Cref{main_th} which is recalled here;
\begin{Th}
    Under~\Cref{assum1}.3 and~\Cref{assum2}, let further assume that the problem given in~\Cref{problem_cvar_smooth} is strictly feasible and there exists $K \in \N$  such that $\x^K$ and $\boldsymbol{\lambda}^K$ are in the domain of attraction of $\x^*$ and $\boldsymbol{\lambda}^*$ with $\boldsymbol{\lambda}^* \in \mathcal{L}^\circ$ respectively. Then, the iterates $(\x^k, \t^k, \boldsymbol{\lambda}^k)$ converge almost surely to a saddle point of the Lagrangian function $L^\beta$ and $\x^*$ is a locally optimal solution for the smooth CVaR-constrained problem given by~\Cref{problem_cvar_smooth}.
\end{Th}
\begin{proof}
   Under the assumptions of the theorem, since $(\x^*, \t^*)$
 is a local minimum of $L^\beta(\x, \t, \boldsymbol{\lambda})$ over the bounded set $(\x, \t) \in \mathcal{X} \times \mathcal{T}$, there exists a $r > 0$ such that
 \begin{equation*}
     L^\beta(\x^*, \t^*, \boldsymbol{\lambda}^*) \leq L^\beta(\x, \t, \boldsymbol{\lambda}^*), \; \forall (\x, \t) \in \mathcal{X} \times \mathcal{T} \cap \mathcal{B}_r(\x^*, \t^*).
 \end{equation*}
 In order to complete the proof, we must show that for all $j \in [1,m]$
 \begin{align}
 \label{proof_ineq}
     &c_j(\x^*, \t^*) := t_j^* + \frac{1}{1-\alpha} \mathbb{E}_{\u, \v, \boldsymbol{\xi}}[(C(\x^* + \beta_1 \u, \boldsymbol{\xi}) - (t_j^* + \beta_2 v_j))^+] \leq 0 \; \text{ and} \\
\label{proof_eq}
     &\lambda_j^* c_j(\x^*, \t^*) = \lambda_j^* \left( t_j^* + \frac{1}{1-\alpha} \mathbb{E}_{\u, \v, \boldsymbol{\xi}}[(C(\x^* + \beta_1 \u, \boldsymbol{\xi}) - (t_j^* + \beta_2 v_j))^+] \right) = 0.
 \end{align}

The proof of the inequality given in~\Cref{proof_ineq} is made by contradiction. Suppose that 
\begin{equation*}
    c_j(\x^*, \t^*) = t_j^* + \frac{1}{1-\alpha} \mathbb{E}_{\u, \v, \boldsymbol{\xi}}[(C(\x^* + \beta_1 \u, \boldsymbol{\xi}) - (t_j^* + \beta_2 v_j))^+] > 0.
\end{equation*}
This implies for $\lambda_j \in \mathcal{L}_j^\circ$ that for any $\eta \in (0, \Bar{\eta}]$
\begin{align*}
    \Pi_{\mathcal{L}} \left[\lambda_j^* - \eta \left(t_j^* + \frac{1}{1-\alpha} \mathbb{E}_{\u, \v, \boldsymbol{\xi}}[(C(\x^* + \beta_1 \u, \boldsymbol{\xi}) - (t_j^* + \beta_2 v_j))^+] \right) \right] &= \Pi_{\mathcal{L}} \left[\lambda_j^* - \eta c_j(\x^*, \t^*) \right] \\
    &= \lambda_j^* - \eta c_j(\x^*, \t^*),
\end{align*}
with $\Bar{\eta}$ sufficiently small. Therefore, it follows that $\Gamma_{\lambda_j}[ (\Psi(\boldsymbol{\lambda}^*))_j ] = c_j(\x^*, \t^*) > 0$, which contradicts the definition of $\boldsymbol{\lambda}^*$ given in ~\Cref{sol_lambda}. Thus, the inequality given in~\Cref{proof_ineq} holds. To show the result given in~\Cref{proof_eq}, it is sufficient to show that $\lambda_j^* = 0$ when $c_j(\x^*, \t^*) < 0$. For $\lambda_j^* \in \mathcal{L}^\circ$, there exists a sufficiently small $\eta > 0$ such that 
\begin{equation*}
    \frac{\Pi_{\mathcal{L}} \left[\lambda_j^* 
+ \eta c_j(\x^*, \t^*)  \right] - \lambda_j^*}{\eta} =  c_j(\x^*, \t^*) < 0.
\end{equation*}
This is again in contradiction with the definition of $\boldsymbol{\lambda}^*$ given in ~\Cref{sol_lambda} and thus the equality in~\Cref{proof_eq} holds. Finally, by the local saddle point theorem, it follows that $\x^*$ is a locally optimal solution for the smooth CVaR-constrained problem given by~\Cref{problem_cvar_smooth}.
\end{proof}

\section{Analytical problems description}
\label{prob-def}
Here are the list of analytical problems considered in~\Cref{sec-num}. \\
\textbf{Steel column problem~\cite{yang2020enriched}}
\begin{itemize}
    \item Dimension: $n = 3$ and $m = 1$.
    \item Original lower bounds: $\b_\ell = (200, 10, 100)$
    \item Original upper bounds: $\b_u = (400, 30, 500)$
    \item Original $\x_0$: $ (200, 10.5, 100)$
    \item Equations:
    \begin{align*}
        &C_0(\x, \boldsymbol{\xi}) = (x_1 + \xi_1)( x_2 + \xi_2) + 5(x_3 + \xi_3),\\
        &C_1(\x, \boldsymbol{\xi}) = F \left(\frac{1}{A_s}  + \frac{\xi_8 e_b}{U_s(e_b - F)}\right) - \xi_4, \\
        &\text{with } A_s = 2(x_1 + \xi_1)(x_2 + \xi_2),  \, U_s = (x_1 + \xi_1)(x_2 + \xi_2)(x_3 + \xi_3), \, e_b = \frac{\pi^2  \xi_9  U_i}{L^2}, \\
        &U_i = \frac{1}{2}(x_1 + \xi_1)(x_2 + \xi_2)(x_3 + \xi_3)^2 \text{ and } F = \xi_5 + \xi_6 + \xi_7.  
    \end{align*}
    \item Uncertainties: $\xi_1 \sim \mathcal{N}(0, 0.1 x_1)$, $\xi_2 \sim \mathcal{N}(0, 0.1 x_2)$, $\xi_3 \sim \mathcal{N}(0, 0.1 x_3)$, $\xi_4 \sim \mathcal{N}(400, 40)$, $\xi_5 \sim~\mathcal{N}(5\times10^5, 5\times 10^4)$,  $\xi_6 \sim~\mathcal{N}(6\times10^5, 6\times 10^4)$, $\xi_7 \sim~\mathcal{N}(6\times10^5, 6\times 10^4)$, $\xi_8 \sim~\mathcal{N}(30, 3)$, $\xi_9 \sim~\mathcal{N}(21000, 2100)$ and $L = 7500$.
    \item Solution in~\cite{yang2020enriched}: $\x^* = (257.7806,  13.5335, 100)$ with $\mathbb{E}[C_0(\x^*, \boldsymbol{\xi})] = 3988.95$ and $\mathbb{P}(C_1(\x^*, \boldsymbol{\xi}) \leq 0) = 0.9947$ (estimated in this work from $10^6$ samples).
\end{itemize}
\textbf{Welded Beam problem~\cite{yang2020enriched}}
\begin{itemize}
    \item Dimension: $n = 4$ and $m = 5$.
    \item Original lower bounds: $\b_\ell = (3.175, 0.0, 0.0, 0.0)$
    \item Original upper bounds: $\b_u = (50.8, 254, 254, 50.8)$
    \item Original $\x_0$: $ (6.208, 157.82, 210.62, 6.208)$
    \item Equations:
    \begin{align*}
        &C_0(\x, \boldsymbol{\xi}) = \kappa_1(x_1 + \xi_1)^2(x_2 + \xi_2) 
        + \kappa_2 (x_3 + \xi_3)(x_4+ \xi_4)(\kappa_3 + x_2 + \xi_2) \\
        &C_1(\x, \boldsymbol{\xi}) = \frac{\tau }{93.77}  - 1 \, \text{ with }\\
        &  \tau = \sqrt{ \tau_1^2 + 2\frac{\tau_1 \tau_2(x_2 + \xi_2)}{2 R} + \tau_2^2}, \, \tau_1 = \frac{\kappa_4}{\sqrt{2}(x_1 + \xi_1)(x_2 + \xi_2)}, \\
        &R = \frac{\sqrt{(x_2 + \xi_2)^2 + (x_1 + \xi_1 + x_3 + \xi_3)^2}}{2}, \, M = \kappa_4 \left( \kappa_3 + \frac{x_2 + \xi_2}{2} \right),\\
        &J = \sqrt{2}(x_1 + \xi_1)(x_2 + \xi_2) \left( \frac{(x_2 + \xi_2)^2}{12} + \frac{(x_1 + \xi_1 + x_3 + \xi_3)^2}{4} \right), \, \tau_2 = \frac{M R}{J},\\
        &C_2(\x, \boldsymbol{\xi}) = \frac{\sigma}{206.85} - 1 \, \text{ with } \sigma = \frac{6 \kappa_4  \kappa_3}{(x_3 + \xi_3)^2(x_4 + \xi_4)},\\
        &C_3(\x, \boldsymbol{\xi}) = \frac{x_1 + \xi_1}{x_4 + \xi_4} - 1, \\
        & C_4(\x, \boldsymbol{\xi}) = \frac{\delta}{6.35} - 1 \, \text{ with } \delta = \frac{4 \kappa_4  (\kappa_3)^3}{2.0685\times10^5(x_3 + \xi_3)^3(x_4 + \xi_4)}, \\
        &C_5(\x, \boldsymbol{\xi}) = 1 - \frac{P}{\kappa_4 } \, \text{ with } P = \frac{4.013(x_3 + \xi_3)(x_4 + \xi_4)^3 \sqrt{\kappa_5  \kappa_6}}{6  (\kappa_3)^2} \left( 1 - \frac{x_3 + \xi_3}{4  \kappa_3} \sqrt{\frac{\kappa_5}{\kappa_6}} \right),
    \end{align*}
    where $\kappa_1 = 6.74135 \times 10^{-5}$, $\kappa_2 = 2.93585 \times 10^{-6}$, $\kappa_3 = 3.556 \times 10^2$, $\kappa_4 = 2.6688\times 10^4$, $\kappa_5 = 2.0685 \times 10^5$ and $\kappa_6 = 8.274 \times 10^4$.
    \item Uncertainties: $\xi_1 \sim \mathcal{U}(-0.1693, 0.1693)$, $\xi_2 \sim \mathcal{U}(-0.1693, 0.1693)$, $\xi_3 \sim \mathcal{U}(-0.0107, 0.0107)$, $\xi_4 \sim \mathcal{U}(-0.0107, 0.0107)$.
    \item Solution in~\cite{yang2020enriched}: $x^* = [5.9188, 181.2849, 210.6114, 6.2253]$ with $\mathbb{E}[C_0(\x^*, \boldsymbol{\xi})] = 2.4948$ and $\forall j \in [1, 5], \, \mathbb{P}(C_j(\x^*, \boldsymbol{\xi}) \leq 0) = 1.0$ (estimated from $10^6$ samples).
\end{itemize}
\textbf{Vehicle Side Impact problem~\cite{yang2020enriched}}
\begin{itemize}
    \item Dimension: $n = 7$ and $m = 10$.
    \item Original lower bounds: $\b_\ell = (0.5, 0.45, 0.5, 0.5, 0.875, 0.4, 0.4)$
    \item Original upper bounds: $\b_u = (1.5, 1.35, 1.5, 1.5, 2.625, 1.2, 1.2)$
    \item Original $\x_0$: $ (1.0, 1.0, 1.0, 1.0, 2.0, 1.0, 1.0)$
    \item Equations:
    \begin{align*}
        C_0(\x, \boldsymbol{\xi}) &= 1.98 + 4.9(x_1 + \xi_1) + 6.67(x_2 + \xi_2) + 6.98(x_3 + \xi_3) + 4.01(x_4 + \xi_4) + 1.78(x_5 + \xi_5)\\
        & \quad + 2.73(x_7 + \xi_7), \\
        C_1(\x, \boldsymbol{\xi}) &= 1.16- 0.3717(x_2 + \xi_2)(x_4 + \xi_4)  - 0.00931(x_2 + \xi_2)\xi_{10} -0.484(x_3 + \xi_3) \xi_9  \\
        & \quad +0.01343(x_6 + \xi_6) \xi_{10} -1,\\
        C_2(\x, \boldsymbol{\xi}) &= 0.261 -0.0159(x_1 + \xi_1)(x_2 + \xi_2) - 0.188(x_1 + \xi_1)\xi_8 - 0.019(x_2 + \xi_2)(x_7 + \xi_7) \\
        & \quad + 0.0144(x_3 + \xi_3)(x_5 + \xi_5) 
        + 0.00087570(x_5 + \xi_5)\xi_{10} + 0.08045(x_6 + \xi_6)\xi_9 \\
        & \quad +  0.00139\xi_{8}\xi_{11} 
         + 1.575\times 10^{-6} \xi_{10}\xi_{11} - 0.32,\\
    C_3(\x, \boldsymbol{\xi}) &= 0.2147 + 0.00817(x_5 + \xi_5) - 0.131(x_1 + \xi_1) \xi_8 -0.0704(x_1 + \xi_1) \xi_9 \\   
    & \quad + 0.03099(x_2 + \xi_2)(x_6 + \xi_6) - 0.018(x_2 + \xi_2)(x_7 + \xi_7)  
    + 0.0208(x_3 + \xi_3) \xi_8 \\
    & \quad + 0.121(x_3 + \xi_3) \xi_9 - 0.00364(x_5 + \xi_5)(x_6 + \xi_6) +0.0007715(x_5 + \xi_5) \xi_{10}\\
    & \quad-0.0005354(x_6 + \xi_6) \xi_{10} 
    + 0.00121 \xi_8 \xi_{11} +0.00184 \xi_9 \xi_{10} \\
    & \quad - 0.02(x_2 + \xi_2)^2 -0.32,\\
    C_4(\x, \boldsymbol{\xi}) &= 0.74 -0.61(x_2 + \xi_2) - 0.163(x_3 + \xi_3)\xi_8 + 0.001232(x_3 + \xi_3)\xi_{10} \\
    & \quad - 0.166(x_7 + \xi_7)\xi_9 + 0.227(x_2 + \xi_2)^2 - 0.32,\\
    C_5(\x, \boldsymbol{\xi}) &= 28.98 + 3.818(x_3 + \xi_3) - 4.2(x_1 + \xi_1)(x_2 + \xi_2) +0.0207(x_5 + \xi_5) \xi_{10} \\
    & \quad + 6.63(x_6 + \xi_6) \xi_9 - 7.77(x_7 + \xi_7) \xi_8 + 0.32 \xi_9 \xi_{10} - 32,\\
    C_6(\x, \boldsymbol{\xi}) &=33.86 +2.95(x_3 + \xi_3) +0.1792 \xi_{10} - 5.057(x_1 + \xi_1)(x_2 + \xi_2) -11(x_2 + \xi_2) \xi_8 \\
    & \quad - 0.0215(x_5 + \xi_5) \xi_{10} - 9.98(x_7 + \xi_7) \xi_8 + 22 \xi_8 \xi_9 - 32,\\
    C_7(\x, \boldsymbol{\xi}) &= 46.36 - 9.9(x_2 + \xi_2) - 12.9(x_1 + \xi_1) \xi_8 
    \\ & \quad + 0.1107(x_3 + \xi_3) \xi_{10} -32,\\
    C_8(\x, \boldsymbol{\xi}) &= 4.72 - 0.54(x_4 + \xi_4) - 0.19(x_2 + \xi_2)(x_3 + \xi_3) - 0.0122(x_4 + \xi_4) \xi_{10} \\
    & \quad + 0.009325(x_6 + \xi_6) \xi_{10} 
    + 0.000191 \xi_{11}^2 - 4,\\
    C_9(\x, \boldsymbol{\xi}) &= 10.58 - 0.674(x_1 + \xi_1)(x_2 + \xi_2) - 1.95(x_2 + \xi_2) \xi_8 + 0.028(x_6 + \xi_6) \xi_{10} \\
    & \quad + 0.02054(x_3 + \xi_3) \xi_{10} -0.0198(x_4 + \xi_4) \xi_{10} - 9.9,\\
    C_{10}(\x, \boldsymbol{\xi}) &= 16.45 - 0.489(x_3 + \xi_3)(x_7 + \xi_7) -0.843(x_5 + \xi_5)(x_6 + \xi_6) +0.0432 \xi_9 \xi_{10} \\
    & \quad- 0.0556 \xi_9 \xi_{11} - 0.000786 \xi_{11}^2 - 15.69.
    \end{align*}
   
    \item Uncertainties: $\forall i \in \{1, 2, 3, 4, 6, 7\}, \, \xi_i \sim \mathcal{N}(0, 0.03)$, $\xi_5 \sim \mathcal{N}(0, 0.05)$, $\xi_8 \sim \mathcal{N}(0.345, 0.006)$, $\xi_9 \sim \mathcal{N}(0.345, 0.006)$, $\xi_{10} \sim \mathcal{N}(0, 10)$ and $\xi_{11} \sim \mathcal{N}(0, 10)$.
    \item Solution in~\cite{yang2020enriched}: $x^* = (0.7872, 1.35  , 0.6887, 1.5   , 1.0706, 1.2   , 0.7284) $ with $\mathbb{E}[C_0(\x^*, \boldsymbol{\xi})] = 29.5585$ and $\forall j \in [1, 10], \, \mathbb{P}(C_j(\x^*, \boldsymbol{\xi}) \leq 0) \geq 0.9982$ (estimated from $10^6$ samples).
\end{itemize}
\textbf{Speed Reducer problem~\cite{chen2013optimal}}
\begin{itemize}
    \item Dimension: $n = 7$ and $m = 11$.
    \item Original lower bounds: $\b_\ell = (2.6, 0.7, 17, 7.3, 7.3, 2.9, 5.0)$
    \item Original upper bounds: $\b_u = (3.6, 0.8, 28, 8.3, 8.3, 3.9, 5.5)$
    \item Original $\x_0$: $  (3.5, 0.7, 17, 7.3, 7.72, 3.35, 5.29)$
    \item Equations:
    \begin{align*}
        C_0(\x, \boldsymbol{\xi}) &= 0.7854 (x_1 + \xi_1) (x_2 + \xi_2)^2(3.3333 (x_3 + \xi_3)^2 + 14.9334 (x_3 + \xi_3)-43.0934)\\
        & \quad- 1.508 (x_1 + \xi_1)( (x_6 + \xi_6)^2+ (x_7 + \xi_7)^2) 
        + 7.477( (x_6 + \xi_6)^3 +  (x_7 + \xi_7)^3) \\
        &\quad + 0.7854( (x_4 + \xi_4) (x_6 + \xi_6)^2 +  (x_5 + \xi_5) (x_7 + \xi_7)^2),\\ 
        C_1(\x, \boldsymbol{\xi}) &= \frac{27}{ (x_1 + \xi_1) (x_2 + \xi_2)^2 (x_3 + \xi_3)} -1,  \\
        C_2(\x, \boldsymbol{\xi}) &= \frac{397.5}{ (x_1 + \xi_1) (x_2 + \xi_2)^2 (x_3 + \xi_3)^2} - 1, \\
    C_3(\x, \boldsymbol{\xi}) &= \frac{1.93 (x_4 + \xi_4)^3}{ (x_2 + \xi_2) (x_3 + \xi_3) (x_6 + \xi_6)^4} - 1,\\
    C_4(\x, \boldsymbol{\xi}) &= \frac{1.93 (x_5 + \xi_5)^3}{ (x_2 + \xi_2) (x_3 + \xi_3) (x_7 + \xi_7)^4} - 1,\\
    C_5(\x, \boldsymbol{\xi}) &=  \frac{\sqrt{ \left(\frac{745 (x_5 + \xi_5)}{ (x_2 + \xi_2) (x_3 + \xi_3)} \right)^2+ 16.9\times10^6 }}{0.1 (x_6 + \xi_6)^3} - 1100,\\
    C_6(\x, \boldsymbol{\xi}) &= \frac{\sqrt{ \left( \frac{745 (x_5 + \xi_5)}{ (x_2 + \xi_2) (x_3 + \xi_3)} \right)^2 + 157.5\times10^6}}{0.1 (x_7 + \xi_7)^3} - 850,\\
    C_7(\x, \boldsymbol{\xi}) &=  (x_2 + \xi_2) (x_3 + \xi_3) - 40,\\
    C_8(\x, \boldsymbol{\xi}) &= 5 - \frac{ (x_1 + \xi_1)}{ (x_2 + \xi_2)}\\
    C_9(\x, \boldsymbol{\xi}) &= \frac{ (x_1 + \xi_1)}{ (x_2 + \xi_2)} -12,\\
    C_{10}(\x, \boldsymbol{\xi}) &= \frac{1.5 (x_6 + \xi_6)+1.9}{ (x_4 + \xi_4)} - 1,\\
    C_{11}(\x, \boldsymbol{\xi}) &= \frac{1.1 (x_7 + \xi_7) + 1.9}{ (x_5 + \xi_5)}  - 1.
    \end{align*}
   
    \item Uncertainties: $\forall i \in [1, 7], \, \xi_i \sim \mathcal{N}(0, 0.005)$.
    \item Solution in~\cite{yang2020enriched}: $x^* = (3.5765, 0.7, 17.0, 7.3, 7.7541, 3.3652,  5.3017)$ with $\mathbb{E}[C_0(\x^*, \boldsymbol{\xi})] = 3038.72$ and $\forall j \in [1, 11], \, \mathbb{P}(C_j(\x^*, \boldsymbol{\xi}) \leq 0) \geq 0.9976$ (estimated from $10^6$ samples).
\end{itemize}

\section{Detailed numerical results}
\label{detailed-res}

This section details the numerical results of~\Cref{sec-comp-a} and~\Cref{sec-comp-b}. In these sections, only the average result over the 100 runs are presented. In this section, boxplots are used to describe the result of all the 100 runs. Each run is represented by a cross, the orange line is the mediane and the bounds of the box are the first and third quartiles. Finally, the circled crosses are the outliers. Here are the results for~\Cref{sec-comp-a}.

\begin{figure}[htb!]
    \hbox{
    \begin{minipage}[c]{.5\linewidth}
        \centering
        \includegraphics[width =  \textwidth]{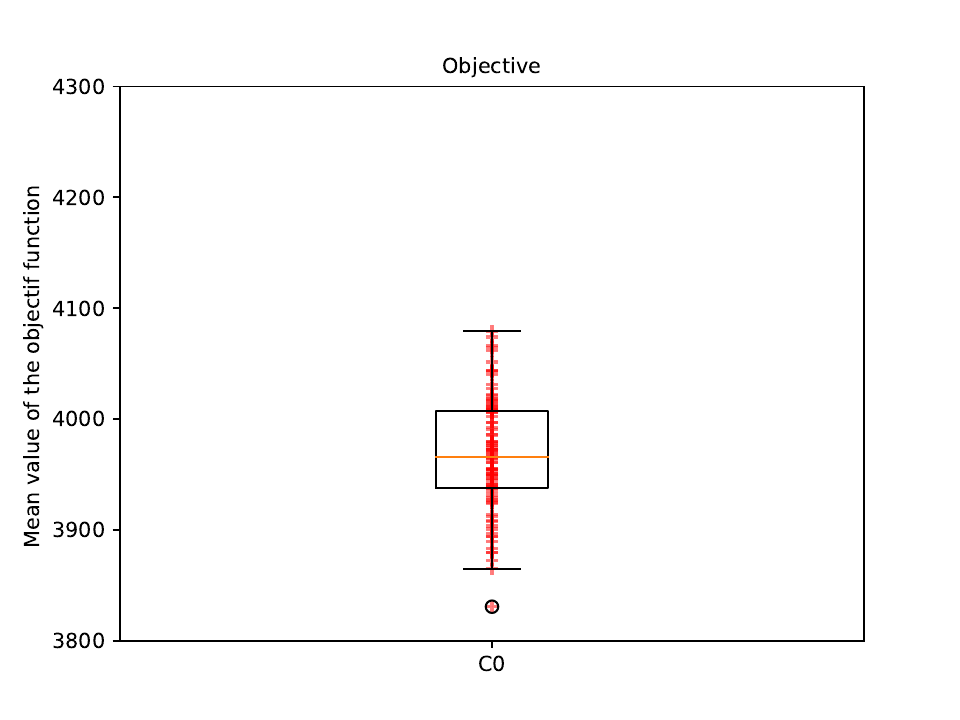}
    \end{minipage}
    \hfill
    \begin{minipage}[c]{.5\linewidth}
        \centering
        \includegraphics[width =  \textwidth]{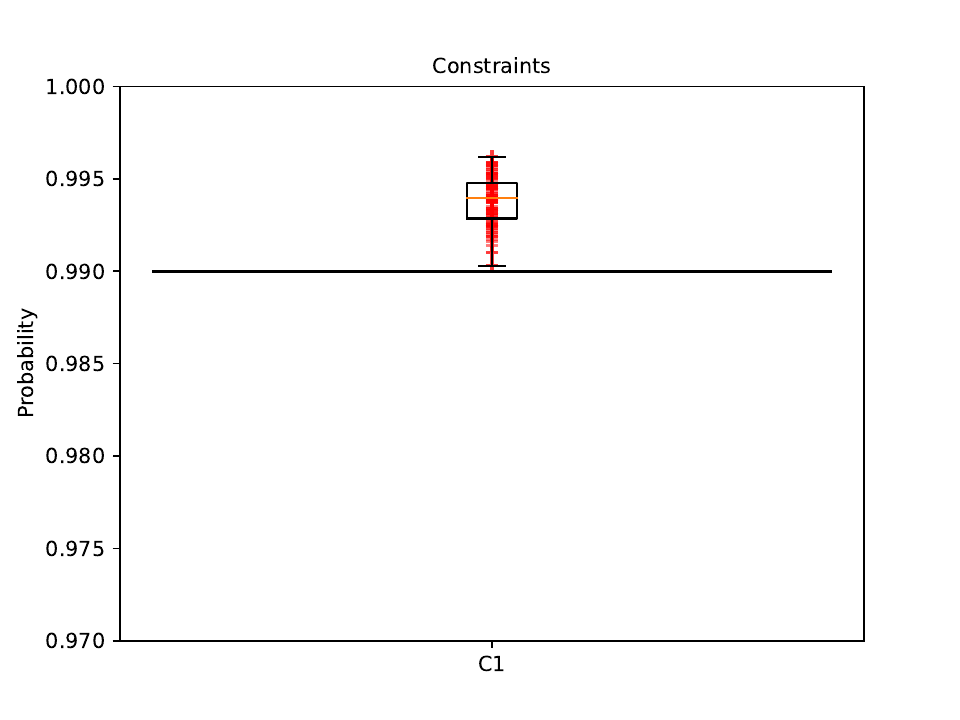}
    \end{minipage}
    }
    \caption{Detail result for Steel Column Design problem with classical Gaussian gradient approximation}
    \label{fig1}
\end{figure}

\begin{figure}[htb!]
    \hbox{
    \begin{minipage}[c]{.5\linewidth}
        \centering
        \includegraphics[width =  \textwidth]{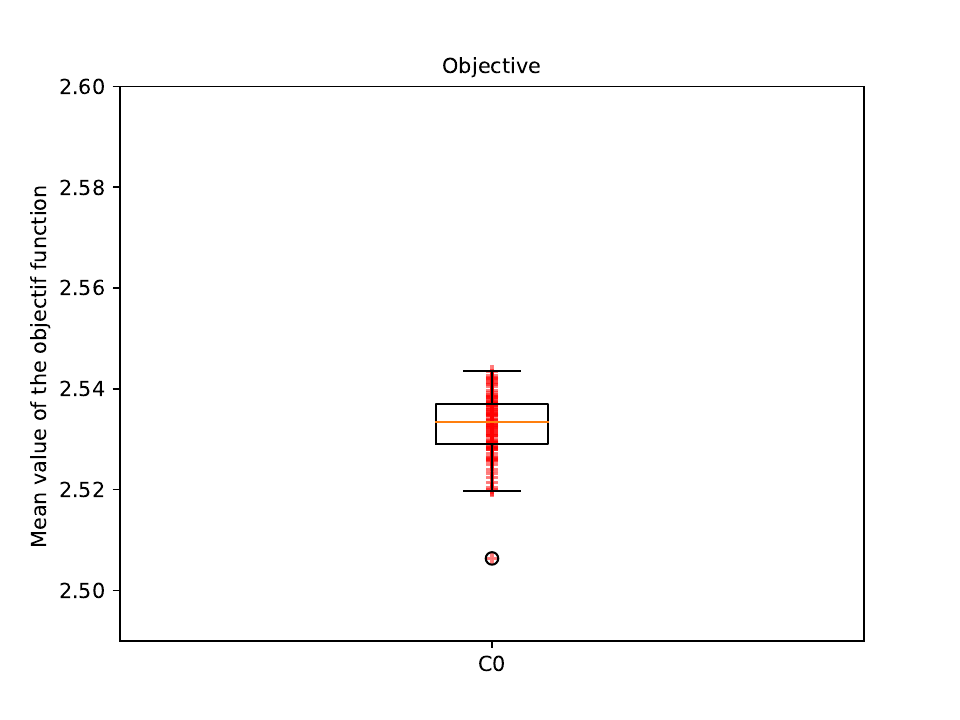}
    \end{minipage}
    \hfill
    \begin{minipage}[c]{.5\linewidth}
        \centering
        \includegraphics[width = \textwidth]{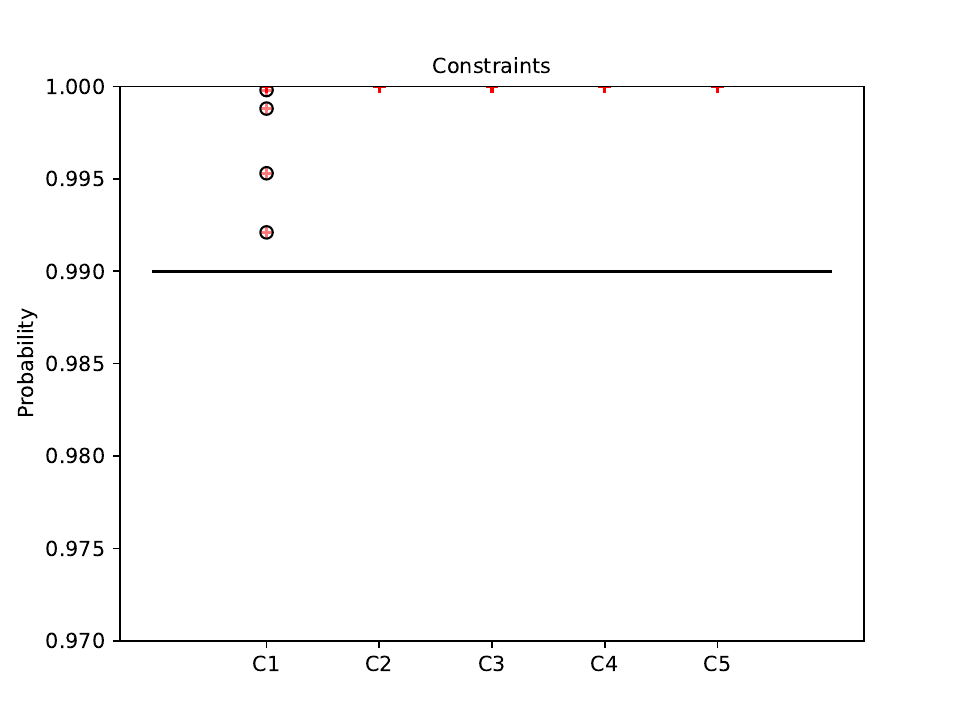}
    \end{minipage}
    }
    \caption{Detail result for Welded Beam Design problem with classical Gaussian gradient approximation}
    \label{fig2}
\end{figure}

\begin{figure}[htb!]
    \hbox{
    \begin{minipage}[c]{.5\linewidth}
        \centering
        \includegraphics[width =  \textwidth]{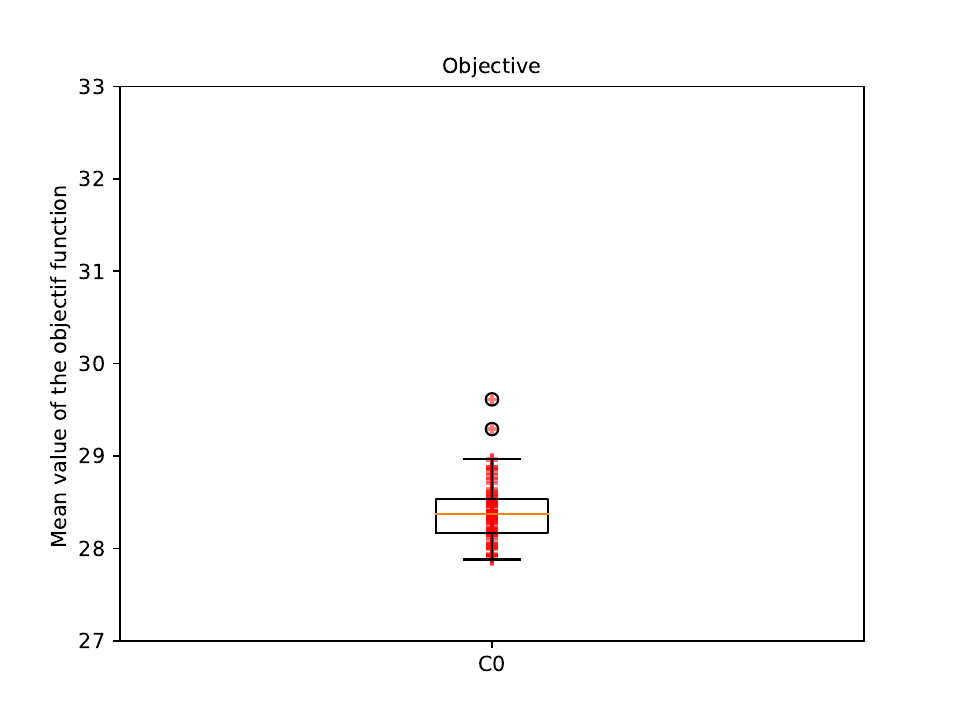}
    \end{minipage}
    \hfill
    \begin{minipage}[c]{.5\linewidth}
        \centering
        \includegraphics[width = \textwidth]{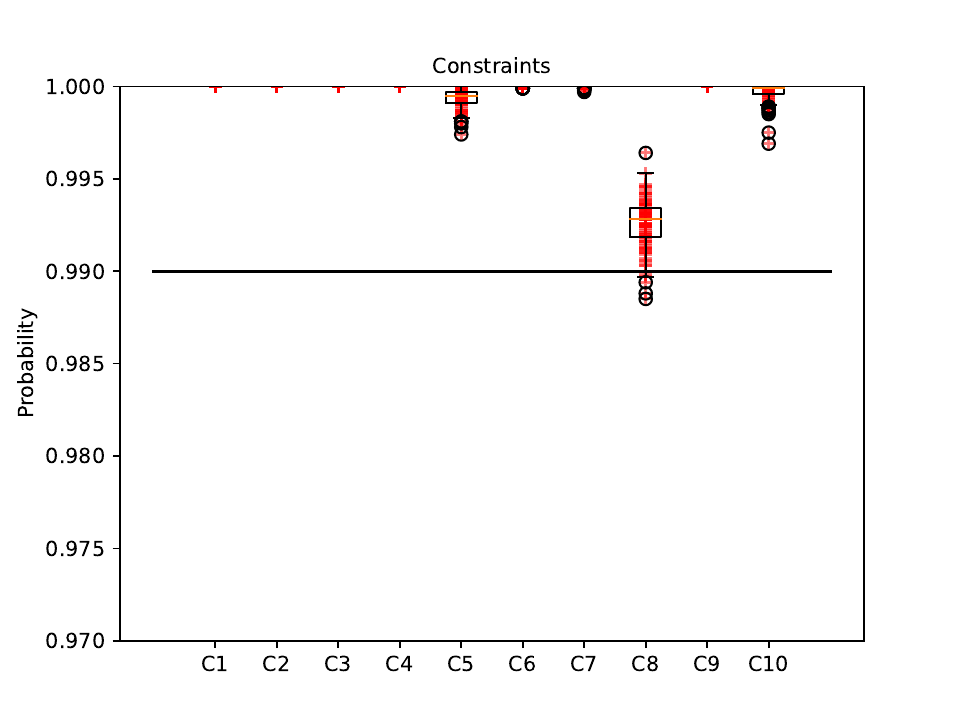}
    \end{minipage}
    }
    \caption{Detail result for Vehicle Side Impact problem with classical Gaussian gradient approximation}
    \label{fig3}
\end{figure}

\begin{figure}[htb!]
    \hbox{
    \begin{minipage}[c]{.5\linewidth}
        \centering
        \includegraphics[width =  \textwidth]{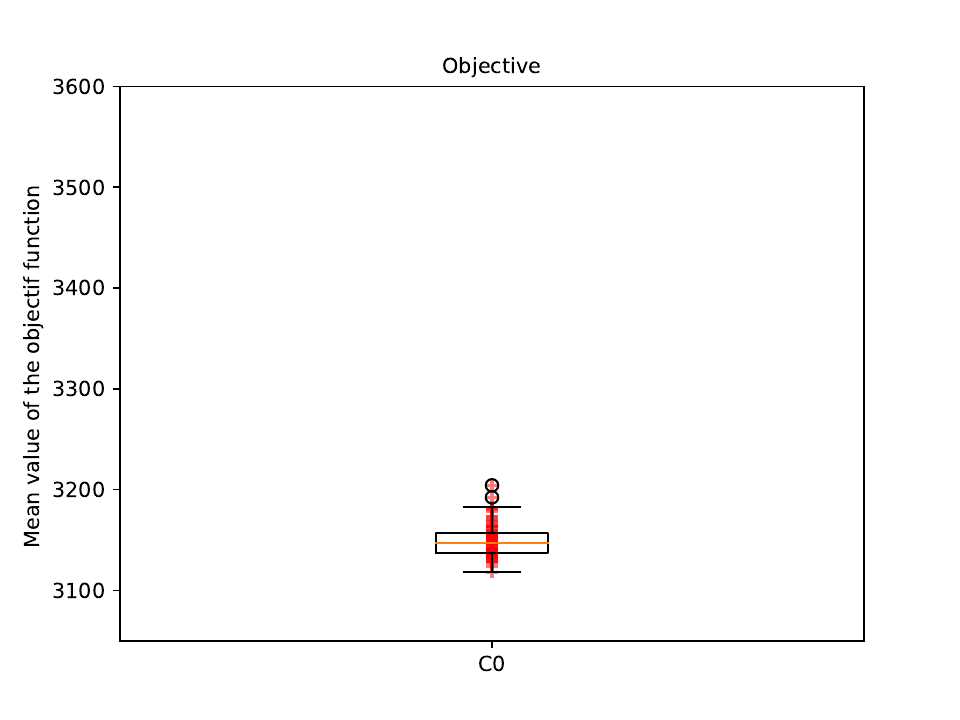}
    \end{minipage}
    \hfill
    \begin{minipage}[c]{.5\linewidth}
        \centering
        \includegraphics[width = \textwidth]{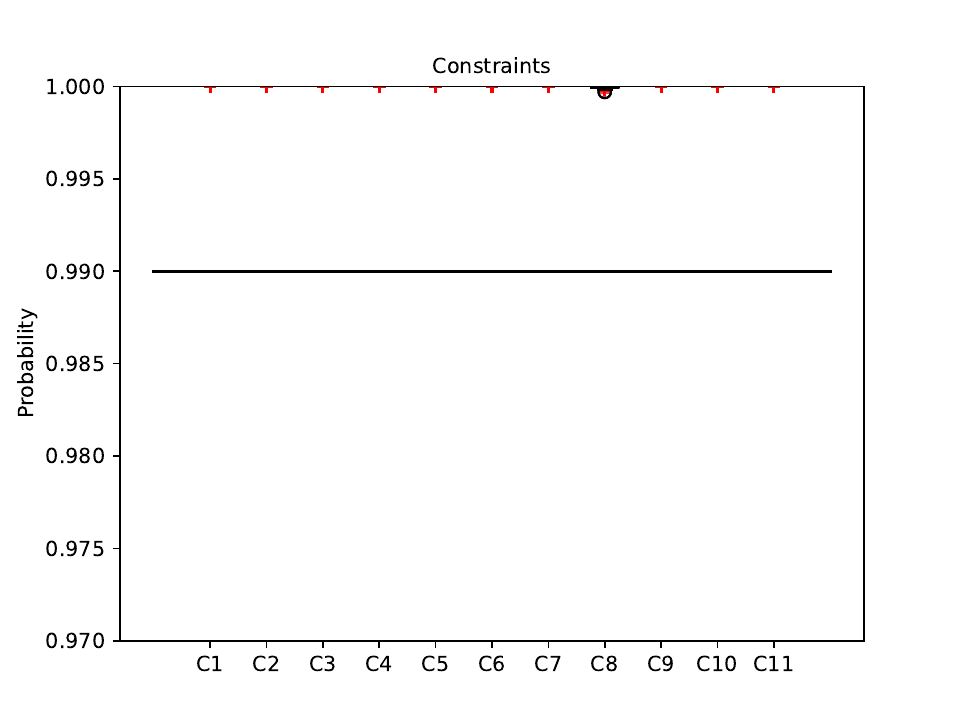}
    \end{minipage}
    }
    \caption{Detail result for Speed Reducer Design problem with classical Gaussian gradient approximation}
    \label{fig4}
\end{figure}
\newpage
Here are the results for~\Cref{sec-comp-b}.
\begin{figure}[htb!]
    \hbox{
    \begin{minipage}[c]{.5\linewidth}
        \centering
        \includegraphics[width = \textwidth]{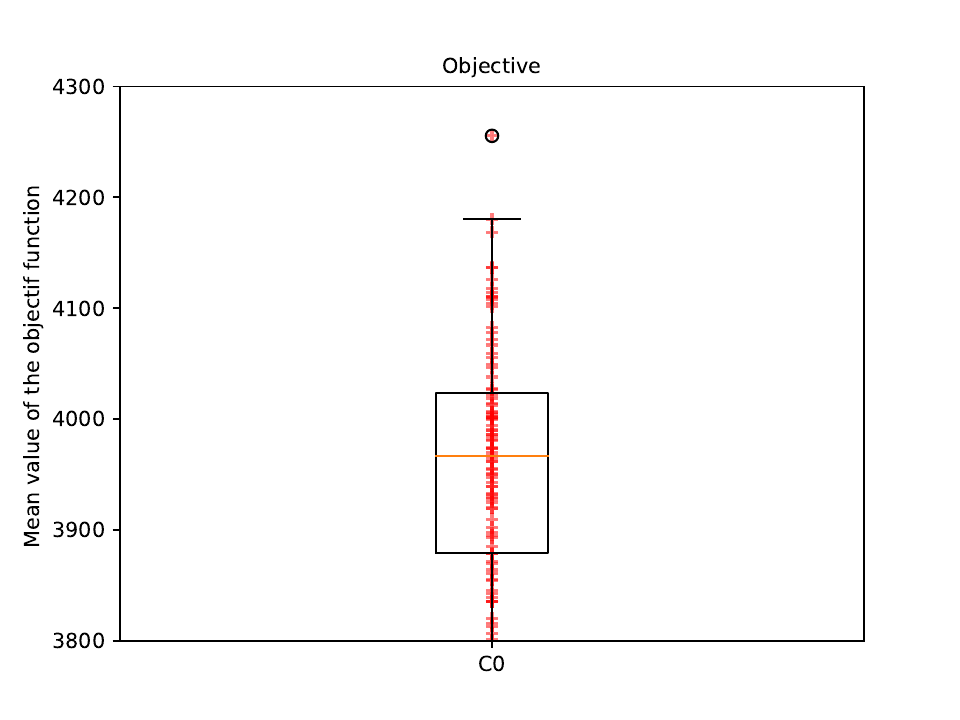}
    \end{minipage}
    \hfill
    \begin{minipage}[c]{.5\linewidth}
        \centering
        \includegraphics[width = \textwidth]{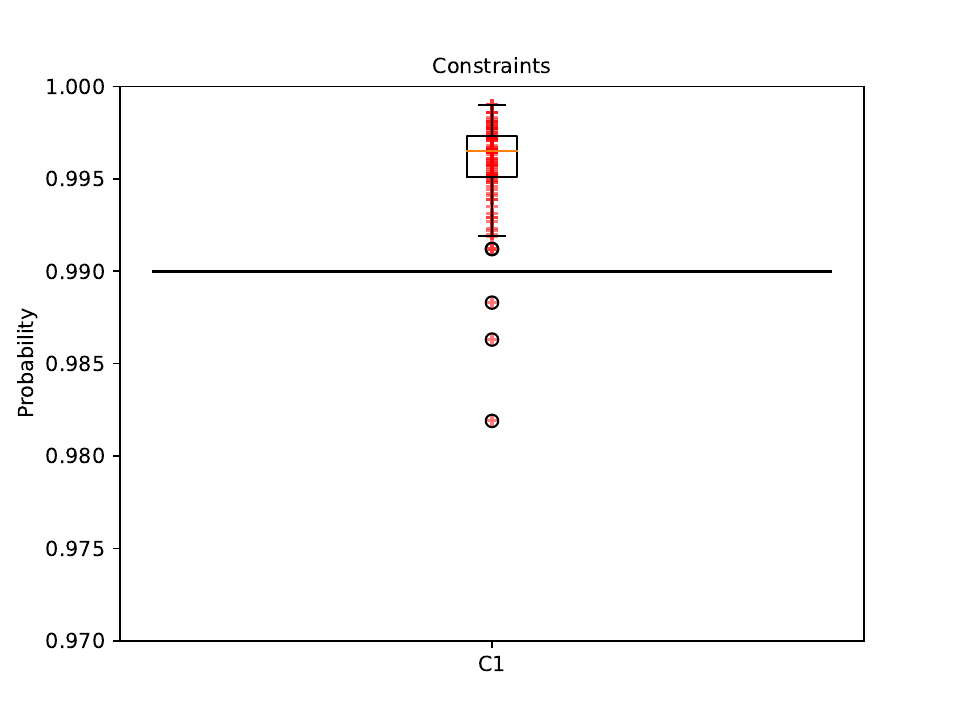}
    \end{minipage}
    }
    \caption{Detail result for Steel Column Design problem with truncated Gaussian gradient approximation}
    \label{fig5}
\end{figure}

\begin{figure}[htb!]
    \hbox{
    \begin{minipage}[c]{.5\linewidth}
        \centering
        \includegraphics[width = \textwidth]{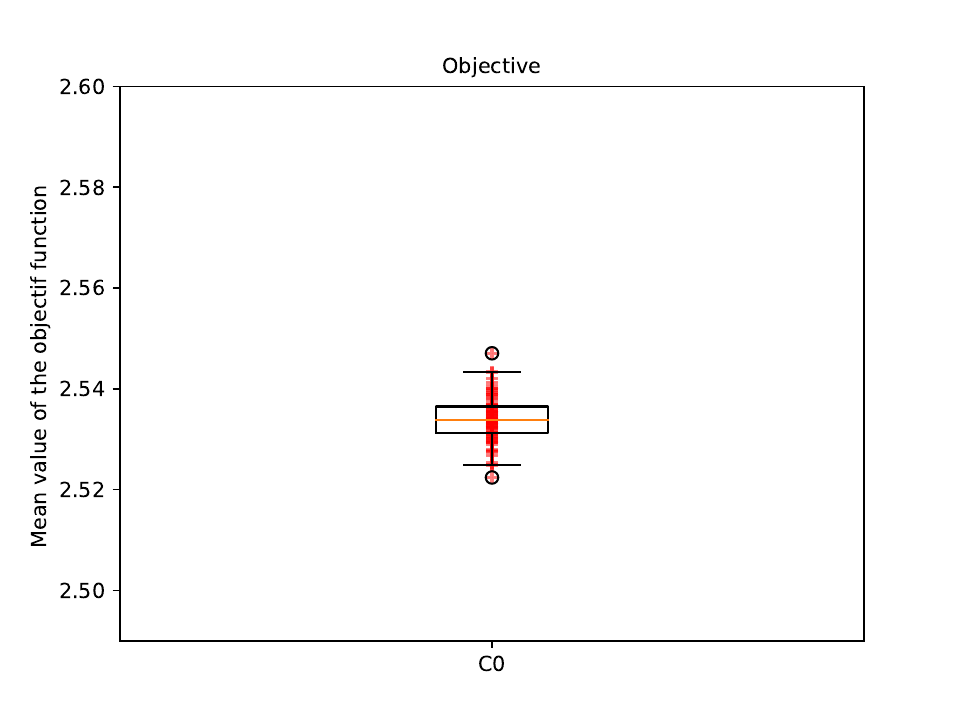}
    \end{minipage}
    \hfill
    \begin{minipage}[c]{.5\linewidth}
        \centering
        \includegraphics[width =\textwidth]{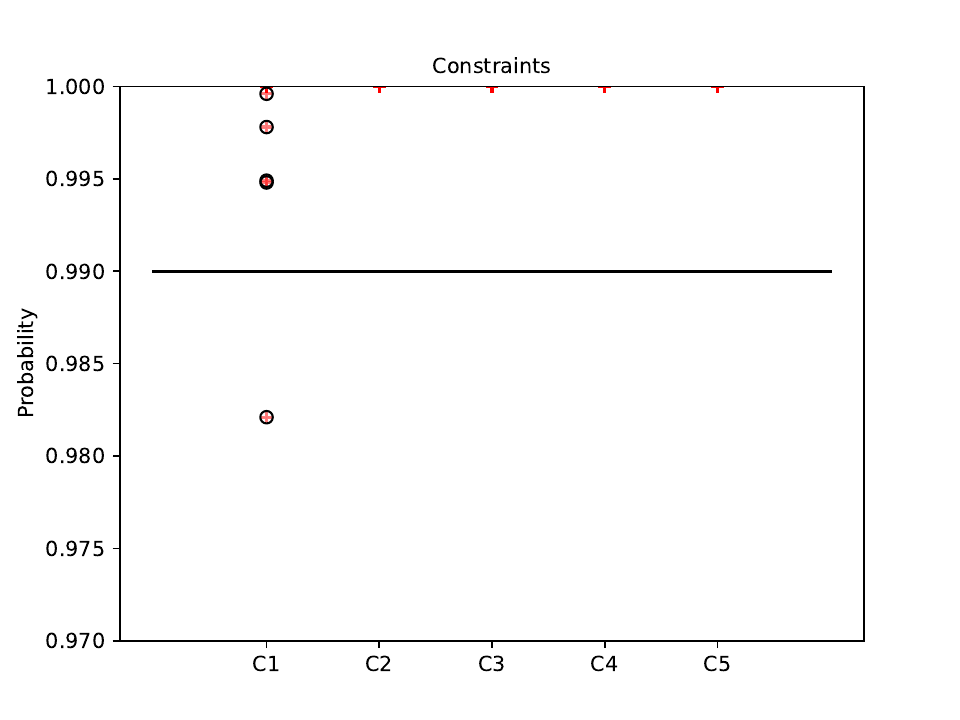}
    \end{minipage}
    }
    \caption{Detail result for Welded Beam Design problem with truncated Gaussian gradient approximation}
    \label{fig6}
\end{figure}

\begin{figure}[htb!]
    \hbox{
    \begin{minipage}[c]{.5\linewidth}
        \centering
        \includegraphics[width = \textwidth]{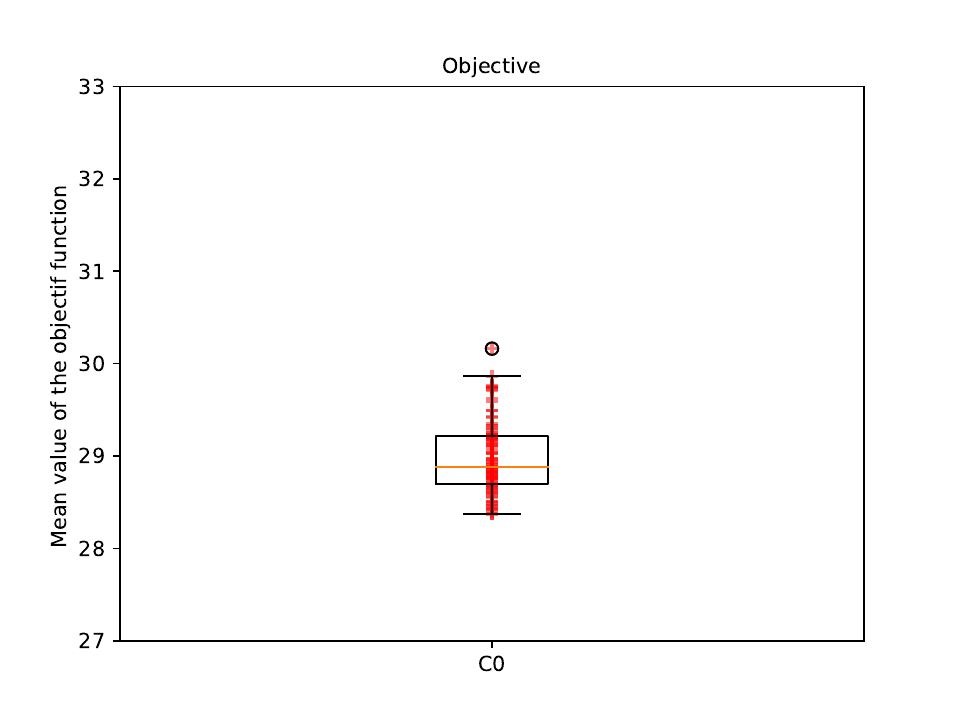}
    \end{minipage}
    \hfill
    \begin{minipage}[c]{.5\linewidth}
        \centering
        \includegraphics[width = \textwidth]{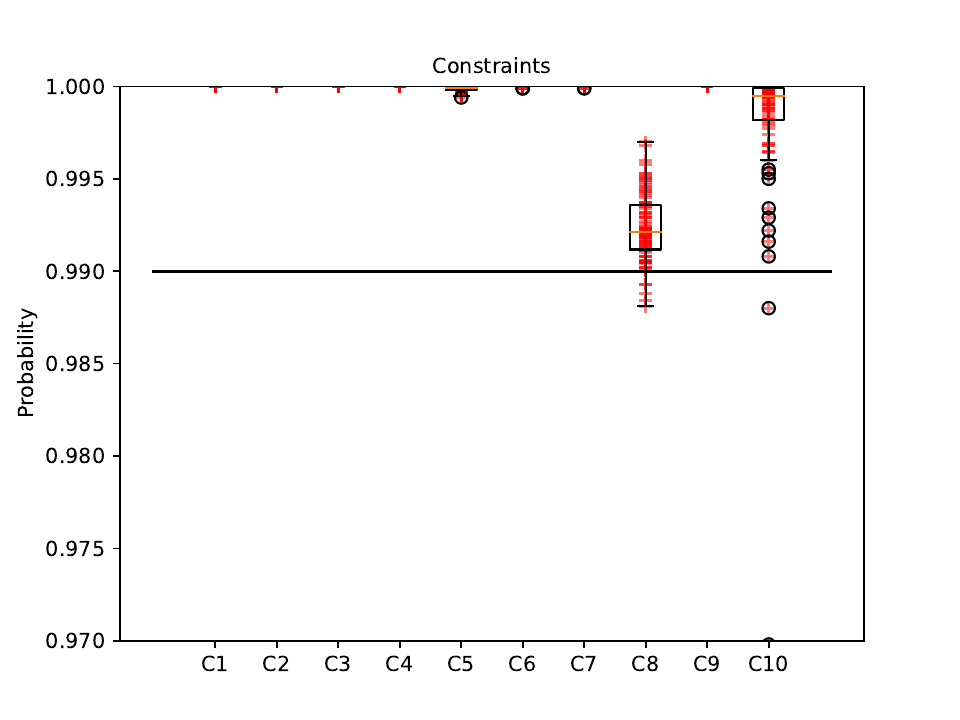}
    \end{minipage}
    }
    \caption{Detail result for Vehicle Side Impact problem with truncated Gaussian gradient approximation}
    \label{fig7}
\end{figure}

\begin{figure}[htb!]
    \hbox{
    \begin{minipage}[c]{.5\linewidth}
        \centering
        \includegraphics[width = \textwidth]{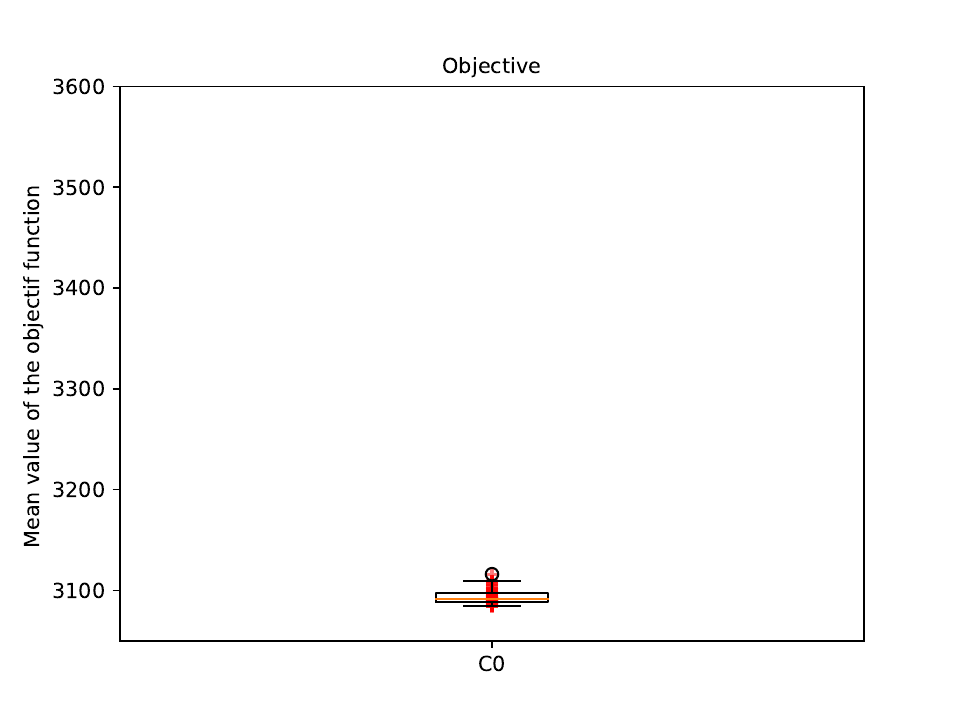}
    \end{minipage}
    \hfill
    \begin{minipage}[c]{.5\linewidth}
        \centering
        \includegraphics[width = \textwidth]{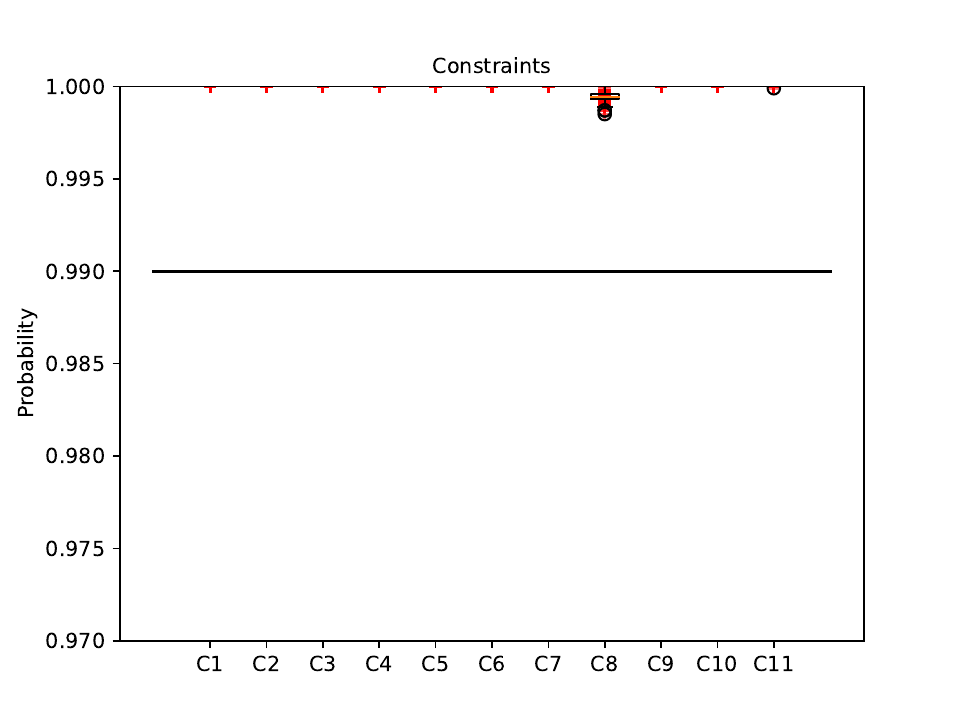}
    \end{minipage}
    }
    \caption{Detail result for Speed Reducer Design problem with truncated Gaussian gradient approximation}
    \label{fig8}
\end{figure}

\end{document}